\documentclass[a4paper,12pt]{elsarticle}
\usepackage[english]{babel}
\addto\captionsenglish{}
\usepackage{mathtools}
\usepackage{dsfont}
\usepackage{color}
\usepackage{amsfonts}
\usepackage{mathrsfs}
\usepackage[table]{xcolor}
\usepackage{amssymb}
\usepackage{geometry}
\usepackage{multirow}
\usepackage{acronym}
\usepackage{xcolor}
\definecolor{red}{RGB}{255,0,0}
\definecolor{blue}{rgb}{0.0, 0.4, 0.65}
\definecolor{orange}{RGB}{255,84,0}
\definecolor{green}{RGB}{154,205,50}
\usepackage{textcomp}
\usepackage{graphics,graphicx,caption,subfig}
\graphicspath{{figures/}}

\newtheorem{remark}{Remark}

\newcommand{\FIN}{\textsc{fin}}
\newcommand{\LL}{\textsc{l}}
\newcommand{\LA}{\textsc{la}}

\renewcommand{\AA}{\textsc{a}}
\newcommand{\SYM}{\mathrm{sym}}
\newcommand{\EXT}{\mathrm{ext}}
\newcommand{\MAX}{\mathrm{max}}
\newcommand{\MIN}{\mathrm{min}}
\newcommand{\BE}{\textsc{be}}
\newcommand{\NS}{\textsc{ns}}
\newcommand{\CH}{\textsc{ch}}
\newcommand{\DIFF}{\mathrm{diff}}
\newcommand{\UPP}{\mathrm{upp}}
\newcommand{\LOW}{\mathrm{low}}
\renewcommand{\eqref}[1]{(\ref{#1})}
\newcommand{\EQ}[1]{\eqref{eq:#1}}
\newcommand{\SEC}[1]{\ref{sec:#1}}
\newcommand{\FIG}[1]{\ref{fig:#1}}
\newcommand{\APP}[1]{\ref{app:#1}}
\newcommand{\TAB}[1]{\ref{tab:#1}}

\journal{Computer Methods in Applied Mechanics and Engineering}

\begin{document}

\begin{frontmatter}

\title{A Robust and Accurate Adaptive Approximation Method for a Diffuse-Interface Model of Binary-Fluid Flows}

\author[TUE1]{T.H.B. Demont}
\author[Evalf]{G.J.~van~Zwieten}
\author[TUE1,UT]{C. Diddens}
\author[TUE1]{E.H.~van~Brummelen}

\address[TUE1]{%
  Eindhoven University of Technology, Department of Mechanical Engineering,
  P.O. Box 513, 5600 MB Eindhoven, The Netherlands}
\address[UT]{%
  University of Twente, Faculty of Science \& Technology,
   P.O. Box 217, 7500 AE Enschede, Netherlands}
\address[Evalf]{%
  Evalf Computing N.V., Burgwal 45, 2611 GG Delft, The Netherlands}

\begin{abstract}
We present an adaptive simulation framework for binary-fluid flows, based on the Abels--Garcke--Gr\"un Navier--Stokes--Cahn--Hilliard (AGG NSCH) diffuse-interface model. The adaptive-refinement procedure is guided by a two-level hierarchical a-posteriori error estimate, and it effectively resolves the spatial multiscale behavior of the diffuse-interface model. To improve the robustness of the solution procedure and avoid severe time-step restrictions for small-interface thicknesses, we introduce an $\varepsilon$\nobreakdash-continuation procedure, in which the diffuse interface thickness ($\varepsilon$) are enlarged on coarse meshes, and the mobility is scaled accordingly. To further accelerate the computations and improve robustness, we apply a modified Backward Euler scheme in the initial stages of the adaptive-refinement procedure in each time step, and a Crank--Nicolson scheme in the final stages of the refinement procedure. To enhance the robustness of the nonlinear solution procedure, we introduce a partitioned solution procedure for the linear tangent problems in Newton's method, based on a decomposition of the NSCH system into its NS and CH subsystems. We conduct a systematic investigation of the conditioning of the monolithic NSCH tangent matrix and of its NS and CH subsystems for a representative 2D model problem. To illustrate the properties of the presented adaptive simulation framework, we present numerical results for a 2D oscillating water droplet suspended in air, and we validate the obtained results versus those of a corresponding sharp-interface model.
\end{abstract}

\begin{keyword}
Navier--Stokes--Cahn--Hilliard equations\sep
diffuse-interface models\sep
binary-fluid flows\sep
adaptive refinement\sep
$\varepsilon$\nobreakdash-continuation\sep
partitioned solution methods
\end{keyword}

\end{frontmatter}

\section{Introduction}

Binary-fluid flows in which the two fluid components are separated by a molecular transition layer are omnipresent in science and engineering. 
An important high-tech application of (typically, aqueous-aereous) binary fluids pertains to inkjet printing, where functionalized minuscule liquid droplets are deposited on a substrate. Inkjet printing is in fact the most widespread technological application of microfluidics~\cite{Lohse:2022vh}. Inkjet-printing technology exhibits exponential growth in a similar manner as Moore's law in semi-conductor technology~\cite{Lohse:2022vh}, and to sustain this technological development, it is imperative to enable the deposition of increasingly many, increasingly small droplets per unit time. Modeling and simulation of binary-fluid flows in inkjet printing are essential to support the necessary technological advancements. Due to the small spatial and temporal length scales of the jetting process, in combination with the complexity of the involved microfluidic processes, such as wetting on fluid-solid interfaces, interface coalescence and break-up, phase transition, surfactant transport and dynamic surface-energy phenomena, etc., advanced numerical simulation techniques form an indispensable research instrument, complementary to experimental investigations.

Mathematical-physical models for binary-fluid flows can generally be subsumed under two categories, viz.\ sharp-interface and diffuse-interface models. In sharp-interface models, the transition layer between the two fluid components is modeled as a manifold of co-dimension one (e.g.\ a 2D surface in a 3D ambient space), at which the interaction of the two fluid components is accounted for by kinematic and dynamic (possibly extended by thermal and thermodynamic) interface conditions. In diffuse-interface models, the transition layer is represented as a thin-but-finite layer, in which the two components are mixed in a proportion that varies continuously and monotonously between the two pure species across the transition layer. Diffuse-interface models are significantly more versatile than their sharp-interface counterpart, in the sense that diffuse-interface models can intrinsically provide a description of phenomena that cannot be intrinsically accounted for in sharp-interface models, such as topological changes of the fluid-fluid interface due to coalescence or break-up~\cite{Lowengrub:1998uq}, and wetting of solid substrates~\cite{Yue:2011uq,Seppecher:1996kx,Jacqmin:2000kx}. Diffuse-interface models therefore provide a modeling paradigm par excellence for complex microfluidic processes, e.g.\ those occurring in inkjet printing. In situations where both models are valid, the sharp-interface model can generally be obtained in the so-called sharp-interface limit, i.e.\ in the limit as the diffuse-interface thickness passes to zero~\cite{Abels:2018ly,Lowengrub:1998uq,Yue:2010hq}. It is to be noted, however, that contemporary theory on sharp-interface limits is still incomplete, and that for instance scaling of model parameters is presently unresolved. By virtue of their intrinsic representation of wetting phenomena, diffuse-interface models have recently also been propounded as a compelling modeling option for binary fluids in elasto-capillary fluid-solid interaction~\cite{Brummelen:2021aw,Shokrpour-Roudbari:2019au,Brummelen:2017sh,Aland:2021vn,Aland:2020nd,Bueno:2018vs}.

Diffuse-interface models for two immiscible incompressible fluid species are generally described by the Navier--Stokes-Cahn--Hilliard (NSCH) equations. The NSCH system is presently a class of models rather than a single model, as various variants exist. Historically, one of the first instances of the NSCH model is the so-called ``model-H'', presented by Hohenberg and Halperin~\cite{Hohenberg:1977hh} in the late 1970s. This model implicitly assumes that the densities of the two fluid species match, which leads to a solenoidal (incompressible) mixture-velocity field. Since its introduction, the Hohenberg--Halperin model has also often been used in non-matching density scenarios. However, in this case the model can violate thermodynamic consistency. An NSCH model that allows for non-matching densities was proposed in the late 1990s by Lowengrub and Truskinovsky~\cite{Lowengrub:1998uq}. This model is based on a barycentric mixture velocity and the mass fraction features as order parameter (phase variable). The resulting NSCH model is quasi-incompressible, i.e.\ the mixture velocity in this model is generally non-solenoidal. A new quasi-incompressible NSCH system in which the volume fraction acts as order parameter, was introduced by Shokrpour et al. in~\cite{Simsek:2018gb}. This model was developed in conjunction with a linearly-implicit, thermodynamically consistent time integration scheme. A deficiency of the quasi-incompressible NSCH models in~\cite{Lowengrub:1998uq,Simsek:2018gb} is that these models are not generally compatible with the underlying Navier--Stokes equations in pure-species regimes: both models insist that the pressure is harmonic if the phase variable is uniform.
This limitation can be avoided by introducing a degenerate mobility, but this comes at the expense of severe complications in analytical and numerical treatments of the model. Abels, Garcke and Gr\"un presented a thermodynamically consistent NSCH model for binary-fluids with non-matching densities, based on volume-averaged mixture velocity~\cite{Abels:2012vn}. By virtue of the volume-averaged mixture-velocity premise, the mixture velocity in the AGG NSCH model is solenoidal. Mass conservation and thermodynamic consistency are inbuilt in the AGG NSCH model via an auxiliary relative mass flux in the balance laws for mass and linear momentum.

Despite the significant progress that has been made in the formulation and analysis of NSCH models for binary-fluid flows, and in the development of numerical techniques for these models, many outstanding challenges remain. One of the main challenges in numerical approximation of the NSCH system, pertains to maintaining robustness and efficiency of the simulations for realistic settings of the model parameters, e.g.\ for binary fluids comprising water and air, due~to:
\begin{enumerate}
\item {\em spatial multiscale behavior\/}:
In most applications, the thickness of the fluid-fluid interface is many orders of magnitude smaller than other relevant length scales in the problem under consideration, e.g.\ the radius of a droplet or of a liquid jet. For instance, x-ray measurements~\cite{Braslau:1988bu} and molecular-dynamics simulations~\cite{Townsend:1991bs} have conveyed that the thickness of a liquid-vapor interface is typically a few \r{A}, while the radius of a pico-liter droplet, which are the smallest droplets used in current technologies, is approximately $10\,$\textmu m. This corresponds to a scale difference of more than $10^4$, which is beyond the resolution capabilities of standard (non-adaptive) approximation methods. For many applications it is unnecessary to adhere to the actual value of the interface thickness, and it suffices that the interface thickness is sufficiently small relative to other characteristic length scales. A ratio of 1:100 is typically indicated as a guideline~\cite{Yue:2010hq}. Nonetheless, the resulting resolution requirements are generally prohibitive if the problem under consideration exhibits additional multiscale behavior, e.g.\ due to the occurrence of satellite droplets~\cite{Notz:2004px}. 
\item {\em temporal multiscale behavior\/}: 
Typically, the diffusion process at the diffuse interface carries a time scale that is significantly smaller than other characteristic time scales in the problem under consideration, e.g.\ the period of oscillation of a droplet. In fact, it is conjectured that in certain cases the diffusive time scale needs to pass to zero in the limit as the diffuse-interface thickness passes to zero~\cite{Abels:2018ly,Yue:2010hq}. Efficiency of a time-integration scheme for NSCH systems therefore requires that such a scheme is uniformly stable with respect to the interface-diffusion time scale. Recently, there has been much progress in the area of unconditionally thermodynamically stable time integration schemes for NSCH systems; see, e.g.,~\cite{Simsek:2018gb,El-Haddad:2022mj} and references therein. These methods however tend to be very dissipative, and their accuracy is inadequate to efficiently solve problems with weakly dissipative dynamics, such as oscillations of low-viscosity droplets. 
\item {\em $\varepsilon$\nobreakdash-conditional stability\/}: 
A further complication related to time integration schemes for NSCH systems, is that the admissible time step is generally restricted by a CFL type condition in relation to the diffuse-interface thickness, i.e.\ the admissible time step is limited by $\lesssim{}\varepsilon/v$, where $\varepsilon$ indicates the diffuse-interface thickness parameter, and $v$ the (normal) velocity of the interface (see Section~\SEC{eps_cont} below). This restriction may emerge from conditional stability of the time integration scheme, or from potential non-robustness of the iterative solution procedure for implicit time-integration methods. This implies that to retain robustness, increasingly small time steps are required as the diffuse-interface thickness decreases, while the (macroscopic) time scale related to the interface motion is invariant. Such a CFL-type condition leads to prohibitive computational expenses in the sharp-interface limit $\varepsilon\to+0$.
\item {\em heterogenous parameters\/}:
In practice, the components of the binary fluid can exhibit significantly different properties, leading to large parameter variations in the NSCH system. For a binary fluid composed of water and air, for instance, the density and viscosity ratios are approximately 1000:1 and 100:1, respectively. This can cause ill-conditioning in the algebraic systems emerging from discrete approximations of the NSCH system. Moreover, in combination with non-compliance of the phase variable, large density and viscosity contrasts can lead to inadmissible negative densities and viscosities of the mixture. Non-compliance refers to the fact that the phase variable can assume values outside of the admissible codomain.  
The phase variable generally represents a difference of fractions (or a fraction), and its codomain is restricted to 
the interval $[-1,1]$ (resp.~$[0,1]$). The codomain of the phase variable must however be extended, thus including inadmissible values, to ensure existence of a solution to the NSCH system; see Ref.~\cite{Grun:2016gi}. With standard constitutive relations for the density and viscosity, these inadmissible values of the phase variable can yield negative mixture densities or viscosities.
\item {\em incompatibility and spurious velocities\/}: Equilibrium solutions of the NSCH system are generally characterized by a particular profile of the phase variable across the diffuse interface (typically, $\tanh(\cdot)$) and vanishing velocity. However, as the particular equilibrium interface profile 
is generally incompatible with the discrete approximation space, discretization errors in the phase variable occur, which manifest as discretization errors in the interface stresses (in particular, the Korteweg stress tensor), and induce discretization errors in the velocity field. These discretization errors in the velocity field are commonly referred to as {\em spurious velocities\/} 
or {\em parasitic currents\/}~\cite{Scardovelli:1999mo}. These spurious velocities scale inversely proportional with the viscosity, 
and can become excessively large if the viscosity in one of the fluid components is small.
\item {\em ill-conditioning of linear systems\/}: Implicit time-integration methods for NSCH systems generally require the solution of a sequence of nonlinear algebraic systems of equations. The solution of each of these nonlinear problems by Newton's method, in turn, requires the solution of a sequence of linear algebraic problems, corresponding to the linear tangent problems in Newton's method. The linear-algebraic systems emerging from the NSCH equations are generally hard, owing to the intrinsic asymmetry of these systems, the aggregation of disparate subsystems (NS and CH), the occurrence of small parameters, and the occurrence of heterogeneous coefficients that may differ by orders of magnitude; see item~4 above. 
\end{enumerate}

In this work, we present an adaptive simulation framework for binary fluids based on the AGG NSCH model, which effectively addresses the aforementioned issues. To resolve the spatial multiscale behavior related to thin diffuse interfaces, we consider adaptive approximations based on high-regularity truncated hierarchical B-splines~\cite{Hughes2005, Cottrell:2009ad}. In particular, we utilize $C^0$ continuous cubic spline bases for all variables except the pressure field, which is approximated by $C^0$ quadratic splines. The velocity-pressure pair thus corresponds to a member of the Taylor--Hood family. The adaptive-refinement procedure is guided by a two-level hierarchical a-posteriori error estimate. The error estimate is used in a support-based marking-and-refinement procedure~\cite{Kuru:2014fk}. The adaptive procedure efficiently controls discretization errors and, hence, intrinsically suppresses spurious velocities. To enhance the robustness of the Newton solution process throughout the iterative-refinement procedure, we introduce an $\varepsilon$\nobreakdash-continuation approach in which the diffuse-interface parameter ($\varepsilon$) and mobility ($m$) are adapted to the resolution of the mesh. This continuation procedure moreover serves to bypass the severe time-step restriction that would otherwise occur for small~$\varepsilon$ (see item~3 above). To determine the appropriate scaling of $m\coloneqq{}m_{\varepsilon}$ with respect to~$\varepsilon$, we derive the characteristic diffusive time scale corresponding to the Cahn--Hilliard equations, which governs the rate of contraction (resp.\ expansion) of the diffuse interface caused by decreasing (resp.\ increasing) of~$\varepsilon$. To further accelerate the computations and improve robustness, we apply a Backward Euler scheme with a second-order contractive/expansive splitting of the double-well potential~\cite{Wu:2014tg} in the initial stages of the adaptive-refinement procedure in each time step, and a Crank--Nicolson scheme in the ultimate stages to restore second-order accuracy in time. Instability due to non-compliance of the phase variable in combination with large density and viscosity contrasts, is circumvented by applying a soft-clipped interpolation of the mixture density~\cite{Grun:2016gi,Bonart:2019re} and an Arrhenius interpolation of the mixture viscosity~\cite{Brummelen:2021aw,Arrhenius:1887xr}. To avoid non-robustness of the Newton solution procedure due to ill-conditioning of the linear tangent problems, we introduce a partitioned solution procedure in each time step and each iteration of the adaptive-refinement process, where we make use of the inherent composition of the NSCH system of NS and CH subsystems. In addition, we present a systematic investigation of the conditioning of the monolithic NSCH tangent matrix and of its NS and CH subsystems, based on numerical computations of an oscillating liquid droplet in an ambient fluid in 2D. For this test case, we also examine the convergence behavior of preconditioned GMRES, regarding two different preconditioners, corresponding to two different partitions of the NSCH system into its NS and CH subsystems. To illustrate the properties of the presented adaptive simulation framework, we finally present numerical results for a 2D oscillating droplet with parameter values representative of a binary fluid composed of water and air, in which the interface-thickness parameter is approximately 70$\times$ smaller than the minimal radius of curvature of the interface, and we validate the AGG NSCH model by comparing the obtained results with those of a corresponding sharp-interface model.

This paper is structured as follows: In Section~\SEC{problem_statement} we present the strong and weak forms of the diffuse-interface AGG NSCH model, together with the applied constitutive relations. In Section~\SEC{error_adaptive} we introduce the approximation setting and the spatial and temporal discretization schemes, followed by a presentation of the error-estimation procedure and the adaptive-refinement process. Section~\SEC{partitioned_monolithic} considers the Newton solution procedure for the nonlinear algebraic systems emerging from the discretization, and introduces the $\varepsilon$-continuation process and the partitioned solution method for the linear tangent problems.  Section~\SEC{numerical_experiments} presents numerical experiments. Finally, in Section~\SEC{concl} we conclude with a summary and a discussion.

%===============================================================================================================%
%===============================================================================================================%
%===============================================================================================================%

\section{Problem statement}
\label{sec:problem_statement}
\begin{figure}[h!!!!!!!!!!]
\begin{center}
\includegraphics[width=\textwidth]{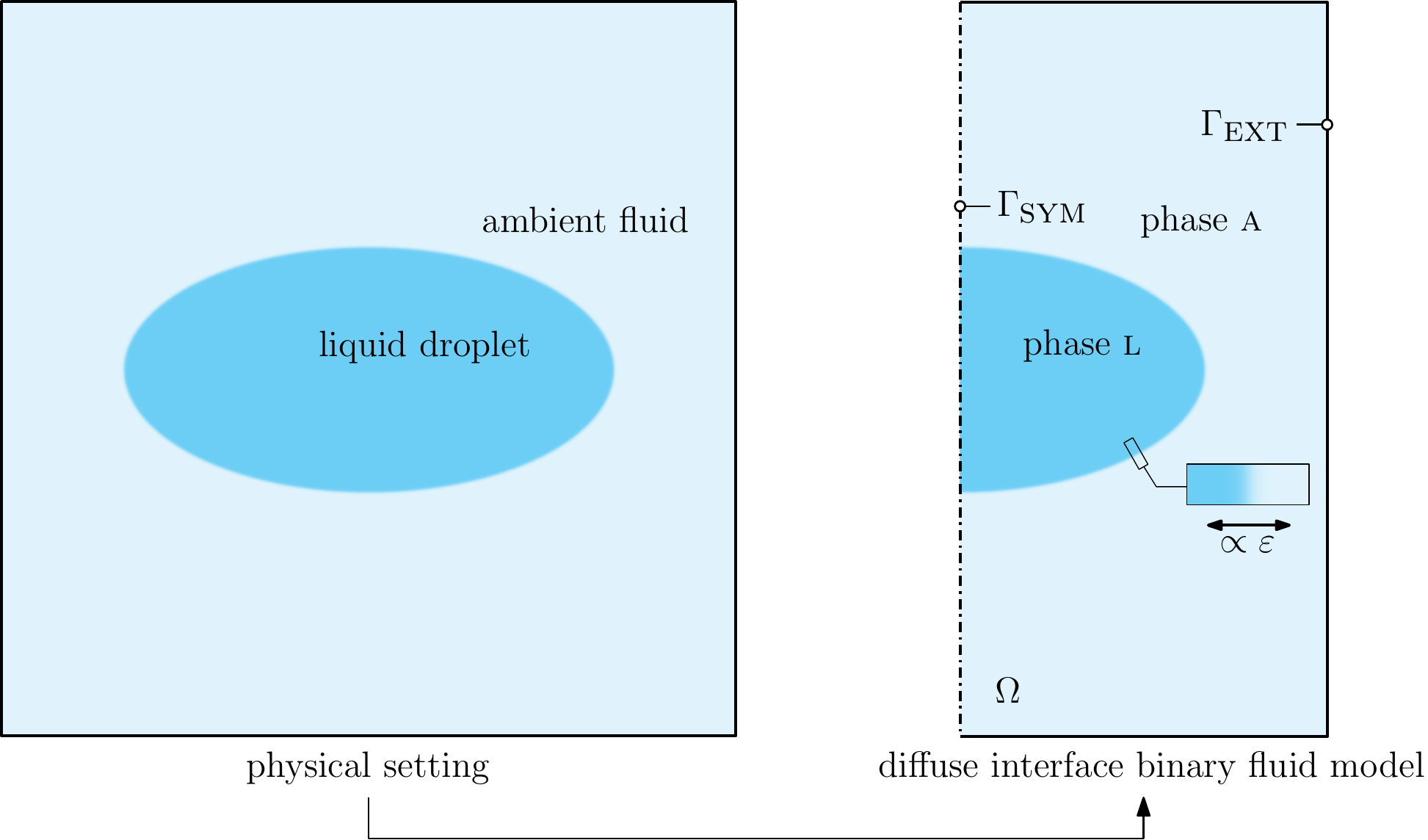}
\end{center}
\caption{On the left, a schematic of the considered physical setting of a liquid droplet submerged in an ambient fluid is displayed. On the right, the binary fluid model is illustrated, separating the fluids by a diffuse interface of thickness $\propto \varepsilon$. The domain is cut in half, making use of one of the symmetry planes through the droplet. In Section~\SEC{numerical_experiments}, simulations are performed on a quarter droplet domain, making use of the second plane of symmetry, not illustrated here.}
\label{fig:droplet_model}
\end{figure}

We consider a binary fluid consisting of two immiscible, incompressible species separated by a finite layer which is comprised of a mixture of those species. The setting for the binary fluid under consideration is a liquid droplet ({\LL}) submerged in an ambient fluid ({\AA}); see the illustration in Figure~\FIG{droplet_model}. To describe the motion of the binary fluid,
we consider the incompressible Navier--Stokes--Cahn--Hilliard (NSCH) model as presented by Abels, Garcke and Gr\"un in~\cite{Abels:2012vn}. This NSCH model is commonly referred as the AGG model, after its originators.
It is noteworthy that alternate NSCH models exist, which can generally be classified into incompressible models~\cite{Hohenberg:1977hh,Ding:2007om,Shen:2010ys} and quasi-incompressible models~\cite{Simsek:2018gb,Lowengrub:1998uq}. Our choice for the AGG model is motivated by its thermodynamic consistency for non-matching species densities and the convenience of implementation related to the incompressibility condition. It is also notable that the AGG model consistently reduces to the underlying single-species Navier--Stokes equations in the pure species setting. Well-posedness of the AGG model in various settings has recenty been established in~\cite{Giorgini:2021bg,Abels:2013bd,Abels:2013pa}. Let us also note the application of the AGG NSCH model in elasto-capillarity in~\cite{Brummelen:2021aw}.

\subsection{NSCH model}
\label{sec:NSCHmodel}
To define a setting for the NSCH model, we regard an open time interval $(0, t_{\FIN}) \subseteq \mathds{R}_{>0}$ and a spatial domain corresponding to a simply connected time-independent subset $\Omega \subset \mathds{R}^d$ ($d=2,3$).
The two fluid species {\LL} and {\AA} are both present according to the volume fractions $\varphi_{\LL}, \varphi_{\AA} \colon (0,t_{\FIN}) \times \Omega \rightarrow [0,1]$ such that $\varphi_{\LL} + \varphi_{\AA} = 1$. The order parameter $\varphi$, identifying the presence of each species, is defined as
\begin{equation}
\label{eq:prob_form_phi}
\varphi \coloneqq \varphi_{\LL} - \varphi_{\AA} \colon (0,t_{\FIN}) \times \Omega \rightarrow [ -1, 1 ],
\end{equation}
from which the identities $\varphi_{\LL} = (1 + \varphi)/2$ and $\varphi_{\AA} = (1 - \varphi)/2$ follow. 
Thus, $\varphi = 1$ and $\varphi = -1$ correspond to pure {\LL}-species and {\AA}-species, respectively, while $\varphi \in (-1,1)$ corresponds to a mixture of the species.

The mixture density $\rho \colon (0,t_{\FIN}) \times \Omega \rightarrow \mathds{R}_{>0}$, where $\rho_{\LL}, \rho_{\AA} > 0$ denote the specific densities of the two species, is related to the volume fractions by $\rho = \rho_{\LL} \varphi_{\LL} + 
\rho_{\AA} \varphi_{\AA}$, and thus
\begin{equation}
\label{eq:rhophi}
\rho \coloneqq \rho \left( \varphi \right) = \frac{ 1 + \varphi }{ 2 } \rho_{\LL} + \frac{ 1 - \varphi }{ 2 } \rho_{\AA}.
\end{equation}
In the AGG NSCH model, the mixture velocity $\boldsymbol u \colon (0,t_{\FIN}) \times \Omega \rightarrow \mathds{R}^d$ is volume-averaged, which leads to a solenoidal mixture velocity. The equation of motion contains an auxiliary mass flux. The AGG NSCH 
diffuse-interface model \cite{Abels:2012vn} is given by:
\begin{equation}
\label{eq:strong}
\left.
\begin{split}
\partial_t \left( \rho \boldsymbol u \right) + \nabla \cdot \left( \rho \boldsymbol u \otimes \boldsymbol u \right) + \nabla \cdot \left( \boldsymbol u \otimes \boldsymbol J \right) + \nabla p - \nabla \cdot \boldsymbol \tau - \nabla \cdot \boldsymbol \zeta & = 0\\
\nabla \cdot \boldsymbol u & = 0\\
\partial_t \varphi + \nabla \cdot \left( \varphi \boldsymbol u \right) - \nabla \cdot \left( m \nabla \mu \right) & = 0\\
\mu + \sigma \varepsilon \Delta \varphi - \frac{ \sigma } { \varepsilon } \Psi' & = 0
\end{split}
\right\}.
\end{equation}
Here,
\begin{equation}
\label{eq:J}
\boldsymbol J \coloneqq m \frac{ \rho_{\AA} - \rho_{\LL} } { 2 } \nabla \mu
\end{equation}
is the auxiliary relative mass flux~\cite{Abels:2012vn}, $\Psi \left( \varphi \right) = \frac{1}{4} \left( \varphi^2 -1 \right)^2$ is the standard polynomial double-well potential for the energy density of the mixture, and the variables $p \colon (0, t_{\FIN}) \times \Omega \rightarrow \mathds{R}$ and $\mu \colon (0, t_{\FIN}) \times \Omega \rightarrow \mathds{R}$ represent the fluid pressure and chemical potential, respectively. We assume the mobility coefficient, $m>0$, to be constant. We note that a degenerate mobility as function of the phase indicator $\varphi$, i.e.\ $m\coloneqq m(\varphi)$ with $m:[-1,1]\mapsto\mathds{R}_{\geq{}0}$ such that $m$ vanishes at and only at $\varphi=\pm{}1$, is preferable from a modeling perspective~\cite{Barrett:1999nx}, as it eliminates the Ostwald-ripening effect in pure species. However, the introduction of a degenerate mobility complicates numerical-approximation procedures~\cite{Barrett:1999nx}, and we therefore opt to retain a non-degenerate constant mobility.

The fluid stress consists of a viscous component, $\boldsymbol \tau$, a capillary component, $\boldsymbol \zeta$, and the isotropic component, $p \boldsymbol I$. The viscous stress tensor is defined by the constitutive relation
\begin{equation}
\label{eq:tau}
\boldsymbol \tau \coloneqq \eta \left( 2 \boldsymbol{\varepsilon} + \Lambda \textrm{tr} \left( \boldsymbol{\varepsilon} \right) \boldsymbol I \right),
\end{equation}
where $\boldsymbol{\varepsilon} \coloneqq \frac{1}{2}( \nabla \boldsymbol u + ( \nabla \boldsymbol u )^T )$ is the symmetric strain-rate tensor and $\textrm{tr} \left( \boldsymbol{\varepsilon} \right) = \nabla \cdot \boldsymbol u$ is its trace. We note that for the AGG NSCH model, the second term in~\EQ{tau} can be omitted on account of~\EQ{strong}$_2$. We have however retained this term (with $\Lambda=-2/d$ in accordance with Stokes' hypothesis; see~\cite{Brummelen:2017sh}) and the term does not vanish for the weakly-solenoidal Taylor--Hood finite-element approximations considered in this work. The mixture viscosity $\eta \coloneqq \eta \left( \varphi \right)$ is defined according to the Arrhenius relation
\[
\log \eta \left( \varphi \right) = \frac{ \left( 1 + \varphi \right) \Lambda_\eta \log \eta_{\LL} + \left( 1 - \varphi \right) \log \eta_{\AA} } { \left( 1 + \varphi \right) \Lambda_\eta + \left( 1 - \varphi \right) }.
\]
$\Lambda_\eta = \frac{ \rho_{\LL} m_{\AA} } { \rho_{\AA} m_{\LL} }$ is the intrinsic liquid-ambient volume ratio, where $m_{\LL}$ and $m_{\AA}$ are the molar densities corresponding to species {\LL} and {\AA}, respectively; see also~\cite{Brummelen:2021aw}. Furthermore, the constitutive relation for the capillary stress tensor is given by
\begin{equation}
\label{eq:zeta}
\boldsymbol \zeta = - \sigma \varepsilon \nabla \varphi \otimes \nabla \varphi + \boldsymbol I \left( \frac{ \sigma \varepsilon } { 2 } | \nabla \varphi |^2 + \frac{ \sigma } { \varepsilon } \Psi \right).
\end{equation}
The parameter $\sigma$ is a rescaling of the liquid-ambient surface tension $\sigma_{\LA}$ according to $2 \sqrt{2} \sigma = 3 \sigma_{\LA}$. The model parameter $\varepsilon > 0$ dictates the transverse length scale (or thickness) of the diffuse interface between the two species. It is noteworthy that in~\EQ{strong} the separate Navier--Stokes and Cahn--Hilliard equations are coupled through the relative mass flux, capillary stress tensor, and advection term $\nabla \cdot \left( \varphi \boldsymbol u \right)$.

Equation~\EQ{prob_form_phi} insists that the codomain of  $\varphi$ is restricted to $[-1,1]$, to enable an unambiguous physical interpretation. We must however desist from this restriction on the codomain of the order parameter to ensure existence of a solution to the NSCH system; see Ref.~\cite{Grun:2016gi}. This implies that~$\varphi$ can assume values outside the interval $[-1,1]$, and accordingly the mixture density in~\EQ{rhophi} can assume negative values. To avoid the occurrence of negative densities, we extend the definition of the mixture density to all of~$\mathds{R}$ via a so-called soft-clipped linear interpolation. As the assignment of the labels $\LL$ and $\AA$ to the fluid species is in principle arbitrary, without loss of generality we assume $\rho_{\LL} > \rho_{\AA}$. Defining the regularization parameter $\delta = \rho_{\AA} / \left( \rho_{\LL} - \rho_{\AA} \right)$, we consider the following extended density function:
\begin{equation}\label{eq:density}
\rho \left( \varphi \right) = \left\{
\def\arraystretch{1.25}
\begin{tabular}{ll}
$\frac{ 1 } { 4 } \rho_{\AA},$ & $\varphi \in ( - \infty , - 1 - 2 \delta ],$\\
$\frac{ 1 } { 4 } \rho_{\AA} + \frac{ 1 } { 4 } \rho_A \delta^{-2} \left( 1 + 2 \delta + \varphi \right)^2,$ & $\varphi \in ( - 1 - 2 \delta , - 1 - \delta),$\\
$\frac{ 1 + \varphi } { 2 } \rho_{\LL} + \frac{ 1 - \varphi } { 2 } \rho_{\AA},$ & $\varphi \in [ - 1 - \delta, 1 + \delta ],$\\
$\rho_{\LL} + \frac{ 3 } { 4 } \rho_{\AA} - \frac{ 1 } { 4 } \rho_{\AA} \delta^{-2} \left( 1 + 2 \delta - \varphi \right)^2,$ & $\varphi \in ( 1 + \delta, 1 + 2 \delta ),$\\
$\rho_{\LL} + \frac{ 3 } { 4 } \rho_{\AA},$ & $\varphi \in [ 1 + 2 \delta, \infty ).$
\end{tabular}
\right.;
\end{equation}
see also~\cite{Grun:2016gi,Bonart:2019re}. The mixture-density definition~\EQ{density} coincides with~\EQ{rhophi} in the interval $[-1,1]$, but is smoothly extended in such a manner that negative densities are prevented. 

A general treatment of auxiliary conditions for the NSCH system~\EQ{strong} is beyond the scope of this work, and we restrict ourselves here to a consideration of the initial and boundary conditions as applied in the numerical experiments in Section~\SEC{numerical_experiments}. The considered test cases pertain to a liquid droplet suspended in an ambient fluid in $\mathds{R}^2$; see the illustration in Figure~\FIG{droplet_model}. We consider configurations that display bilateral symmetry along two axes and, accordingly, only part of the domain is modeled and the symmetry is accounted for by corresponding boundary conditions. The boundary $\Gamma\coloneqq\partial\Omega$ of the domain $\Omega$ is partitioned into two complementary disjoint open subsets, corresponding to the symmetry boundary $\Gamma_{\SYM} \subset \Gamma$ and the external boundary $\Gamma_{\EXT}\coloneqq \operatorname{int}( \Gamma \setminus \Gamma_{\SYM})$. On the symmetry boundaries, we impose a free-slip condition for the velocity field and homogeneous Neumann conditions for the order parameter and the chemical potential:
\begin{equation}
\label{eq:SYM}
\left.
\begin{split}
\partial_n \mu & = 0\\
\partial_n \varphi & = 0\\
u_n & = 0\\
p\boldsymbol n - \boldsymbol \tau \boldsymbol n - \boldsymbol \zeta \boldsymbol n & = 0
\end{split}
\right\} \quad\textrm{on } \Gamma_{\SYM}.
\end{equation}
Let us note that the homogeneous Neumann condition on the order parameter $\varphi$ implies that intersections of the diffuse interface with the symmetry boundary occur at a 90$^\circ$ contact angle. On the lateral boundary $\Gamma_{\EXT}$, we 
impose a no-slip condition, along with the same homogeneous Neumann boundary conditions for $\mu$ and $\varphi$:
\begin{equation}
\label{eq:EXT}
\left.
\begin{split}
\partial_n \mu & = 0\\
\partial_n \varphi & = 0\\
\boldsymbol u & = 0
\end{split}
\right\} \quad\textrm{on } \Gamma_{\EXT}.
\end{equation}
We assume that the domain is large enough to render the effect of the boundary conditions imposed on $\Gamma_{\EXT}$ on the droplet dynamics negligible compared to those induced by the viscous and capillary forces. For the combination of boundary conditions in~\EQ{SYM}--\EQ{EXT}, the pressure variable $p$ is only determined up to a constant. We therefore in addition impose the condition that $p$ vanishes on average, i.e.\ $\langle{}p\rangle\coloneqq\int_\Omega p\, \textrm{d} \boldsymbol x = 0$.

NSCH equations such as~\EQ{strong} intrinsically constitute {\em mixture\/} models, describing the evolution of a binary fluid composed of a mixture of two immiscible constituents. However, NSCH models are usually applied as diffuse-interface models for binary fluids, by regarding only phase-segregated configurations. Conforming to this perspective, we consider initial conditions that correspond to a liquid volume set in an ambient background fluid, separated by a phase-transition 
with a prescribed profile. Specifically, denoting the liquid volume by $\Omega_{\LL}\subset\Omega$ and the corresponding liquid--ambient interface by $\Gamma_{\LA} \coloneqq \partial \Omega_{\LL} \setminus \partial \Omega$, the initial condition for the order parameter is given by:
\begin{equation}
\label{eq:varphi0}
\varphi (0,\boldsymbol x ) = \varphi_0(\boldsymbol x ) \coloneqq \tanh \left( \frac{ d_{\pm} ( \boldsymbol x , \Gamma_{\LA} ) } { \sqrt{2} \varepsilon } \right)
\qquad\text{in }\Omega,
\end{equation}
where $d_{\pm}(\boldsymbol x ,\Gamma_{\LA})$ represents the signed (Euclidian) distance from~$\boldsymbol x$ to~$\Gamma_{\LA}$, positive for $\boldsymbol{x}\in\Omega_{\LL}$. The function 
$s\mapsto\tanh (s/\sqrt{2} \varepsilon )$ represents an equilibrium solution for the phase field of the Cahn--Hilliard equations in one spatial dimension (see, e.g.\ \cite{Cahn:1958fk}, \cite[\S8.4.3]{Brummelen:2017sh}) and, accordingly, the phase field~\EQ{varphi0} is meta-stable if the radius of curvature of $\partial\Omega_{\LL}$ is sufficiently large 
relative to~$\varepsilon$. The initial conditions for~\EQ{strong} are completed by a specification of the mixture velocity. We restrict ourselves in this work to homogeneous initial conditions,
\begin{equation}
\label{eq:u0}
\boldsymbol{u}(0,\boldsymbol x )=\boldsymbol{u}_0(\boldsymbol{x})\coloneqq 0
\qquad\text{in }\Omega,
\end{equation}
implying that the binary-fluid is initially stationary and that its dynamics are induced by the unbalance of the 
capillary and pressure forces corresponding to the initial configuration of the liquid volume encoded in~\EQ{varphi0}.
 
\subsection{Weak formulation of the NSCH system}
\label{sec:Weak}
In this section we present a weak formulation of the NSCH problem \EQ{strong}. This weak formulation also forms the basis for the considered adaptive finite-element approximation. We will not attempt to specify the function spaces for the weak formulation precisely, and instead proceed formally. We define the test space for the weak formulation of~\EQ{strong} as:
\[
\mathds{U} = \left\{ \boldsymbol u \colon \Omega \rightarrow \mathds{R}^d, p \colon \Omega \rightarrow \mathds{R}, \varphi \colon \Omega \rightarrow \mathds{R}, \mu \colon \Omega \rightarrow \mathds{R} \right\}.
\]
The auxiliary condition $\langle{}p\rangle=0$ will be accounted for in the weak formulation by means of the Lagrange-multiplier formalism. Accordingly, the test space is extended by a constant function, which we identify with an element of $\mathds{R}$.  The trial space is a linear space of $\mathds{U}$-valued functions on the considered time interval $(0,t_{\FIN})$, 
denoted by $\mathscr{L}( \left( 0 , t_{\FIN} \right) ; \mathds{U})$. The trial space for the Lagrange multiplier is denoted by
$\mathscr{L}( ( 0 , t_{\FIN} ) ; \mathds{R})$. The aggregated NSCH problem can then be condensed into the following weak formulation:
\begin{multline}
\label{eq:prob_weak_form}
\textrm{\emph{Find }} U \in \big\{( \boldsymbol u , p , \varphi , \mu ) \in 
\mathscr{L} ( ( 0 , t_{\FIN} ) ; \mathds{U} ):\boldsymbol{u}(0,\cdot)=\boldsymbol{u}_0,\varphi(0,\cdot)=\varphi_0\big\},
\chi \in \mathscr{L} ( ( 0 , t_{\FIN} ) ; \mathds{R} )
\textrm{\emph{ such that:}}
\\
\frac{ \textrm{\emph{d}} } { \textrm{\emph{d}} t } a_0 \left( U \left( t \right) ; V \right) + a_1 \left( U \left( t \right) ; V \right) + a_2 \left( U \left( t \right) ; V \right) 
+ \alpha\langle{p}(t)\rangle
+ \chi(t)\langle{q}\rangle 
= l ( V )
\\ \forall \{ V \in \big\{ \left( \boldsymbol v , q , \omega , \lambda \right) \in \mathds{U} : v_n=0 \textrm{ \emph{on} } \Gamma_{\SYM}, \boldsymbol v = 0 \textrm{ \emph{on} } \Gamma_{\EXT} \big\}, \alpha\in\mathds{R}\}, \textrm{ a.e. } t \in ( 0, t_{\FIN} ).
\end{multline}
The semi-linear forms $a_0, a_1, a_2 \colon \mathds{U} \times \mathds{U} \rightarrow \mathds{R}$ are given by
\begin{subequations}
\begin{align}
a_0 \left( U ; V \right) &=  \int_\Omega \left( \rho \boldsymbol u \cdot \boldsymbol {v} + \varphi \omega \right) \textrm{d}\boldsymbol x,\\
a_1 \left( U ; V \right) &=  \int_\Omega \left( \vphantom{\frac{\sigma}{\varepsilon}} \nabla \cdot \left( \rho \boldsymbol u \otimes \boldsymbol u \right) \cdot \boldsymbol {v} -( \boldsymbol u \otimes \boldsymbol J ) : \nabla\boldsymbol {v} + \left( \boldsymbol \tau + \boldsymbol \zeta \right) : \nabla \boldsymbol {v} \right.\notag\\
& \left. +  m \nabla \mu\cdot \nabla \omega + (\boldsymbol u\cdot\nabla \varphi)\omega    + \mu \lambda - \sigma \varepsilon \nabla \varphi \cdot \nabla \lambda - \frac{ \sigma } { \varepsilon } \Psi' \left( \varphi \right) \lambda \right) \textrm{d} \boldsymbol x, 
\label{eq:a_1}\\
a_2 \left( U ; V \right) &=  \int_\Omega \left( q \nabla \cdot \boldsymbol u - p \nabla \cdot \boldsymbol {v} \right) \textrm{d} \boldsymbol x,
\label{eq:a_2}
\end{align}
\end{subequations}
and the linear form $l \colon \mathds{U} \rightarrow \mathds{R}$ contains the exogenous data. To disambiguate the tensor notation, we note that $\nabla\cdot(\boldsymbol{u}\otimes\boldsymbol{v})\cdot\boldsymbol{w}=w_i\partial_j(u_iv_j)$ and $( \boldsymbol u \otimes \boldsymbol v ) : \nabla\boldsymbol {w}=u_iv_j\partial_jw_i$ in Cartesian components (with summation on repeated indices). The semi-linear form~$a_0$ contains the two time-derivative terms. The first three terms of~$a_1$, together with the second term in~$a_2$, correspond to the equation of motion in~(\ref{eq:strong})$_1$. The fourth and fifth terms of $a_1$ represent the convective Cahn--Hilliard equation~(\ref{eq:strong})$_3$. Let us note that $\nabla\cdot(\varphi\boldsymbol{u})$ has been replaced by~$\boldsymbol{u}\cdot\nabla\varphi$, which is a consistent modification for solenoidal~$\boldsymbol{u}$. The remaining terms in~$a_1$ define the chemical potential in~(\ref{eq:strong})$_4$. The first term in~$a_2$ accounts for the solenoidality of the volume-weighted mixture velocity~(\ref{eq:strong})$_2$. 

%===============================================================================================================%
%===============================================================================================================%
%===============================================================================================================%

\section{Error estimation and adaptive approximation}
\label{sec:error_adaptive}
In most applications of the NSCH equations it is necessary to select the transverse length scale of the diffuse 
interface, $\varepsilon$, significantly smaller than other characteristic length scales in the problem under consideration, e.g.\ the radius of a liquid droplet or the size of the domain. To account for the resulting spatial multi-scale behavior of the NSCH equations, we consider adaptive approximations of~\EQ{prob_weak_form}. Specifically, we consider adaptive approximations based on Truncated Hierarchical B-splines~ \cite{Giannelli:2012rr} (THB-splines) analogous to~\cite{Brummelen:2021aw}. In this work, we restrict ourselves however to conventional $C^0$ continuous approximation, so that the obtained results in Section~\SEC{numerical_experiments} extend immediately to conventional finite-element approximations. The adaptive-refinement procedure follows the standard recursive SEMR (\texttt{Solve} $\rightarrow$ \texttt{Estimate} $\rightarrow$ \texttt{Mark} $\rightarrow$ \texttt{Refine}) process \cite{Bertoluzza:2012kx, Dorfler:1996uq}. To enhance robustness and efficiency of the adaptive-approximation procedure, we apply a first order Backward Euler time-integration scheme on coarse meshes and a second-order Crank--Nicolson scheme on the finest meshes.

\subsection{Approximation of the NSCH system}
\label{sec:err_adap_approx}
We regard a discretization of the NSCH system~\EQ{prob_weak_form} with respect to the spatial dependence based on THB splines. The details of the THB-splines are omitted here, and the reader is referred to~\cite[\S3.1]{Brummelen:2021aw} for an overview. We define the hierarchical approximation space as:
\[
\mathds{U}_L = \left\{ \boldsymbol v \in T^\# \left( \textrm{span} \left[ \mathds{H}_L^{3,0} \right]^d \right), p \in T^\# \left( \textrm{span} \mathds{H}_L^{2,0} \right), \varphi \in T^\# \left( \textrm{span} \mathds{H}_L^{3,0} \right), \mu \in T^\# \left( \textrm{span} \mathds{H}_L^{3,0} \right) \right\},
\]
where $\mathds{H}^{k, \alpha}_L$ is the piece-wise polynomial THB spline space of order $k$ and regularity $\alpha$. The subscript~$L=0,1,\ldots,L_{\MAX}$ denotes the refinement level, and the spaces $\mathds{U}_L$ are hierarchical in the sense that
$\mathds{U}_L\subset\mathds{U}_{L+1}$. The operator $T^\#$ is the pullback of the transformation $T \colon \Omega \rightarrow (0,1)^d$ from the physical domain to the unit square. Note that with the indicated choice of $k$ and $\alpha$, we have $C^0$ continuous spline bases of degree~3 for all but the (adjusted) pressure field~$p$, the latter's degree being 1 lower. Accordingly, the velocity-pressure pair is a member of the Taylor--Hood family.

To enhance the stability of the discrete approximation, the convective term in the Navier--Stokes equations in~\EQ{a_1} is replaced by a skew-symmetric formulation according to Ref.~\cite{Layton:2008fk}, by adding the following semi-linear form to the formulation:
\begin{equation}
\label{eq:s_skew}
\begin{aligned}
s \left( U ; V \right) & = 
- \frac{1}{2} \int_\Omega \nabla \cdot \left( \rho \boldsymbol u \otimes \boldsymbol u \right) \cdot \boldsymbol {v}\,\textrm{d} \boldsymbol x 
- \frac{1}{2} \int_\Omega \nabla \cdot \left( \rho \boldsymbol {v} \otimes \boldsymbol u \right) \cdot \boldsymbol u\,\textrm{d} \boldsymbol x\\
& \phantom{=} 
+ \frac{1}{2} \int_\Omega \boldsymbol u \cdot \boldsymbol {v}\, \boldsymbol u \cdot \nabla \rho\,\textrm{d} \boldsymbol x 
+ \frac{1}{2} \int_{\partial\Omega} \rho u_n\, \boldsymbol u \cdot \boldsymbol {v}\,\textrm{d} \boldsymbol x.
\end{aligned}
\end{equation}
By means of the product rule and Gauss' theorem, one can establish that $s$ according to~\EQ{s_skew} forms a partition of zero for solenoidal~$\boldsymbol{u}$. Augmenting the weak formulation with~$s$ causes the volumetric terms associated with convective transport, except for the term corresponding to $\nabla\rho$, to cancel if $\boldsymbol{v}$ is replaced by~$\boldsymbol{u}$, irrespective of solenoidality of~$\boldsymbol{u}$:
\begin{equation}
\int_{\Omega}\nabla \cdot \left( \rho \boldsymbol u \otimes \boldsymbol u \right) \cdot \boldsymbol {u}\,\textrm{d} \boldsymbol x
+s(U,U)=\frac{1}{2}\int_{\Omega}|\boldsymbol{u}|^2\boldsymbol u \cdot \nabla \rho\,\textrm{d} \boldsymbol x + \frac{1}{2} \int_{\partial\Omega} \rho u_n |\boldsymbol u |^2\,\textrm{d} \boldsymbol x.
\end{equation}
In pure-species regions (i.e.\ regions where $\nabla\rho=0$), this modification eliminates the artificial energy production that could otherwise occur on account of inexact solenoidality of the Taylor--Hood approximation of the velocity.

To discretize the weak form in time, we introduce a partition of the time interval~$ \left( 0 , t_{\FIN} \right)$ into subintervals of length $\tau_{n} = 2^{-l}\tau_{\MAX}$ for some $l\coloneqq l_n \in \mathds{Z}_{\geq 0}$ and fixed $\tau_{\MAX}$. The time-step-reduction factor~$l$ is determined a-posteriori, based on the convergence of the Newton procedure; see Section~\SEC{eps_cont} below. To enhance the stability of the formulation, on coarse meshes we apply a first order Backward Euler (BE) scheme to discretize the equations in the time dependence, with a second order contractive-expansive splitting of the double-well potential with stabilization~\cite{Wu:2014tg}:%\cite{Eyre:1998qe,Elliott:1993sf}:
\begin{multline}
\label{eq:error_system_BE}
\textrm{\emph{Find }} \{U_L^{n}\coloneqq(\boldsymbol{u},p,\varphi,\mu)^{n}_L \in \mathds{U}_L^{n}, \chi^{n}\in\mathds{R} \}\textrm{\emph{ such that:}}\\
\frac{ 1 } { \tau_{n} } a_0 \left( U_L^{n} ; V_L^{n} \right) + a_1 \left( U_L^{n} ; V_L^{n} \right) + a_2 \left( U_L^{n} ; V_L^{n} \right) +  s \left( U_L^{n} ; V_L^{n} \right)
+s_+ \left( \varphi^n_L,\Pi\underline{\varphi}^{n-1}_{L_{\MAX}},\lambda^n_L \right)
\\
+ \beta \tau_n \tilde{s} \left( \varphi^n_L, \lambda^n_L \right)
+\alpha^n\langle{p^n}\rangle
+\chi^n\langle{q^n}\rangle
= \beta \tau_n \tilde{s} \left( \underline{\varphi}^{n-1}_{L_{\MAX}}, \lambda^n_L \right)
+ l \left( V_L^{n} \right) + \frac{ 1 } { \tau_{n} } a_0 \left( \Pi \underline{U}_{L_{\MAX}}^{n-1} ; V_L^{n} \right) %+s_+(\varphi^{n-1}_{L_{\MAX}},\lambda^n_L)
\\
\forall\{V_L^{n}\coloneqq(\boldsymbol{v},q,\omega,\lambda)^{n}_L \in \mathds{U}_L^{n},\alpha^n\in\mathds{R}\}.
\end{multline}
Here, $s_+$ and $\tilde{s}$ are given by
\[
\begin{split}
s_+ \left(\varphi, \tilde{\varphi}, \lambda \right) & = \int_\Omega \left( \frac{\sigma}{\varepsilon} \Psi'(\varphi) \lambda - \frac{\sigma}{2 \varepsilon} \left( - \tilde{\varphi} - \varphi - \varphi^3 + 3 \tilde{\varphi} \varphi^2 \right) \lambda \right) \textrm{d} \boldsymbol x,\\%, \, \textrm{and}\\
\tilde{s} \left( \varphi, \lambda \right) & = - \int_\Omega \left( \nabla \lambda \cdot \nabla \varphi \right) \textrm{d} \boldsymbol x.
\end{split}
\]
$s_+$ replaces the derivative of the double-well potential by a second order contractive-expansive splitting, while the artificial-diffusion term $\tilde{s}$ acts as stabilization, where the coefficient $\beta$ needs to satisfy conditions as stated in~\cite[Thm.1]{Wu:2014tg} for unconditional energy stability. 
In this work, we select $\beta \coloneqq 2m$, which allows some leniency if $| \varphi | > 1$%, chosen s.t.\ $\beta \geq 2 m$.%; see also~\cite{}
. We use $\Pi$ to denote a generic projection operator. In particular, in~\EQ{error_system_BE}, the projection %$\Pi \underscore{\varphi}^{n-1}_{L_{\MAX}}$ 
$\Pi \underline{U}_{L_{\MAX}}^{n-1}$ (incl. $\Pi \underline{\varphi}^{n-1}_{L_{\MAX}}$) corresponds to an $L^2$\nobreakdash-projection of the super-refined approximation (see Section~\SEC{error_adapt}) at the maximum level of refinement in the previous time step to the level\nobreakdash-$L$ approximation space in time step~$n$; see Remark~\ref{rem:projection} below for further elaboration.
On fine meshes in the sequence of adaptive refinements, we apply a second order Crank--Nicolson (CN) scheme to discretize the equations in the temporal dependence:
\begin{multline}
\label{eq:error_system_CN}
\textrm{\emph{Find }} \{U_L^{n}\coloneqq(\boldsymbol{u},p,\varphi,\mu)^{n}_L \in \mathds{U}_L^{n}, \chi^{n}\in\mathds{R} \}\textrm{\emph{ such that:}}
\\
\frac{ 1 } { \tau_{n} } a_0 \left( U_L^{n} ; V_L^{n} \right) 
+ \frac{1}{2} a_1 \left( U_L^{n} ; V_L^{n} \right) 
+ a_2 \left( U_L^{n} ; V_L^{n} \right) 
+ \frac{1}{2} s \left( U_L^{n} ; V_L^{n} \right)
+\alpha^n\langle{p^n}\rangle
+\chi^n\langle{q^n}\rangle
\\
= 
l \left( V_L^{n} \right) 
+\frac{ 1 } { \tau_{n} } a_0 \left( \underline{U}_{L_{\MAX}}^{n-1} ; V_L^{n} \right) 
-\frac{1}{2} a_1 \left( \underline{U}_{L_{\MAX}}^{n-1} ; V_L^{n} \right)
-\frac{1}{2} s \left( \underline{U}_{L_{\MAX}}^{n-1} ; V_L^{n} \right)  
\\
\forall\{V_L^{n}\coloneqq(\boldsymbol{v},q,\omega,\lambda)^{n}_L \in \mathds{U}_L^{n},\alpha^n\in\mathds{R}\}.
\end{multline}
It is to be noted that in~\EQ{error_system_CN}, the bilinear form~$a_2$ according to~\EQ{a_2} is treated fully implicitly, as are the terms related to the constraint on the average pressure, $\alpha\langle{p}\rangle$, and the corresponding constraint on the test space, $\chi\langle{q}\rangle$. Such a fully implicit treatment of constraints and corresponding Lagrange multipliers is standard in the CN scheme; see, e.g., \cite{John:2006qi}. Note also that in~\EQ{error_system_BE} and~\EQ{error_system_CN}, terms involving the approximation in the previous time step, i.e.\ with index~${n-1}$, are evaluated at the super-refinement of the maximum adaptive-refinement level.

\begin{remark}
\label{rem:projection}
Both~\EQ{error_system_BE} and~\EQ{error_system_CN} contain terms which depend on the approximation obtained in the previous time step at the super-refinement of the maximum level of refinement, $\underline{U}^{n-1}_{L_{\MAX}}$. The assembly of such terms carries a relatively high computational cost, even if the operative approximation space, $\mathds{U}^n_L$, is still coarse, i.e.\ if~$L<L_{\MAX}$.  In the right member of the CN scheme~\EQ{error_system_CN}, $\underline{U}^{n-1}_{L_{\MAX}}$ appears in the bilinear form~$a_0$, and in the semi-linear forms $a_1$ and~$s$. In a standard BE scheme (i.e.\ according to~\EQ{error_system_BE} without projections), $\underline{U}^{n-1}_{L_{\MAX}}$ only appears
in the bilinear form~$a_0$ and in the stabilization terms $s_+$ and~$\tilde{s}$. Because the semi-linear forms $a_1$ and, to a lesser degree, $s$ are significantly more complex than~$a_0$,~$\tilde{s}$ and~$s_+$, on coarse meshes the computational expense of the CN scheme is considerably higher than that of a BE scheme. From the perspective that for $L<L_{\MAX}$ each approximation $U^n_L$ merely serves to guide the adaptive algorithm to construct a suitably refined approximation space at level $L+1$ (and to provide an initial estimate to the Newton iteration procedure to solve for $U^n_{L+1}$; see Section~\SEC{eps_cont}), we opt to use the BE scheme~\EQ{error_system_BE} at the coarsest refinement levels $L=0,1,2,\ldots,L_{\BE}$. On the subsequent refinement levels $L_{\BE}+1,\ldots,L_{\MAX}$, we switch to the CN approximation~\EQ{error_system_CN} to recover second-order temporal accuracy. Application of the BE scheme (with an expansive-contractive decomposition of the double-well potential with stabilization) on coarse meshes, carries the additional advantage that this scheme is generally more stable, which enhances the robustness and efficiency of the Newton solution procedure on the coarse meshes; see Section~\SEC{eps_cont}.

To further reduce the computational complexity of the BE scheme at the coarse levels, instead of using~$\underline{U}^{n-1}_{L_{\MAX}}$ directly, we use its $L^2$\nobreakdash-projection $\Pi\underline{U}^{n-1}_{L_{\MAX}}$ onto~$\mathds{U}^n_L$. The projection needs to be evaluated only once at each level of refinement, and it generally provides a sufficiently accurate approximation for the~$a_0$ and~$s_+$ terms in~\EQ{error_system_BE} for the purpose of guiding the adaptive-refinement procedure. We retain $\underline{\varphi}^{n-1}_{L_{\MAX}}$ in the right-hand side of~\EQ{error_system_BE} because this term involves the gradient of the phase field instead of the phase field in algebraic form, and it needs to be evaluated only once. Especially in the nonlinear term~$s_+$ in~\EQ{error_system_BE}, the use of the projection~$\Pi\underline{\varphi}^{n-1}_{L_{\MAX}}$ instead of~$\underline{\varphi}^{n-1}_{L_{\MAX}}$ is beneficial, because the terms that emanate from~$s_+$ in the nonlinear solution procedure (see Section~\SEC{eps_cont}) require frequent re-evaluation.
\end{remark}

\begin{remark}
In accordance with~\cite{John:2006qi}, in practice we observe that the CN scheme provides significantly better accuracy than the BE scheme. The inaccuracy of the BE scheme is actually prohibitive, in the sense that the time step which is required by the BE scheme to adequately resolve the dynamics of the interface motion, is generally many orders of magnitude smaller than the characteristic time scale of the motion itself. The aforementioned switch between BE and CN within the adaptive refinement procedure however relies on the assumption that within a {\em single\/} time step, the BE approximation is sufficiently accurate to guide the refinement process on the coarse levels, i.e.\ the diffuse-interface features corresponding to the CN approximation on level $L_{\BE}+1$ should be adequately resolved by the super-refinement $\underline{\mathds{U}}^n_{L_{\BE}}$. Even if the overall accuracy of the BE scheme is too limited, it can generally serve this purpose.  
\end{remark}

The systems~\EQ{error_system_BE} and~\EQ{error_system_CN} are considered for consecutive time steps $n=1,2,\ldots$. For $n=1$, the approximation $\underline{U}^0_{L_{\MAX}}$ that appears in the right members of~\EQ{error_system_BE} and~\EQ{error_system_CN}, is replaced by the initial conditions. With reference to Section~\SEC{NSCHmodel}, we note that the initial-boundary-value problem for the NSCH equations is generally furnished with initial conditions for the order parameter $\varphi$ and for the mixture velocity~$\boldsymbol{u}$. Because the right-hand side of the BE scheme~\EQ{error_system_BE} only depends on $a_0$ and $s_+$ and, hence, only on~$\varphi$ and~$\boldsymbol{u}$, the prescribed initial conditions suffice to initialize the BE scheme. The CN scheme~\EQ{error_system_CN} depends on $a_0$, $a_1$ and $s$ and, therefore, on~$\varphi$, $\boldsymbol{u}$ and~$\mu$. To start the CN scheme, the initial chemical potential must thus be specified in addition to the initial order parameter and mixture velocity. The chemical potential can be consistently initialized based on the initial datum for the order parameter, $\varphi_0$, by invoking its definition~(\ref{eq:strong})$_4$. In practice, we accomplish this by an $L^2(\Omega)$ projection onto a highly resolved hierarchical-refined approximation space, $\tilde{V}$, according to:
\begin{equation}
\label{eq:mu0}
\mu_0\in{}\tilde{V}:
\qquad
\int_{\Omega}
\mu_0\lambda\,\textrm{d} \boldsymbol x
=
\int_{\Omega}
\Big(\sigma\varepsilon\nabla\varphi_0\cdot\nabla\lambda+\frac{\sigma}{\varepsilon}\Psi'(\varphi_0)\,\lambda\Big)\,\textrm{d} \boldsymbol x
\qquad
\forall\lambda\in\tilde{V}.
\end{equation}
It is to be noted that the gradient of~$\varphi_0$ that appears in the right member of~\EQ{mu0} cannot generally be computed exactly for initial data of the form~\EQ{varphi0}. However, for the considered test cases, we have access to very accurate series expansions of~$\varphi_0$ and, hence, its gradient.

\subsection{Error estimation and adaptive spatial refinement}
\label{sec:error_adapt}

The diffuse interface introduces features into the solution with a length scale $\varepsilon$ that is generally significantly smaller than other characteristic length scales in the problem under consideration. To efficiently resolve the various features of the solution, including those in the vicinity of the diffuse interface, within each time step in the time-integration process we construct a sequence of adaptively refined meshes, following the standard recursive SEMR process. The adaptive-refinement procedure is directed by a hierarchical a-posteriori error estimate. The error-estimation-and-refinement procedure is described in detail in~\cite{Brummelen:2021aw}, but its main elements are repeated here for completeness.

We employ a so-called {\em two-level hierarchical error estimate\/} wherein, for each field, the approximation in $\mathds{U}_L$ is compared to the approximation in its super-refinement $\underline{\mathds{U}}_L$ to estimate the error. More precisely, for each time step $n$ and refinement level $L$, we define the super-refinement $\underline{\mathds{U}}^n_L$ of the approximation space $\mathds{U}^n_L$ by uniform subdivision of each knot span (element) of $\mathds{U}^n_L$ into $2^d$ equal parts. We denote by~$\underline{U}^n_L\in\underline{\mathds{U}}^n_L$ the super-refined approximation obtained from~\EQ{error_system_BE} or~\EQ{error_system_CN} (depending on the refinement level $L$). Considering appropriate norms for each of the fields $f\in\mathds{F}\coloneqq\{\boldsymbol{u},p,\varphi,\mu\}$, we can then estimate the error in $U^n_L$ as:
\begin{equation}
\label{eq:estimate}
\|\underline{U}_L^n - U_L^n \| \coloneqq \| \underline{\boldsymbol u}_L^n - \boldsymbol u_L^n \|_{H^1 \left( \Omega, \mathds{R}^d \right)} + \| \underline{p}_L^n - p_L^n \|_{L^2 \left( \Omega \right)} + \| \underline{\varphi}_L^n - \varphi_L^n \|_{H^1 \left( \Omega \right)} + \| \underline{\mu}_L^n - \mu_L^n \|_{H^1 \left( \Omega \right)}.
\end{equation}
The applied marking strategy is based on support-wise error indicators~\cite{Kuru:2014fk,Brummelen:2021aw,Brummelen:2017rr,Richter:2015zr}. Specifically, denoting by $\underline{b}^f_L$ any basis function of the super-refined approximation space for field~$f$, we define the error indicators:
\begin{equation}
\begin{aligned}
\imath^{\boldsymbol{u}}(\underline{b}^{\boldsymbol{u}}_L)&=
\| \underline{\boldsymbol u}_L^n - \boldsymbol u_L^n \|_{H^1 \left( \operatorname{supp}^{\circ}(\underline{b}^{\boldsymbol{u}}_L), \mathds{R}^d \right)}
\\
\imath^{p}(\underline{b}^{p}_L)&=
\| \underline{p}_L^n - p_L^n \|_{L^2 \left( \operatorname{supp}^{\circ}(\underline{b}^{p}_L)\right)}
\end{aligned}
\qquad
\begin{aligned}
\imath^{\varphi}(\underline{b}^{\varphi}_L)&=
\| \underline{\varphi}_L^n - \varphi_L^n \|_{H^1 \left( \operatorname{supp}^{\circ}(\underline{b}^{\varphi}_L)\right)}
\\
\imath^{\mu}(\underline{b}^{\mu}_L)&=
\| \underline{\mu}_L^n - \mu_L^n \|_{H^1 \left( \operatorname{supp}^{\circ}(\underline{b}^{\mu}_L)\right)}
\end{aligned}
\end{equation}
where $\operatorname{supp}^{\circ}(\underline{b}^{f}_L)$ denotes the open support of the basis function $\underline{b}^{f}_L$.
By virtue of the overlap in the basis functions, the sum of the support-wise error indicators yields an upper bound to the error estimate~\EQ{estimate}. For each field $f \in \mathds{F}$, we then mark a minimal set of basis functions such that the sum of the corresponding error indicators exceeds a prescribed fraction $a^f\in[0,1)$ of the total sum:
\begin{equation*}
\sum_{\underline{b}^{f}_L\in\operatorname{Marked}}
\imath^{f}(\underline{b}^{f}_L)
\geq{}
a^f
\sum_{\underline{b}^{f}_L\in\operatorname{All}}
\imath^{f}(\underline{b}^{f}_L).
\end{equation*}
Note that we may set $a^f = 0$ for some $f \in \mathds{F}$, in which case the fields $f$ in question are not considered in the adaptive-refinement process. In the refinement step, all basis functions $b^f_L$ of which the support intersects with the support of the marked set will be replaced by their corresponding hierarchical refinements; see~\cite{Brummelen:2021aw} for details. In practice, some marginal auxiliary refinements are introduced in a so-called {\em completion step\/} to retain a regular structure.

%===============================================================================================================%
%===============================================================================================================%
%===============================================================================================================%

\section{Solution methods}
\label{sec:partitioned_monolithic}
The nonlinear systems~\EQ{error_system_BE} and~\EQ{error_system_CN} translate into systems of nonlinear algebraic equations for the coefficients relative to the THB-spline basis functions. Devising a robust solution procedure for the nonlinear algebraic systems emerging throughout the adaptive-refinement procedure and time-integration process, is nontrivial. We consider a solution procedure for the nonlinear problems based on a Newton method with line search. To enhance the robustness of the Newton procedure, we introduce a continuation approach in which the diffuse-interface parameter~$\varepsilon$ and mobility~$m$ are adapted to the resolution of the mesh. This continuation procedure also serves to bypass the severe time-step restriction that could otherwise occur for small~$\varepsilon$. Furthermore, a time-step-reduction procedure is implemented as a fail-safe to the Newton solution method. The latter also aids in retaining robustness of the solution procedure in the presence of the wide range of dynamical features that generally occurs in the evolution of the NSCH equations for small~$\varepsilon$. Motivated by the severe ill-conditioning that can occur in the linear tangent problems in Newton's method for the coupled NSCH system for small $\varepsilon$ (see Section~\SEC{num_exp_cond}), we regard both a monolithic (i.e.\ fully coupled) solution procedure and an iterative solution procedure based on a partitioning of the NSCH system into its composing NS and CH subsystems. 

\subsection{Newton solution procedure with $\varepsilon$-continuation}
\label{sec:eps_cont}

The nonlinear systems~\EQ{error_system_BE} and~\EQ{error_system_CN} give rise to a system of nonlinear algebraic equations of the generic form $\boldsymbol{R}^n_{L}(\boldsymbol{U}^n_L)=0$, where $\boldsymbol{U}^n_L$ contains the coefficients of 
$U^n_L\in\mathds{U}^n_L$ relative to the applied THB-spline bases. Below, we will generally suppress the affixes~$n$ and~$L$, unless these are relevant for the exposition. Let us note that the super-refinements in the a-posteriori error-estimation process are essentially of the same form, and the discussion below applies mutatis mutandis to these super-refinements. The basic Newton procedure with line search can be summarized as follows: given an initial estimate $\boldsymbol{U}^0$, repeat for $k=0,1,2,\ldots$:
\begin{equation}
\label{eq:Newton}
\begin{aligned}
\boldsymbol{A}\,\boldsymbol{\delta}&=-\boldsymbol{R}(\boldsymbol{U}^k)
\\
\boldsymbol{U}^{k+1}&=\boldsymbol{U}^{k}+S\boldsymbol{\delta}
\end{aligned}
\end{equation}
where $\boldsymbol{A}\coloneqq\boldsymbol{A}(\boldsymbol{U}^{k})$ denotes the derivative of the residual $\boldsymbol{R}$.
The step size (or relaxation factor)~$S$ is determined from minimizing a cubic interpolation of the residual norm between $\boldsymbol{U}^k$ and $\boldsymbol{U}^{k+1}$, based on the residual and the tangent matrix at these approximations; see~\cite{nutils} for details on the implementation. The iteration~\EQ{Newton} terminates when the norm of the residual drops below a certain prescribed tolerance.

The robustness and efficiency of the Newton procedure~\EQ{Newton} depend sensitively on the error in the initial estimate, on the nonlinearity of the residual function, and on the properties of the tangent matrix, notably, its condition number in the vicinity of the solution. The Newton procedure is generally non-robust for the nonlinear problems that emerge for small values of the diffuse-interface parameter~$\varepsilon$ on coarse meshes. Specifically, the under-resolved thin diffuse-interface in the first few iterations of the adaptive algorithm can cause divergence of the Newton procedure, or may lead to prohibitively small steps in the line search. 

An additional problem in the solution of NSCH systems, pertains to the fact that the solution from a previous time step fails to represent an accurate approximation to the solution in the time step under consideration if the interface motion within a time step is large in comparison with the transverse length scale, $\varepsilon$. This problem is in fact universal, independent of the adaptive-refinement process. The adaptive-refinement process however enables much smaller values of~$\varepsilon$ than would generally be feasible on uniform meshes and, hence, it prominently exposes the problem. To illustrate this issue, we consider the uniform translation of a 1D diffuse interface described by the usual tangent-hyperbolic profile, with length scale~$\varepsilon$:
\begin{equation}\label{eq:tangent_hyperbolic_profile}
\varphi(t,x)=\tanh\bigg(\frac{x-vt}{\sqrt{2}\varepsilon}\bigg),
\end{equation}
where $v>0$ denotes the translation velocity of the diffuse interface; see the illustration in Figure~\FIG{interface_translation}. Considering a time step~$\tau$, it holds that:

\begin{figure}
\begin{center}
\includegraphics[width=0.8\textwidth]{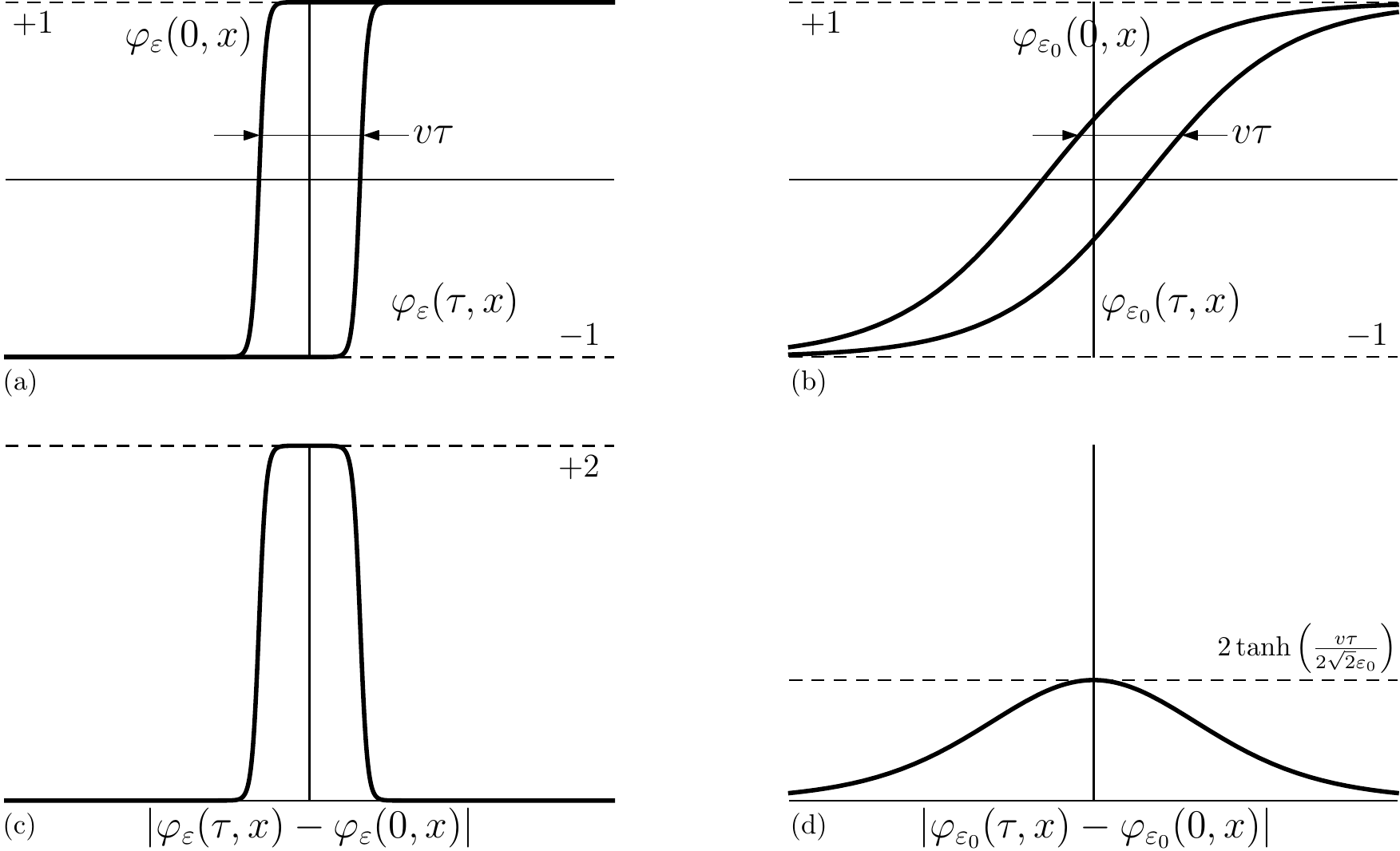}
\end{center}
\caption{A 1D diffuse interface at times $t=0,\tau$ given by Equation~\EQ{tangent_hyperbolic_profile}, uniformly translated over a distance $v \tau$ with length scales $\varepsilon$ (a) and $\varepsilon_0$ (b). The absolute differences of the functions before and after translation are presented in (c) and (d), respectively.}
\label{fig:interface_translation}
\end{figure}

\begin{equation}
\label{eq:O1error}
\max_{x\in\mathds{R}}|\varphi(\tau,x)-\varphi(0,x)|=2\tanh\bigg(\frac{v\tau}{2\sqrt{2}\varepsilon}\bigg)\approx
\frac{v\tau}{\sqrt{2}\varepsilon}
\end{equation}
where the approximation holds to reasonable accuracy in the interval $0\leq{}v\tau/\varepsilon\leq{}1$. Equation~\EQ{O1error} conveys that in order for $\varphi(0,x)$ to provide a suitable approximation to~$\varphi(\tau,x)$, it is necessary that~$v\tau$
is sufficiently small compared to~$\varepsilon$. This implies that if the diffuse-interface thickness~$\varepsilon$ is reduced, the time step in the time-integration procedure has to be decreased proportionally. Otherwise, the initial approximation provided to the Newton procedure is too far off from the actual solution, leading to divergence of the Newton process or excessively small steps in the line-search procedure.

To improve the robustness of the Newton procedure on coarse meshes and to avoid the severe time-step restriction that can emerge for small~$\varepsilon$, we introduce a so-called {\em $\varepsilon$-continuation\/} process, viz.\ we enlarge the parameter~$\varepsilon$ for the first iterations of the adaptive refinement process. As opposed to~\cite{Brummelen:2021aw}, in conjunction with the interface-thickness parameter~$\varepsilon$, we rescale the mobility parameter~$m$, so that the continuation effectively pertains to~$(\varepsilon,m)$. This rescaling of the mobility is imperative to retain robustness if significant interface displacements occur within a time step. To implement the $\varepsilon$\nobreakdash-continuation procedure, we select $K\leq{}L_{\MAX}$ and set $\varepsilon_0\coloneqq 2^K \varepsilon$, such that the diffuse interface with thickness corresponding to~$\varepsilon_0$, is adequately resolved by the baseline approximation space $\mathds{U}^n_0$. Correspondingly, we introduce a mobility parameter $m_0=C_{\tau}\,2^{3K} m$, with $C_{\tau}$ a sufficiently large constant such that the diffuse interface is appropriately equilibrated within a time step; see Remark~\ref{rem:epscont} below. Note that, accordingly, the constant $C_{\tau}$ depends on the time step in the time-integration process. In iteration $L = 0, 1,\ldots, K-1$ of the refinement process, we solve problem~\EQ{error_system_BE} or~\EQ{error_system_CN} with interface-thickness parameter $\varepsilon_L = 2^{-L} \varepsilon_0$ and mobility $m_L = 2^{-3L} m_0$, i.e.\ the interface thickness is decreased during the adaptive refinement and the mobility is scaled correspondingly. For $L = K, \ldots, L_{\MAX}$, we solve \EQ{error_system_BE} or \EQ{error_system_CN} with the original parameters, $\varepsilon_L = \varepsilon$ and $m_L = m$.  
For all refinement levels except the coarsest, the initial approximation for the Newton procedure is obtained from a projection of the super-refinement at the previous refinement level i.e., symbolically, $(\boldsymbol{U}^n_L)^0=\Pi\underline{\boldsymbol{U}}^n_{L-1}$. The projection $\Pi$ is based on a Stokes projector for velocity to ensure solenoidality of the initial approximation, and an $L^2(\Omega)$ projection for the order parameter and the chemical potential.
The initial approximation for the super-refinements at all levels is obtained via the same projection of the approximation itself, i.e.\ $(\underline{\boldsymbol{U}}^n_L)^0=\Pi\boldsymbol{U}^n_{L}$. The initial approximation for the Newton process at the coarsest level is provided by a projection of the coarsest-level solution of the previous time step, i.e.\ $(\boldsymbol{U}^n_0)^0=\Pi\boldsymbol{U}^{n-1}_0$. By virtue of the $\varepsilon$\nobreakdash-continuation process, the interface is adequately resolved on all meshes that occur in the sequence of adaptive refinements. Moreover, the error in the initial approximation for each Newton process, except at the coarsest level, is essentially independent of the motion of the interface,
and it depends only on the continuation error, i.e.\ the error induced by the transition from length scale~$\varepsilon_{L-1}$ 
to~$\varepsilon_L$. At the coarsest level, the initial error does depend on the motion of the interface but, with reference to~\EQ{O1error}, the error is only proportional to $v\tau/\varepsilon_0$. The time step restriction can therefore be 
effectively diminished by selecting the diffuse-interface thickness on the coarsest mesh, $\varepsilon_0$, sufficiently large.

\begin{remark}
In the transition from refinement~$K$ to~$K+1$, the mobility decreases by the factor $C_{\tau}$. Our experience is that the robustness of the solution procedure is insensitive to~$C_{\tau}$. However, alternate, more gradual continuation schemes can be devised if necessary.
\end{remark}

\begin{remark}
\label{rem:epscont}
On all levels of the adaptive-refinement procedure and, accordingly, for all values of~$\varepsilon_L$, the 
phase field at the previous time step that appears in the right-hand side of~\EQ{error_system_BE} and~\EQ{error_system_CN} 
corresponds to the original interface-thickness parameter,~$\varepsilon$. This implies that in the initial stages of the
$\varepsilon$\nobreakdash-continuation process, where $\varepsilon_L>\varepsilon$, the diffuse interface essentially has to expand from length scale~$\varepsilon$ to length scale~$\varepsilon_L$ within a single time step, $\tau$. To accommodate this expansion, the mobility parameter, $m$, has to be chosen sufficiently large. 
The characteristic diffusive time scale of the Cahn--Hilliard equations, $T_{\DIFF}$, is given by
\begin{equation*}
T_{\DIFF} = \frac{\varepsilon^3}{\sigma m};
\end{equation*}
see~\APP{time_scale}. To maintain the same interfacial expansion rate throughout the different stages of the $\varepsilon$-continuation process, we therefore introduce a scaling of the mobility according to~$m_L\propto{}\varepsilon_L^3$. The diffuse interface is generally suitably equilibrated for $\tau\geq{}T_{\DIFF}$; see~\APP{interface_expansion}.
\end{remark}

\subsection{Direct solution method}
\label{sec:monolithic}
Central to the Newton iteration procedure, is the solution of the linear tangent problems in~(\ref{eq:Newton})$_1$. It appears that the properties of tangent problems corresponding to the NSCH equations have not been systematically investigated previously. An important objective of this work is therefore to elucidate the properties of these tangent problems, by considering the condition number of the tangent matrix~$\boldsymbol{A}$ appearing in~(\ref{eq:Newton})$_1$ and of its corresponding submatrices; see Section~\SEC{num_exp_cond}. Foregoing a detailed investigation, one can however anticipate that inversion of the tangent matrix~$\boldsymbol{A}$ is generally hard, based on the characteristic features of the NSCH system, viz., the intrinsic asymmetry of the system, the aggregation of disparate subsystems (NS and CH), the occurrence of small parameters (notably, $\varepsilon$ and~$m$), and the occurrence of heterogeneous coefficients that may differ by orders of magnitude (notably, the  species densities $\rho_{\LL}$ and~$\rho_{\AA}$, and species viscosities~$\eta_{\LL}$ and~$\eta_{\AA}$). 

Considering the complexity of solving~(\ref{eq:Newton})$_1$, a robust and general solution procedure is required. Direct solvers are generally regarded as the most robust solution methods for problems of the form~(\ref{eq:Newton})$_1$. The robustness of direct solvers is however not uniform, and depends on underlying strategies such as reordering, rescaling and pivoting. 
In this work, we apply the PARDISO sparse direct solver~\cite{Schenk:2001qo,Schenk:2004eu,Bollhofer:2020hb}, which is among the most robust direct solvers available; see, e.g.,~\cite{Gould:2007mr}. Specifically, we use either the Intel oneMKL PARDISO 
solver~\cite{intelMKL}, or the PARDISO 6.0 solver~\cite{pardiso-7.2a,pardiso-7.2b,pardiso-7.2c}.

\subsection{Partitioned solution method}
\label{sec:partitioned}
Motivated by the complexity of solving the coupled problem~(\ref{eq:Newton})$_1$ directly, we also consider a partitioned solution strategy, in which the tangent matrix is subdivided into parts corresponding to the NS and CH subsystems. Let us note that the Lagrange multipliers which have been introduced to impose the auxiliary condition $\langle{}p{}\rangle=0$ (see Section~\SEC{Weak}) are included in the NS subsystem. We opt to introduce the partitioning of the NSCH system at the level of the linear tangent problems~(\ref{eq:Newton})$_1$, instead of at the level of the (nonlinear) problem formulation. The partitioning strategy can hence be conceived of as a preconditioning approach for the linear tangent problems.

To specify the partitioned solution strategy, we note that in correspondence with the composition of the NSCH system~\EQ{strong} of the NS and CH subsystems, we can introduce a decomposition of the coefficients, 
$\boldsymbol{U}\coloneqq(\boldsymbol{U}_{\NS},\boldsymbol{U}_{\CH})$ and, similarly, of the residual, $\boldsymbol{R}\coloneqq(\boldsymbol{R}_{\NS},\boldsymbol{R}_{\CH})$. Each linear tangent problem~(\ref{eq:Newton})$_1$ in the Newton procedure can then be decomposed as:
\begin{equation}
\label{eq:partitionedtangent}
\begin{pmatrix}
\boldsymbol{A}_{\NS\text{-}\NS} & \boldsymbol{A}_{\NS\text{-}\CH} \\
\boldsymbol{A}_{\CH\text{-}\NS} & \boldsymbol{A}_{\CH\text{-}\CH}
\end{pmatrix}
\begin{pmatrix}
\boldsymbol{\delta}_{\NS}\\
\boldsymbol{\delta}_{\CH}
\end{pmatrix}
=
-
\begin{pmatrix}
\boldsymbol{R}_{\NS}\\
\boldsymbol{R}_{\CH}
\end{pmatrix}.
\end{equation}
In the setting of~\EQ{partitionedtangent}, the partitioned solution strategy consists of omitting one of the off-diagonal blocks, $\boldsymbol{A}_{\NS\text{-}\CH}$ or~$\boldsymbol{A}_{\CH\text{-}\NS}$, in the actual solution step. The resulting system is then block upper or lower triangular, and inversion only requires inversion of the subsystem matrices, $\boldsymbol{A}_{\NS\text{-}\NS}$ and~$\boldsymbol{A}_{\CH\text{-}\CH}$. Instead of using the partitioned method directly as an iterative solution method, we apply it as a preconditioner in a GMRES method~\cite{Saad:1986ak} for solving~(\ref{eq:Newton})$_1$. In particular, instead of solving~(\ref{eq:Newton})$_1$, we solve $\boldsymbol{P}\boldsymbol{A}\,\boldsymbol{\delta}=-\boldsymbol{P}\boldsymbol{R}$ by means of GMRES, where $\boldsymbol{P}\in\{\boldsymbol{P}_{\UPP},\boldsymbol{P}_{\LOW}\}$ and
\begin{equation}
\label{eq:Pupp}
\boldsymbol{P}^{-1}_{\UPP}
=
\begin{pmatrix}
\boldsymbol{A}_{\NS\text{-}\NS} & \boldsymbol{A}_{\NS\text{-}\CH} \\
\boldsymbol{0} & \boldsymbol{A}_{\CH\text{-}\CH}
\end{pmatrix}
\quad\Rightarrow\quad
\boldsymbol{P}_{\UPP}
=
\begin{pmatrix}
\boldsymbol{A}_{\NS\text{-}\NS}^{-1} & -\boldsymbol{A}_{\NS\text{-}\NS}^{-1}\boldsymbol{A}_{\NS\text{-}\CH}\boldsymbol{A}_{\CH\text{-}\CH}^{-1} \\
\boldsymbol{0} & \boldsymbol{A}_{\CH\text{-}\CH}^{-1}
\end{pmatrix}
\end{equation}
if the lower triangular $\boldsymbol{A}_{\CH\text{-}\NS}$ block is omitted and, correspondingly, the upper-triangular part is retained, and
\begin{equation}
\label{eq:Plow}
\boldsymbol{P}^{-1}_{\LOW}
=
\begin{pmatrix}
\boldsymbol{A}_{\NS\text{-}\NS} & \boldsymbol{0} \\
\boldsymbol{A}_{\CH\text{-}\NS} & \boldsymbol{A}_{\CH\text{-}\CH}
\end{pmatrix}
\quad\Rightarrow\quad
\boldsymbol{P}_{\LOW}
=
\begin{pmatrix}
\boldsymbol{A}_{\NS\text{-}\NS}^{-1} & \boldsymbol{0} \\
-\boldsymbol{A}_{\CH\text{-}\CH}^{-1}\boldsymbol{A}_{\CH\text{-}\NS}\boldsymbol{A}_{\NS\text{-}\NS}^{-1} & \boldsymbol{A}_{\CH\text{-}\CH}^{-1}
\end{pmatrix}
\end{equation}
if the upper triangular $\boldsymbol{A}_{\NS\text{-}\CH}$ block is omitted. In the considered partitioned solution 
procedure, the subsystem solves corresponding to the inverse matrices $\boldsymbol{A}_{\NS\text{-}\NS}^{-1}$ and
$\boldsymbol{A}_{\CH\text{-}\CH}^{-1}$ are again performed by means of the PARDISO sparse direct solver. The principal 
advantage of the partitioned solution method relative to the direct solution method, is that the subsystem matrices
$\boldsymbol{A}_{\NS\text{-}\NS}$ and~$\boldsymbol{A}_{\CH\text{-}\CH}$ are much better conditioned and more amenable to rescaling than the aggregated tangent matrix~$\boldsymbol{A}$; see Section~\SEC{num_exp_cond}.

\begin{remark}
\label{rem:PuppPlow}
The partitioned solution methods and preconditioners associated with~\EQ{Pupp} and~\EQ{Plow} are fundamentally different. Essentially, in the upper-triangular preconditioner~\EQ{Pupp} the effect of the transport velocity $\boldsymbol{u}$ in the CH equations is treated explicitly and the effect of the order parameter $\varphi$ and the chemical potential~$\mu$ (via the auxiliary mass flux~$\boldsymbol{J}$ in~\EQ{J} and capillary stress tensor~$\boldsymbol{\zeta}$ in~\EQ{zeta}) in the NS equations is treated implicitly, while for the lower-triangular preconditioner~\EQ{Plow} these effects are treated vice versa. The terms $\nabla\cdot(\boldsymbol{u}\otimes\boldsymbol{J})$ and $\nabla\cdot\boldsymbol{\zeta}$ in~(\ref{eq:strong})$_1$ contain higher-order derivatives of~$\varphi$ and~$\mu$, while~$\boldsymbol{u}$ appears only in algebraic (i.e.\ non-differentiated) form 
in~(\ref{eq:strong})$_3$ and as a second-order derivative in~(\ref{eq:strong})$_1$ via the divergence of the viscous-stress tensor, $\nabla\cdot\boldsymbol{\tau}$. This suggests that $\boldsymbol{P}^{-1}_{\UPP}$ is a much better approximation of the principal part of the NSCH system than~$\boldsymbol{P}^{-1}_{\LOW}$ and, accordingly, that $\boldsymbol{P}_{\UPP}$ is much more effective as a preconditioner than~$\boldsymbol{P}_{\LOW}$. This conjecture is confirmed by the numerical experiments in Section~\SEC{GMRES_conv}.

Let us also note that the $\boldsymbol{A}_{\NS\text{-}\CH}$ block that is retained in the upper-triangular preconditioner (and discarded in the lower-triangular preconditioner) carries the fluid-fluid surface tension, via the capillary tensor $\boldsymbol{\zeta}$. In the upper-triangular preconditioner, surface tension is therefore treated implicitly in the solution procedure, while in the lower-triangular preconditioner it is treated explicitly. 
\end{remark}

%===============================================================================================================%
%===============================================================================================================%
%===============================================================================================================%
%

\section{Numerical experiments}
\label{sec:numerical_experiments}
To establish the main properties of the AGG NSCH model~\EQ{strong}, the adaptive approximation method as outlined in Section~\SEC{error_adaptive}, and the solution procedures presented in Section~\SEC{partitioned_monolithic}, we conduct numerical experiments for an oscillating liquid droplet suspended in ambient fluid. We restrict ourselves here to a two-dimensional setting. We consider an initially ellipsoidal droplet; cf. Figure~\FIG{droplet_model}. We exploit the symmetry of the configuration along the principal axes of the ellipse and accordingly consider only a quarter of the liquid--ambient domain. Specifically, the considered computational domain is a square with sides of 50\,\textmu{}m, i.e.\ $\Omega \coloneqq (0,50)^2\,\text{\textmu m}^2$, and prescribe symmetry conditions according to~\EQ{SYM} on
\begin{equation}
\Gamma_\SYM \coloneqq \big\{(x_1,x_2)\in\partial\Omega:\{x_1=0\,\text{\textmu{}m}\}\cup\{x_2=0\,\text{\textmu{}m}\}\big\}
\end{equation}
The complementary part of the domain boundary is referred to as an external boundary, and is denoted by~$\Gamma_\EXT \coloneqq \partial \Omega \setminus \Gamma_\SYM$.
On this part of the boundary, we impose boundary conditions according to~\EQ{EXT}.

The initial configuration of the droplet corresponds to an ellipse with foci $( \pm 10\,\text{\textmu m}, 0\,\text{\textmu m})$:
\begin{equation}
\label{eq:quartellipse}
\Omega_{\LL} \coloneqq \bigg\{ (x_1,x_2) \in \Omega :  \Big( \frac{x_1}{20\,\text{\textmu m}} \Big)^2 
+ \Big( \frac{x_2}{10\,\textrm{\textmu m}} \Big)^2 \leq 1 \bigg\}.
\end{equation} 
The ellipsoidal initial configuration of the liquid droplet is represented by an initial condition conforming to~\EQ{varphi0} for the order parameter. For the mixture-velocity field, we impose the homogeneous initial condition~\EQ{u0}, with $\Gamma_{\LA}$ corresponding to $\Omega_{\LL}$ in~\EQ{quartellipse}. Let us note that the initial configuration does not correspond to an eigenmode of the (linearized)
oscillating droplet: for infinitesimal radial perturbations, an ellipsoidal shape formally corresponds to the second mode of oscillation, but the corresponding velocity field does not vanish at any moment. An overview of the model-parameter values for the various considered test cases is presented in Table~\TAB{cases}.

\begin{table}[]
\caption{Overview of physical and numerical parameter values for the considered test cases. Entries marked with the symbol $*$ indicate a range of values, which will be specified in the text.}
\label{tab:cases}
\resizebox{\columnwidth}{!}{
\begin{tabular}{|c|cc|cc|ccc|cccc|}
\cline{2-12}
\multicolumn{1}{c|}{}& \multicolumn{2}{c|}{Liquid}&\multicolumn{2}{c|}{Ambient}&\multicolumn{3}{c|}{Interface}&\multicolumn{4}{c|}{Numerical approximation} \\
\multicolumn{1}{c|}{}& $\rho_{\LL}$ & $\eta_{\LL}$ & $\rho_{\AA}$  & $\eta_{\AA}$  & $\varepsilon$  & $\sigma_{\LA}$       & $m$       & $\tau_{\MAX}$             & $h_0$ & $L_{\MAX}$ & $K$ \\
\hline
 TC\#    & $\frac{\textrm{kg}}{\textrm{m}^d}$ &$\frac{\textrm{kg}\,\textrm{m}^{2-d}}{\textrm{s}}$& $\frac{\textrm{kg}}{\textrm{m}^d}$  & $\frac{\textrm{kg}\,\textrm{m}^{2-d}}{\textrm{s}}$ & $10^2$ \textmu m & $\frac{\textrm{kg}\,\textrm{m}^{3-d}}{\textrm{s}^2}$ & $\frac{\textrm{m}^d\,\textrm{s}}{\textrm{kg}}$ & $10$ \textmu s & \textmu m & --- & --- \\ \hline
1A              & $10^3$   & $10^{-3}$ & $10^{3}$ & $10^{-3}$                & $\ast$              & $7.28 \!\times\! 10^{-2}$ & $2.4390 \!\times\! 10^{-10}$ & $2^{-6}$ & $\ast$ & ---    & ---    \\
1B              & $10^3$   & $10^{-3}$ & $10^{3}$ & $10^{-3}$                & $\ast$              & $7.28 \!\times\! 10^{-2}$ & $2.4390 \!\times\! 10^{-10}$ & $2^{-6}$ & $5$    & $\ast$ & $\ast$ \\
1C              & $10^3$   & $10^{-3}$ & $10^{3}$ & $10^{-3}$                & $2^{-7\phantom{0}}$ & $7.28 \!\times\! 10^{-2}$ & $\ast$                       & $2^{-6}$ & $5$    & 3      & 3      \\
1D              & $10^3$   & $10^{-3}$ & $10^{3}$ & $10^{-3}$                & $2^{-7\phantom{0}}$ & $7.28 \!\times\! 10^{-2}$ & $2.4390 \!\times\! 10^{-10}$ & $\ast$   & $5$    & 3      & 3      \\
1E              & $\ast$   & $10^{-3}$ & $\ast$   & $10^{-3}$                & $2^{-7\phantom{0}}$ & $7.28 \!\times\! 10^{-2}$ & $2.4390 \!\times\! 10^{-10}$ & $2^{-6}$ & $5$    & 3      & 3      \\
1F              & $10^3$   & $\ast$    & $10^{3}$ & $\ast$                   & $2^{-7\phantom{0}}$ & $7.28 \!\times\! 10^{-2}$ & $2.4390 \!\times\! 10^{-10}$ & $2^{-6}$ & $5$    & 3      & 3      \\
1G              & $10^3$   & $10^{-3}$ & $10^{3}$ & $10^{-3}$                & $\ast$              & $7.28 \!\times\! 10^{-2}$ & $\ast$                       & $2^{-6}$ & $5$    & $\ast$ & $\ast$ \\
2\phantom{A}    & $10^3$   & $10^{-3}$ & $1$      & $1.813\!\times\!10^{-5}$ & $2^{-10}$           & $7.28 \!\times\! 10^{-2}$ & $4.7636 \!\times\! 10^{-13}$ & $2^{-7}$ & $5$    & 7      & 5      \\
\hline
\end{tabular}
}
\end{table}

Section \SEC{num_exp_cond} considers the condition number of the linear (sub)systems that occur in the tangent problems in Newton's method for the coupled NSCH system, and its dependence on relevant system parameters. In Section~\SEC{GMRES_conv} we analyse the convergence behavior of the GMRES method with the $\boldsymbol P_{\UPP}$ and $\boldsymbol P_{\LOW}$ preconditioners, corresponding to the two variants of the partitioned solution method discussed in Section~\SEC{partitioned}. Finally, in Section \SEC{num_thin_int},  we present a simulation with a thin diffuse interface and high spatiotemporal resolution, and validate the AGG NSCH model by comparing the obtained results with those of a corresponding sharp -interface model.

\subsection{Conditioning of the tangent matrices}
\label{sec:num_exp_cond}
To assess the complexity of solving the linear tangent problems in~\EQ{Newton} by means of the direct solution method in Section~\SEC{monolithic} and the partitioned solution methods in Section~\SEC{partitioned}, we investigate the dependence of the 1-norm condition number of the tangent matrix $\boldsymbol{A}$ and of the submatrices $\boldsymbol{A}_{\NS\text{-}\NS}$ and~$\boldsymbol{A}_{\CH\text{-}\CH}$ on the parameters of the NSCH model and its discretization. The condition number is an important characteristic in assessing the complexity of linear systems but, evidently, it provides only partial information. To account for the fact that most modern direct solution methods apply an {\em equilibration\/} (reordering and rescaling) of the matrix before entering the solution procedure, we consider both the condition number of the original tangent matrices, and of corresponding equilibrated matrices. In this work, we apply the PARDISO sparse direct solver~\cite{Schenk:2001qo,Schenk:2004eu,Bollhofer:2020hb}. Because the equilibrated systems as generated and used by the PARDISO solver are not accessible externally, we instead consider the systems obtained from MATLAB's \texttt{equilibrate\/} function to assess the effect of equilibration. In view of the size of the linear systems, instead of computing the condition numbers directly, we utilize MATLAB's built-in \texttt{condest} function, which is based on the algorithms described in~\cite{hager1984condition} and~\cite{higham2000block}, to compute a repeatable\footnote{by resetting MATLAB's random number generator seed with the \texttt{rng(`default')} command before every operation that requests pseudo-random numbers} lower bound on the 1-norm condition number.

The analysis is performed on the aforementioned oscillating-droplet test case. We vary the interface thickness parameter~$\varepsilon$ and (minimal) mesh width $h$ such that $\varepsilon / h \geq 1/2$, thus ensuring adequate resolution in the vicinity of the diffuse interface. In addition, we limit ourselves to $\varepsilon \geq 2^{-10}\!\times\!10^{2}$\,\textmu{}m, to restrict the size of the  linear systems to a feasible range. We regard configurations with parameter values as presented for test cases 1A--G in Table~\ref{tab:cases}. For each test case, the entries represented by the $\ast$ symbol indicate parameter(s) to be varied. We restrict our presentation to results for the Crank--Nicolson time discretization according to~\EQ{error_system_CN}. However, the Backward Euler approximation~\EQ{error_system_BE} yields similar behavior (results not displayed). The presented condition-number estimates pertain to the tangent matrices in the last Newton iteration at refinement level $L_{\MAX}$ in the first time step. We have verified that the presented condition-number estimates are representative, by probing several linear systems in later time steps and during the Newton iterations, for various cases (results not displayed).

Table~\TAB{condition_1} presents an overview of the estimated lower bound on the 1-norm condition numbers for test case 1A, exploring the $(\varepsilon, h)$\nobreakdash-parameter dependence. In this test case, we consider uniform meshes comprising $(5\!\times\!5)\,2^k$ elements ($k=0,1,2,3$). The number of DOFs corresponding to the finest mesh exceeds $65\!\times\!10^3$. Preliminary tests for $k=4$ conveyed that the condition-number estimates are unreliable for $k>3$ (results not displayed). Let us note that the empty entries in Table~\TAB{condition_1} correspond to $(\varepsilon, h)$ combinations for which the diffuse interface is not adequately resolved. These entries therefore have no significance. The results  in Table~\TAB{condition_1} convey that the condition number of the aggregated NSCH system is generally significantly larger than that of the NS and CH subsystems separately, especially for small $\varepsilon$ and~$h$. Considering the $h$\nobreakdash-dependence of the condition numbers of the original matrices, one can observe that the condition numbers of the CH and NS subsystem are essentially independent of~$h$, while the condition number of the aggregated NSCH system deteriorates under $h$\nobreakdash-refinement. Considering the effect of equilibration, we note that equilibration is generally very effective, reducing the condition number by orders of magnitude. However, the results indicate that the effectiveness of equilibration decreases under $h$\nobreakdash-refinement. For the CH and NS subsystems, the condition number of the equilibrated matrix increases, while the condition number of the original matrices is essentially independent of~$h$. For the aggregated NSCH system, the ratio of the condition number of the original matrix to the condition number of the equilibrated matrix generally decreases under $h$\nobreakdash-refinement. In general, the condition number of the equilibrated NSCH system deteriorates more strongly under $h$\nobreakdash-refinement than that of the NS and CH subsystems.
Considering the $\varepsilon$\nobreakdash-dependence of the condition numbers, we observe that the condition number of the aggregated NSCH system deteriorates monotonously
as~$\varepsilon$ decreases. This dependence is however weak and, approximately, $\operatorname{cond}(\boldsymbol{B})=O(\varepsilon^{1/2})$ (as $\varepsilon\to+0$).
The condition number of the CH subsystem varies weakly and non-monotonously in relation to~$\varepsilon$. The results in Table~\TAB{condition_1} convey that the condition number of the CH subsystem varies at most by a factor of~$2$ across the full range of considered values of~$\varepsilon$. The estimated condition number of the NS subsystem is essentially independent of~$\varepsilon$, in accordance with the fact that the NS subsystem does not depend on~$\varepsilon$ in the matched-density scenario considered in test case 1A. For sufficiently small~$h$, equilibration leads to monotonous decay of the condition number of the aggregated NSCH system as~$\varepsilon$ decreases. Also in this case, the dependence of the condition number on~$\varepsilon$ is weak. For the CH subsystem, equilibration leads to monotonous decay of the condition number as~$\varepsilon$ decreases. Evidently, the condition number of the equilibrated NS subsystem is essentially independent of~$\varepsilon$, analogous to the condition number of the original NS subsystem. 
\begin{table}
\caption{Lower bounds of the 1-norm condition numbers as computed by MATLAB's \texttt{condest} function for test case 1A in Table~\TAB{cases}. Parameters of interest are the interface thickness $\varepsilon$ and mesh width $h$. No local mesh refinement has been applied. The labels `\textrm{orig}' and `\textrm{equi}' indicate the condition numbers of the matrices before and after equilibration, respectively. The monolithic system is represented by `{full}'. The Cahn--Hilliard and Navier--Stokes subsystems are indicated by `{CH}' and `{NS}', respectively. The blue shaded entries correspond to $(\varepsilon,h)$ combinations for which $\varepsilon/h=2^{-3}\!\times\!10$, and the diffuse interface is marginally resolved.
\label{tab:condition_1}}
\resizebox{\columnwidth}{!}{
\begin{tabular}{ll|cccccccc|}
\cline{3-10}\multicolumn{1}{c}{}
 &  & \multicolumn{2}{c|}{$h=2^0\!\times\!10\,\text{\textmu m}$} & \multicolumn{2}{c|}{$h=2^{-1}\!\times\!10\,\text{\textmu m}$} & \multicolumn{2}{c|}{$h=2^{-2}\!\times\!10\,\text{\textmu m}$} & \multicolumn{2}{c|}{$h=2^{-3}\!\times\!10\,\text{\textmu m}$} \\
 &  & \multicolumn{1}{c}{\textrm{orig}}                          & \multicolumn{1}{c|}{\textrm{equi}}                            & \multicolumn{1}{c}{\textrm{orig}}                             & \multicolumn{1}{c|}{\textrm{equi}}                            & \multicolumn{1}{c}{\textrm{orig}} & \multicolumn{1}{c|}{\textrm{equi}} & \multicolumn{1}{c}{\textrm{orig}} & \multicolumn{1}{c|}{\textrm{equi}} \\ \hline
\multicolumn{1}{|l}{}                                                    & {full} & \cellcolor{blue!10}3.01E+11 & \multicolumn{1}{l|}{\cellcolor{blue!10}1.95E+04} & 1.76E+11 & \multicolumn{1}{l|}{5.27E+04} & 4.15E+11 & \multicolumn{1}{l|}{5.87E+05} & 1.48E+12 & \multicolumn{1}{l|}{3.29E+07} \\
\multicolumn{1}{|c}{$\varepsilon=2^{-3}\!\times\!100\,\text{\textmu m}$} & {CH}   & \cellcolor{blue!10}2.99E+08 & \multicolumn{1}{l|}{\cellcolor{blue!10}3.54E+02} & 3.58E+07 & \multicolumn{1}{l|}{6.77E+02} & 3.16E+07 & \multicolumn{1}{l|}{6.64E+03} & 5.05E+07 & \multicolumn{1}{l|}{9.76E+04} \\
\multicolumn{1}{|l}{}                                                    & {NS}   & \cellcolor{blue!10}3.73E+10 & \multicolumn{1}{l|}{\cellcolor{blue!10}1.91E+04} & 2.39E+10 & \multicolumn{1}{l|}{5.05E+04} & 2.08E+10 & \multicolumn{1}{l|}{1.04E+05} & 2.45E+10 & \multicolumn{1}{l|}{9.50E+05} \\ \hline
\multicolumn{1}{|l}{}                                                    & {full} & 5.34E+11 & \multicolumn{1}{l|}{1.95E+04} &  \cellcolor{blue!10}3.12E+11 & \multicolumn{1}{l|}{\cellcolor{blue!10}5.53E+04} & 5.43E+11 & \multicolumn{1}{l|}{2.78E+05} & 1.73E+12 & \multicolumn{1}{l|}{1.67E+07} \\
\multicolumn{1}{|c}{$\varepsilon=2^{-4}\!\times\!100\,\text{\textmu m}$} & {CH}   & 5.35E+08 & \multicolumn{1}{l|}{5.31E+02} &  \cellcolor{blue!10}4.86E+07 & \multicolumn{1}{l|}{\cellcolor{blue!10}4.96E+02} & 2.73E+07 & \multicolumn{1}{l|}{3.40E+03} & 3.52E+07 & \multicolumn{1}{l|}{4.97E+04} \\
\multicolumn{1}{|l}{}                                                    & {NS}   & 3.73E+10 & \multicolumn{1}{l|}{1.94E+04} &  \cellcolor{blue!10}2.39E+10 & \multicolumn{1}{l|}{\cellcolor{blue!10}5.05E+04} & 2.08E+10 & \multicolumn{1}{l|}{1.04E+05} & 2.45E+10 & \multicolumn{1}{l|}{9.50E+05} \\ \hline
\multicolumn{1}{|l}{}                                                    & {full} &          & \multicolumn{1}{l|}{}         & 4.69E+11 & \multicolumn{1}{l|}{5.61E+04} & \cellcolor{blue!10}8.04E+11 & \multicolumn{1}{l|}{\cellcolor{blue!10}2.08E+05} & 2.22E+12 & \multicolumn{1}{l|}{8.16E+06} \\
\multicolumn{1}{|c}{$\varepsilon=2^{-5}\!\times\!100\,\text{\textmu m}$} & {CH}   &          & \multicolumn{1}{l|}{}         & 8.18E+07 & \multicolumn{1}{l|}{4.52E+02} & \cellcolor{blue!10}2.65E+07 & \multicolumn{1}{l|}{\cellcolor{blue!10}1.84E+03} & 2.97E+07 & \multicolumn{1}{l|}{2.69E+04} \\
\multicolumn{1}{|l}{}                                                    & {NS}   &          & \multicolumn{1}{l|}{}         & 2.39E+10 & \multicolumn{1}{l|}{5.05E+04} & \cellcolor{blue!10}2.08E+10 & \multicolumn{1}{l|}{\cellcolor{blue!10}1.04E+05} & 2.45E+10 & \multicolumn{1}{l|}{9.50E+05} \\ \hline
\multicolumn{1}{|l}{}                                                    & {full} &          & \multicolumn{1}{l|}{}         &          & \multicolumn{1}{l|}{}         & 1.18E+12 & \multicolumn{1}{l|}{1.08E+05}   & \cellcolor{blue!10}2.70E+12 & \multicolumn{1}{l|}{\cellcolor{blue!10}3.76E+06} \\
\multicolumn{1}{|c}{$\varepsilon=2^{-6}\!\times\!100\,\text{\textmu m}$} & {CH}   &          & \multicolumn{1}{l|}{}         &          & \multicolumn{1}{l|}{}         & 3.22E+07 & \multicolumn{1}{l|}{1.18E+03}   & \cellcolor{blue!10}3.21E+07 & \multicolumn{1}{l|}{\cellcolor{blue!10}1.62E+04} \\
\multicolumn{1}{|l}{}                                                    & {NS}   &          & \multicolumn{1}{l|}{}         &          & \multicolumn{1}{l|}{}         & 2.08E+10 & \multicolumn{1}{l|}{1.03E+05}   & \cellcolor{blue!10}2.40E+10 & \multicolumn{1}{l|}{\cellcolor{blue!10}9.50E+05} \\ \hline
\multicolumn{1}{|l}{}                                                    & {full} &          & \multicolumn{1}{l|}{}         &          & \multicolumn{1}{l|}{}         &          & \multicolumn{1}{l|}{}           & 5.06E+12 & \multicolumn{1}{l|}{3.17E+06} \\
\multicolumn{1}{|c}{$\varepsilon=2^{-7}\!\times\!100\,\text{\textmu m}$} & {CH}   &          & \multicolumn{1}{l|}{}         &          & \multicolumn{1}{l|}{}         &          & \multicolumn{1}{l|}{}           & 5.37E+07 & \multicolumn{1}{l|}{1.32E+04} \\
\multicolumn{1}{|l}{}                                                    & {NS}   &          & \multicolumn{1}{l|}{}         &          & \multicolumn{1}{l|}{}         &          & \multicolumn{1}{l|}{}           & 2.40E+10 & \multicolumn{1}{l|}{9.50E+05} \\ \hline
\end{tabular}
}
\end{table}

In actual computations, the mesh width~$h$ must be reduced if the interface-thickness parameter~$\varepsilon$ is decreased, in order to ensure an appropriate resolution of the diffuse interface. The blue-shaded diagonal entries of Table~\TAB{condition_1} correspond to combinations of $(\varepsilon,h)$ for which~$\varepsilon/h=2^{-3}\!\times\!10$ and the diffuse-interface is marginally resolved. It must be noted, however, that in a typical application of the NSCH system, one would select the mobility~$m$ proportional to $\varepsilon^3$, such that the diffuse time scale is fixed; see test case 1G and~\APP{time_scale}. These variations in~$m$ also affect the condition number; see e.g.\ Figure~\FIG{results_m} below. Therefore, we restrict ourselves here to a basic characterization of the behavior of the equilibrated condition number under simultaneous $\varepsilon$ and $h$ refinement. From Table~\TAB{condition_1}, one can infer that the condition number of the equilibrated NSCH tangent matrix scales  approximately as $\varepsilon^{-5/2}$ under simultaneous $\varepsilon$ and $h\propto\varepsilon$ refinement, and the condition numbers of the equilibrated NS and CH subsystems scale approximately as~$\varepsilon^{-2}$. 

Table~\TAB{condition_2} presents the estimated lower bound on the 1-norm condition numbers for test case~1B in Table~\TAB{cases}. Test case 1B pertains to an adaptively refined approximation, in which the interface-thickness parameter $\varepsilon$ and the number of refinement levels $L_{\MAX}$ are varied. The coarsest mesh comprises 
$10\!\times\!10$ square elements with mesh width $h_0=5\,\text{\textmu m}$. For completeness, we mention that we apply $K=L_{\MAX}$ continuation steps, and set (here and throughout) the refinement factor to $a^f=0.95$ for $f\in\{\varphi,\boldsymbol{u}\}$ and $a^f=0$ for $f\in\{p,\mu\}$; see Section~\SEC{error_adaptive}. Table~\TAB{condition_2} reports the condition numbers corresponding to the super-refinement of the maximum refinement level.
To relate the adaptive-approximation results in Table~\TAB{condition_2} to the uniform-mesh results in Table~\TAB{condition_1}, we note that for $L_{\MAX}=1$ the minimal element size is $2^{-2}\!\times\!5\,\text{\textmu m}$, which coincides with the element size in the right-most column of Table~\TAB{condition_1}. Comparing the condition numbers in left-most column of Table~\TAB{condition_2} to those in the right-most column of Table~\TAB{condition_1} (for identical values of~$\varepsilon$), we observe that the adaptive approximations generally exhibits moderately worse condition numbers for the original tangent matrices, and moderately better condition numbers for the equilibrated 
matrices. The condition numbers of the aggregated NSCH system, and of the NS and CH subsystems, generally deteriorate as the number of refinement levels $L_{\MAX}$ increases.
In addition, the effectiveness of equilibration diminishes as $L_{\MAX}$ increases, in the sense that the ratio of the condition number of the original matrix to the condition number of the equilibrated matrix decreases. The condition numbers of the equilibrated NSCH system and the NS and CH subsystems increase approximately proportional to $2^{\frac{7}{2}L_{\MAX}}$, $2^{L_{\MAX}}$, and $2^{\frac{7}{2}L_{\MAX}}$, respectively. 
As opposed to the case of uniform refinements, the estimated condition number for the NS subsystem exhibits a (weak) dependence 
on~$\varepsilon$. This is caused by the fact that the adaptively-refined meshes implicitly depend on~$\varepsilon$. Considering the $\varepsilon$\nobreakdash-dependence of the equilibrated condition numbers, one can observe that the condition numbers of the equilibrated NSCH system and of the equilibrated CH subsystem display a weak monotonous decay as~$\varepsilon$ decreases, approximately proportional to $\varepsilon^{1/2}$, while the condition number of the equilibrated NS system is essentially independent of~$\varepsilon$. 
\begin{table}
\caption{Lower bounds of the 1-norm condition numbers as computed by MATLAB's \texttt{condest} function for test case 1B in Table~\TAB{cases}. Parameters of interest are the interface thickness $\varepsilon$, and the number of refinements $L_{max}$. The labels `\textrm{orig}' and `\textrm{equi}' indicate the condition numbers of the matrices before and after equilibration, respectively. The monolithic system is represented by `{full}'. The Cahn--Hilliard and Navier--Stokes subsystems are indicated by `{CH}' and `{NS}', respectively. The blue shaded entries correspond to combinations of $(\varepsilon,L_{\MAX})$ for which the ratio $\varepsilon/h$ of $\varepsilon$ to the minimal element size of the super-refinement $h=2^{-(L_{\MAX}+1)}h_0$ equals $2^{-3}\!\times\!10$, and the diffuse interface is marginally resolved. The yellow highlighted entry corresponds to settings that are repeated in all test cases 1B--F.
\label{tab:condition_2}}
\resizebox{\columnwidth}{!}{
\begin{tabular}{ll|cccccccc|}
\cline{3-10}
 &  & \multicolumn{2}{c|}{$L_{\MAX}=1$} & \multicolumn{2}{c|}{$L_{\MAX}=2$}  & \multicolumn{2}{c|}{$L_{\MAX}=3$} & \multicolumn{2}{c|}{$L_{\MAX}=4$}  \\
 &  & \multicolumn{1}{c}{\textrm{orig}} & \multicolumn{1}{c|}{\textrm{equi}} & \multicolumn{1}{c}{\textrm{orig}} & \multicolumn{1}{c|}{\textrm{equi}} & \multicolumn{1}{c}{\textrm{orig}} & \multicolumn{1}{c|}{\textrm{equi}} & \multicolumn{1}{c}{\textrm{orig}} & \multicolumn{1}{c|}{\textrm{equi}} \\ \hline
\multicolumn{1}{|l}{}                                                     & {full} & 2.12E+12 & \multicolumn{1}{l|}{2.33E+06} & 1.70E+15 & \multicolumn{1}{l|}{2.41E+09} & 5.68E+16 & \multicolumn{1}{l|}{3.10E+10} & 8.33E+13 & \multicolumn{1}{l|}{7.87E+09} \\
\multicolumn{1}{|c}{$\varepsilon=2^{-5}\!\times\!100\,\text{\textmu m}$}  & {CH}   & 9.41E+07 & \multicolumn{1}{l|}{3.42E+04} & 3.58E+08 & \multicolumn{1}{l|}{5.35E+05} & 1.37E+09 & \multicolumn{1}{l|}{9.25E+06} & 3.62E+09 & \multicolumn{1}{l|}{8.84E+07} \\
\multicolumn{1}{|l}{}                                                     & {NS}   & 4.52E+10 & \multicolumn{1}{l|}{2.00E+05} & 1.39E+11 & \multicolumn{1}{l|}{5.99E+05} & 4.43E+11 & \multicolumn{1}{l|}{2.52E+06} & 9.89E+11 & \multicolumn{1}{l|}{4.89E+06} \\ \hline
\multicolumn{1}{|l}{}                                                     & {full} & \cellcolor{blue!10}3.04E+12 & \multicolumn{1}{l|}{\cellcolor{blue!10}1.21E+06} & 9.04E+12 & \multicolumn{1}{l|}{1.42E+07} & 6.60E+16 & \multicolumn{1}{l|}{6.83E+09} & 1.28E+14 & \multicolumn{1}{l|}{8.37E+09} \\
\multicolumn{1}{|c}{$\varepsilon=2^{-6}\!\times\!100\,\text{\textmu m}$}  & {CH}   & \cellcolor{blue!10}1.11E+08 & \multicolumn{1}{l|}{\cellcolor{blue!10}1.95E+04} & 4.16E+08 & \multicolumn{1}{l|}{3.34E+05} & 1.61E+09 & \multicolumn{1}{l|}{5.72E+06} & 5.35E+09 & \multicolumn{1}{l|}{7.67E+07} \\
\multicolumn{1}{|l}{}                                                     & {NS}   & \cellcolor{blue!10}4.26E+10 & \multicolumn{1}{l|}{\cellcolor{blue!10}1.80E+05} & 1.05E+11 & \multicolumn{1}{l|}{3.51E+05} & 3.75E+11 & \multicolumn{1}{l|}{1.14E+06} & 4.67E+18 & \multicolumn{1}{l|}{2.73E+06} \\ \hline
\multicolumn{1}{|l}{}                                                     & {full} & 5.58E+12 & \multicolumn{1}{l|}{8.40E+05} & \cellcolor{blue!10}1.51E+13 & \multicolumn{1}{l|}{\cellcolor{blue!10}1.07E+07} & \cellcolor{yellow!10}4.99E+13 & \multicolumn{1}{l|}{\cellcolor{yellow!10}1.55E+08} & 1.62E+14 & \multicolumn{1}{l|}{3.76E+10} \\
\multicolumn{1}{|c}{$\varepsilon=2^{-7}\!\times\!100\,\text{\textmu m}$}  & {CH}   & 1.96E+08 & \multicolumn{1}{l|}{1.36E+04} & \cellcolor{blue!10}6.61E+08 & \multicolumn{1}{l|}{\cellcolor{blue!10}3.00E+05} & \cellcolor{yellow!10}2.30E+09 & \multicolumn{1}{l|}{\cellcolor{yellow!10}4.19E+06} & 8.56E+09 & \multicolumn{1}{l|}{6.31E+07} \\
\multicolumn{1}{|l}{}                                                     & {NS}   & 4.20E+10 & \multicolumn{1}{l|}{1.29E+05} & \cellcolor{blue!10}9.46E+10 & \multicolumn{1}{l|}{\cellcolor{blue!10}3.20E+05} & \cellcolor{yellow!10}2.81E+11 & \multicolumn{1}{l|}{\cellcolor{yellow!10}5.22E+05} & 1.73E+16 & \multicolumn{1}{l|}{2.04E+06} \\ \hline
\multicolumn{1}{|l}{}                                                     & {full} &          & \multicolumn{1}{l|}{}         & 3.08E+13 & \multicolumn{1}{l|}{9.36E+06}   & \cellcolor{blue!10}8.06E+13 & \multicolumn{1}{l|}{\cellcolor{blue!10}1.08E+08}   & 2.54E+14 & \multicolumn{1}{l|}{1.43E+09} \\
\multicolumn{1}{|c}{$\varepsilon=2^{-8}\!\times\!100\,\text{\textmu m}$}  & {CH}   &          & \multicolumn{1}{l|}{}         & 1.38E+09 & \multicolumn{1}{l|}{2.04E+05}   & \cellcolor{blue!10}3.79E+09 & \multicolumn{1}{l|}{\cellcolor{blue!10}3.23E+06}   & 1.31E+10 & \multicolumn{1}{l|}{4.27E+07} \\
\multicolumn{1}{|l}{}                                                     & {NS}   &          & \multicolumn{1}{l|}{}         & 8.03E+10 & \multicolumn{1}{l|}{4.68E+05}   & \cellcolor{blue!10}2.12E+11 & \multicolumn{1}{l|}{\cellcolor{blue!10}8.96E+05}   & 7.26E+11 & \multicolumn{1}{l|}{1.03E+06} \\ \hline
\multicolumn{1}{|l}{}                                                     & {full} &          & \multicolumn{1}{l|}{}         &          & \multicolumn{1}{l|}{}           & 1.58E+14 & \multicolumn{1}{l|}{9.81E+07}   & \cellcolor{blue!10}4.01E+14 & \multicolumn{1}{l|}{\cellcolor{blue!10}9.95E+08} \\
\multicolumn{1}{|c}{$\varepsilon=2^{-9}\!\times\!100\,\text{\textmu m}$}  & {CH}   &          & \multicolumn{1}{l|}{}         &          & \multicolumn{1}{l|}{}           & 8.46E+09 & \multicolumn{1}{l|}{2.17E+06}   & \cellcolor{blue!10}2.03E+10 & \multicolumn{1}{l|}{\cellcolor{blue!10}2.99E+07} \\
\multicolumn{1}{|l}{}                                                     & {NS}   &          & \multicolumn{1}{l|}{}         &          & \multicolumn{1}{l|}{}           & 1.92E+11 & \multicolumn{1}{l|}{8.32E+05}   & \cellcolor{blue!10}5.52E+11 & \multicolumn{1}{l|}{\cellcolor{blue!10}1.42E+06} \\ \hline
\multicolumn{1}{|l}{}                                                     & {full} &          & \multicolumn{1}{l|}{}         &          & \multicolumn{1}{l|}{}           &          & \multicolumn{1}{l|}{}           & 7.55E+14 & \multicolumn{1}{l|}{9.24E+08} \\
\multicolumn{1}{|c}{$\varepsilon=2^{-10}\!\times\!100\,\text{\textmu m}$} & {CH}   &          & \multicolumn{1}{l|}{}         &          & \multicolumn{1}{l|}{}           &          & \multicolumn{1}{l|}{}           & 4.78E+10 & \multicolumn{1}{l|}{2.39E+07} \\
\multicolumn{1}{|l}{}                                                     & {NS}   &          & \multicolumn{1}{l|}{}         &          & \multicolumn{1}{l|}{}           &          & \multicolumn{1}{l|}{}           & 4.97E+11 & \multicolumn{1}{l|}{1.19E+06} \\ \hline
\end{tabular}
}
\end{table}

Regarding the scenario in which $\varepsilon$ is reduced and $L_{\MAX}$ is simultaneously increased to ensure proper resolution of the diffuse interface, as indicated by the blue-shaded entries in Table~\TAB{condition_2}, one can infer that the condition number of the equilibrated NSCH tangent matrix scales approximately as~$\varepsilon^{-3}$, the condition number of the equilibrated NS subsystem scales approximately as~$\varepsilon^{-1}$, and the condition number of the equilibrated CH subsystem scales approximately as~$\varepsilon^{-3}$ under simultaneous $\varepsilon$ and $h\propto\varepsilon$ refinement, asymptotically as $\varepsilon\to+0$. These proportionalities are consistent with the previously established behavior of the equilibrated condition numbers with respect to~$\varepsilon$ and~$L_{\MAX}$ separately.

In test case 1C we vary the mobility $m$ to determine its effect on the conditioning of the tangent matrices. Other parameters for this test case are listed in Table~\TAB{cases}. Figure~\FIG{results_m} plots the condition numbers of the aggregated NSCH system and of the NS and CH subsystems, and of their equilibrated counterparts, versus the mobility parameter. For $m \leq 4\!\times\! 10^{-13}\,\textrm{m}^d\,\textrm{s}\,\textrm{kg}^{-1} $, the diffusive time scale associated with the 
interface-equilibration process exceeds the period of oscillation of the droplet; see Section~\SEC{num_thin_int}  and~\APP{time_scale}. Such small values of $m$ therefore have limited significance in this setting. On the other hand, for $m \geq 10^{-7} \textrm{m}^d\,\textrm{s}\,\textrm{kg}^{-1}$, the diffusive interface-equilibration process corresponding to the Cahn--Hilliard equations is excessively dissipative, leading to strong artificial damping of the droplet oscillation. Such large values of~$m$ are therefore also irrelevant. From Figure~\FIG{results_m} one can observe that the condition numbers of the original NSCH tangent matrix and of the CH subsystem are essentially independent of~$m$. The condition number of the original NS subsystem is completely independent of~$m$, as is the condition number of the equilibrated NS subsystem, in agreement with the fact that NS subsystem does not contain~$m$. Equilibration causes the condition numbers of the NSCH system and of the CH subsystem to reduce as~$m$ decreases. In particular, Figure~\FIG{results_m} indicates that for $4\!\times\! 10^{-13}\,\textrm{m}^d\,\textrm{s}\,\textrm{kg}^{-1}\leq{}m\leq 10^{-7}\,\textrm{m}^d\,\textrm{s}\,\textrm{kg}^{-1}$, the condition numbers of the equilibrated NSCH system and CH subsystem are approximately proportional to~$m$. This proportionality ceases for even smaller values of the mobility but, as explained previously, such small values of~$m$ have limited relevance.

\begin{figure}[h!!!!!!!!!!]
\centering
\includegraphics[width=\textwidth, trim = 0mm 1mm 12mm 9mm, clip]{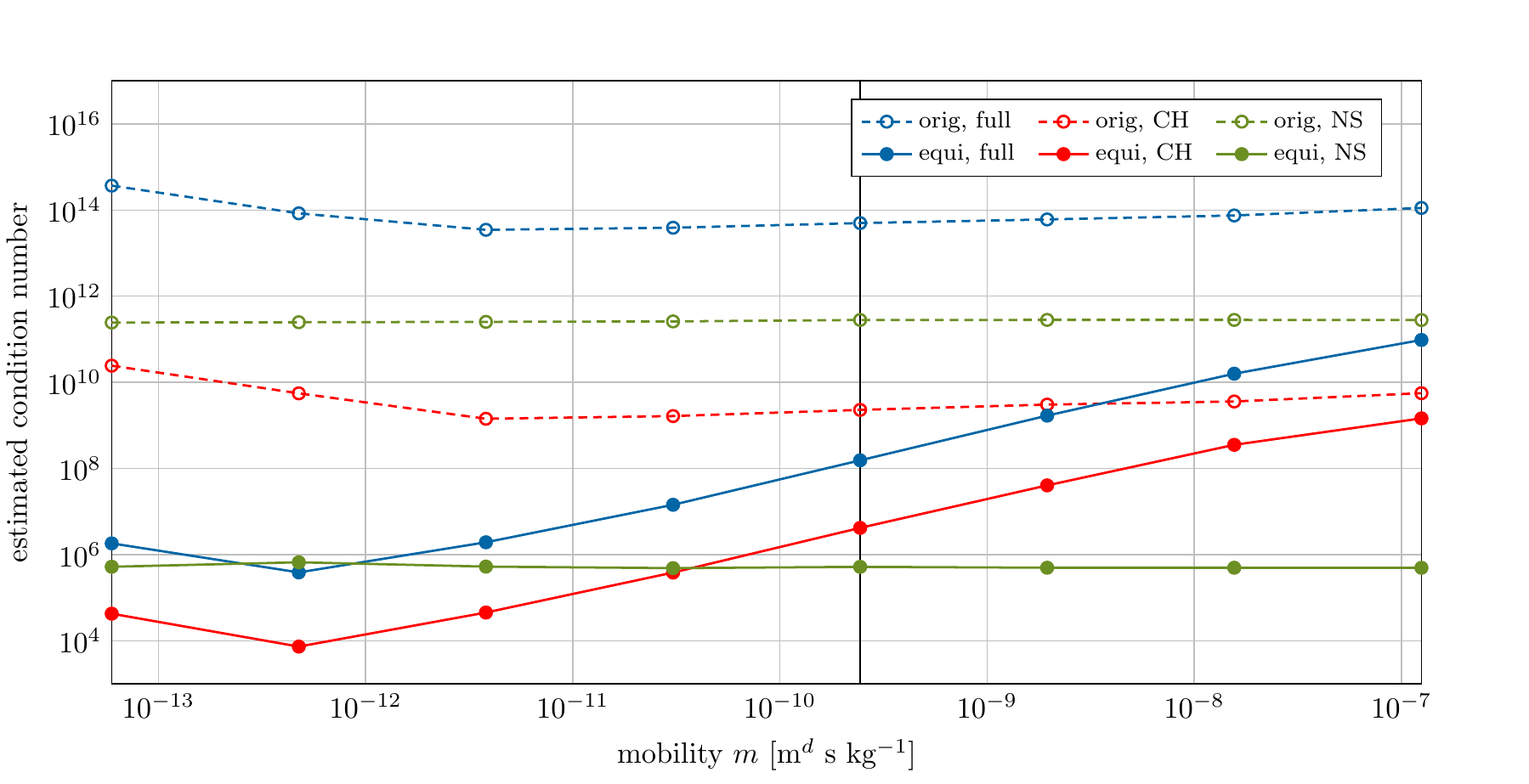}
\caption{Lower bounds of the 1-norm condition numbers as computed by MATLAB's \texttt{condest} function for test case \#1C in Table~\TAB{cases}. Parameter of interest is the mobility coefficient $m$. Matrix naming and computational application are as before in Tables~\TAB{condition_1} and \TAB{condition_2}. The vertical line corresponds to settings that are repeated in all test cases 1B--F.}
\label{fig:results_m}
\end{figure}

\begin{figure}[h!!!!!!!!!!]
\centering
\includegraphics[width=\textwidth, trim = 0mm 1mm 10mm 9mm, clip]{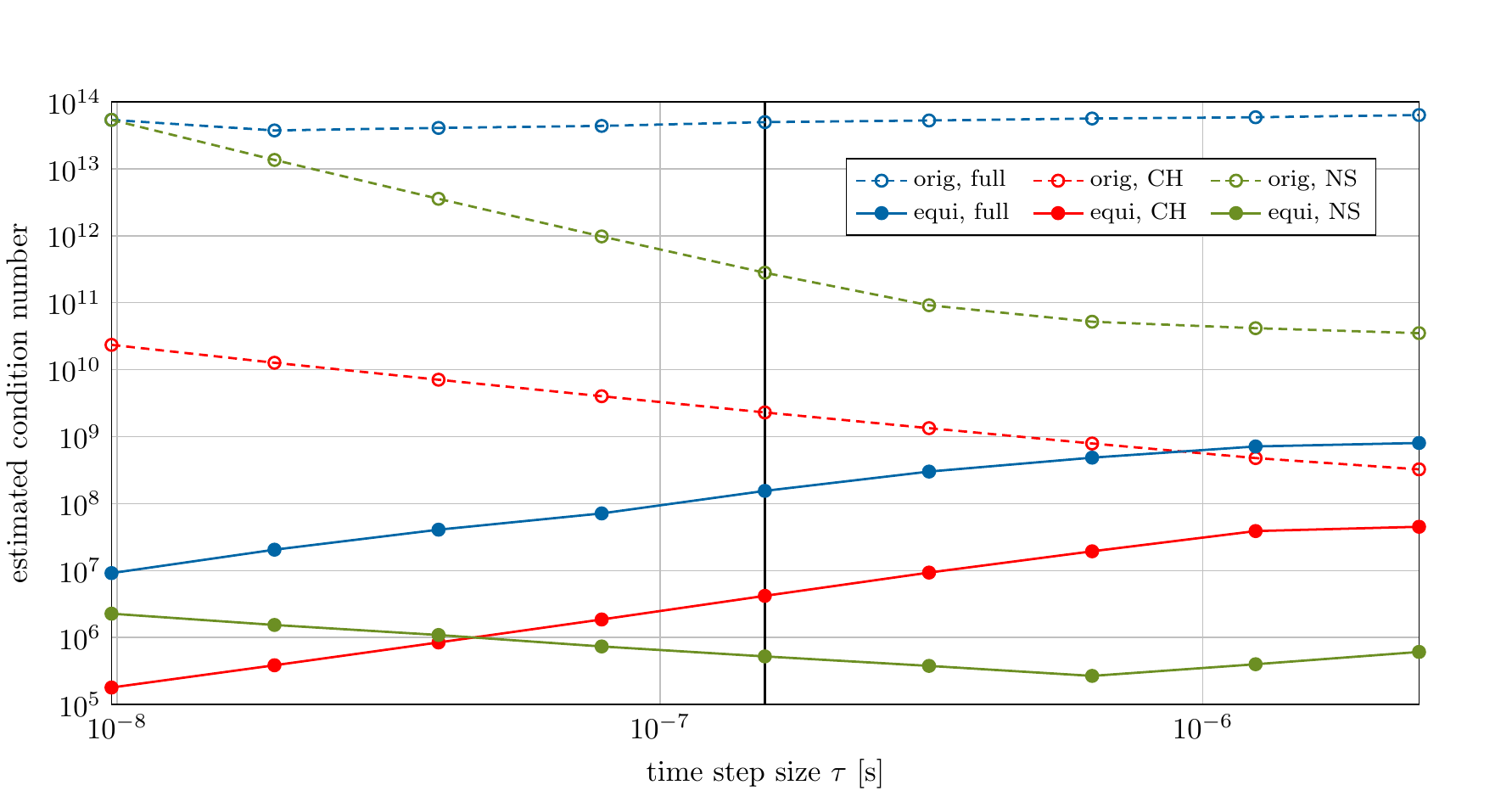}
\caption{Lower bounds of the 1-norm condition numbers as computed by MATLAB's \texttt{condest} function for test case \#1D in Table~\TAB{cases}. Parameter of interest is the time step size $\tau$. Matrix naming and computational application are as before in Tables~\ref{tab:condition_1} and \ref{tab:condition_2}, and Figure~\FIG{results_m}. The vertical line corresponds to settings that are repeated in test cases 1B--F.}
\label{fig:results_Dt}
\end{figure}

Test case 1D is concerned with the dependence of the condition numbers on the time-step size $\tau$; see Table~\TAB{cases}. The corresponding results are presented in Figure~\FIG{results_Dt}. One can observe that the condition number of the original aggregated NSCH tangent matrix is essentially independent of~$\tau$. The condition number of the original NS subsystem increases as the time step is reduced, approximately proportional to~$\tau^{-2}$. For the CH subsystem, the condition number increases approximately proportional to~$\tau^{-1}$ under time-step refinement. Equilibration effectively improves the conditioning of the NSCH tangent matrix and of the NS and CH subsystems. 
The condition numbers of the equilibrated NSCH tangent matrix and of the equilibrated CH subsystem decrease approximately proportional to~$\tau$. The condition number of the equilibrated NS subsystem increases approximately as~$\tau^{-1/2}$ as~$\tau$ decreases.

Test case 1E pertains to the dependence of the condition numbers on variations in the density. The results of test case 1E are presented in Table~\TAB{condition_3}. 
We consider liquid and ambient densities in the physically relevant range $10^0\,\mathrm{kg}\,\mathrm{m}^{-d}\leq\rho_{\LL},\rho_{\AA}\leq{}10^4\,\mathrm{kg}\,\mathrm{m}^{-d}$. The results in Table~\TAB{condition_3} indicate that the condition numbers of the aggregated NSCH tangent matrix and of the NS and CH subsystems generally depend only weakly on~$\rho_{\LL}$ and~$\rho_{\AA}$, both for the original matrices and for their equilibrated counterparts. Because the densities do not appear explicitly in the CH subsystem, one would expect that the condition number of the CH subsystem is entirely independent of~$\rho_{\LL}$ and~$\rho_{\AA}$. Nonetheless, the condition numbers of both the original and the equilibrated CH subsystem vary by a factor of approximately~10 across the considered range of~$(\rho_{\LL},\rho_{\AA})$.
These variations are indirect, induced by e.g.\ differences in the adaptively-refined meshes.
\begin{table}[!htbp]
\caption{Lower bounds of the 1-norm condition numbers as computed by MATLAB's \texttt{condest} function for test case 1E in Table~\TAB{cases}. Parameters of interest are the liquid and ambient densities $\rho_{\LL}$ and $\rho_{\AA}$.}
\label{tab:condition_3}
\resizebox{\columnwidth}{!}{
\begin{tabular}{ll|cccccccccc|}
\cline{3-12}
 & & \multicolumn{2}{c|}{$\rho_{\AA}=10^0\,\frac{\text{kg}}{\text{m}^d}$}  & \multicolumn{2}{c|}{$\rho_{\AA}=10^1\,\frac{\text{kg}}{\text{m}^d}$}  & \multicolumn{2}{c|}{$\rho_{\AA}=10^2\,\frac{\text{kg}}{\text{m}^d}$}  & \multicolumn{2}{c|}{$\rho_{\AA}=10^3\,\frac{\text{kg}}{\text{m}^d}$}  & \multicolumn{2}{c|}{$\rho_{\AA}=10^4\,\frac{\text{kg}}{\text{m}^d}$} \\
 & & orig & \multicolumn{1}{c|}{equi}                                      & orig & \multicolumn{1}{c|}{equi}                                      & orig & \multicolumn{1}{c|}{equi}                                      & orig & \multicolumn{1}{c|}{equi}                                      & orig & \multicolumn{1}{c|}{equi}                                 \\ \hline
\multicolumn{1}{|l}{}                                                & full & 4.81E+13 & \multicolumn{1}{l|}{1.95E+08} & 4.82E+13 & \multicolumn{1}{l|}{1.44E+08} & 4.74E+13 & \multicolumn{1}{l|}{1.39E+08} & 4.86E+13 & \multicolumn{1}{l|}{1.55E+08} & 5.08E+13 & 1.31E+08 \\
\multicolumn{1}{|c}{$\rho_{\LL}=10^0\,\frac{\text{kg}}{\text{m}^d}$} & CH   & 2.47E+09 & \multicolumn{1}{l|}{2.22E+06} & 2.40E+09 & \multicolumn{1}{l|}{2.78E+06} & 2.34E+09 & \multicolumn{1}{l|}{3.84E+06} & 2.38E+09 & \multicolumn{1}{l|}{4.48E+06} & 2.35E+09 & 4.26E+06 \\
\multicolumn{1}{|l}{}                                                & NS   & 2.53E+10 & \multicolumn{1}{l|}{6.64E+06} & 2.92E+10 & \multicolumn{1}{l|}{3.15E+06} & 3.77E+10 & \multicolumn{1}{l|}{1.23E+06} & 2.55E+11 & \multicolumn{1}{l|}{6.73E+05} & 2.06E+13 & 2.59E+06 \\ \hline
\multicolumn{1}{|l}{}                                                & full & 4.74E+13 & \multicolumn{1}{l|}{1.21E+08} & 4.90E+13 & \multicolumn{1}{l|}{1.36E+08} & 4.89E+13 & \multicolumn{1}{l|}{1.40E+08} & 5.00E+13 & \multicolumn{1}{l|}{1.40E+08} & 5.12E+13 & 1.32E+08 \\
\multicolumn{1}{|c}{$\rho_{\LL}=10^1\,\frac{\text{kg}}{\text{m}^d}$} & CH   & 2.41E+09 & \multicolumn{1}{l|}{3.29E+06} & 2.36E+09 & \multicolumn{1}{l|}{3.43E+06} & 2.33E+09 & \multicolumn{1}{l|}{3.84E+06} & 2.38E+09 & \multicolumn{1}{l|}{4.38E+06} & 2.35E+09 & 4.26E+06 \\
\multicolumn{1}{|l}{}                                                & NS   & 2.84E+10 & \multicolumn{1}{l|}{3.80E+06} & 2.91E+10 & \multicolumn{1}{l|}{1.88E+06} & 3.79E+10 & \multicolumn{1}{l|}{9.39E+05} & 2.53E+11 & \multicolumn{1}{l|}{6.26E+05} & 2.07E+13 & 2.63E+06 \\ \hline
\multicolumn{1}{|l}{}                                                & full & 4.70E+13 & \multicolumn{1}{l|}{1.32E+08} & 4.69E+13 & \multicolumn{1}{l|}{1.38E+08} & 4.83E+13 & \multicolumn{1}{l|}{1.62E+08} & 5.05E+13 & \multicolumn{1}{l|}{1.57E+08} & 5.09E+13 & 1.31E+08 \\
\multicolumn{1}{|c}{$\rho_{\LL}=10^2\,\frac{\text{kg}}{\text{m}^d}$} & CH   & 2.32E+09 & \multicolumn{1}{l|}{3.69E+06} & 2.30E+09 & \multicolumn{1}{l|}{3.67E+06} & 2.30E+09 & \multicolumn{1}{l|}{4.27E+06} & 2.36E+09 & \multicolumn{1}{l|}{4.41E+06} & 2.35E+09 & 4.26E+06 \\
\multicolumn{1}{|l}{}                                                & NS   & 2.58E+10 & \multicolumn{1}{l|}{1.99E+06} & 2.61E+10 & \multicolumn{1}{l|}{1.18E+06} & 3.51E+10 & \multicolumn{1}{l|}{4.79E+05} & 2.58E+11 & \multicolumn{1}{l|}{5.10E+05} & 2.05E+13 & 2.56E+06 \\ \hline
\multicolumn{1}{|l}{}                                                & full & 4.89E+13 & \multicolumn{1}{l|}{1.63E+08} & 4.94E+13 & \multicolumn{1}{l|}{1.62E+08} & 4.84E+13 & \multicolumn{1}{l|}{1.55E+08} & \cellcolor{yellow!10}4.99E+13 & \multicolumn{1}{l|}{\cellcolor{yellow!10}1.55E+08} & 4.99E+13 & 1.29E+08 \\
\multicolumn{1}{|c}{$\rho_{\LL}=10^3\,\frac{\text{kg}}{\text{m}^d}$} & CH   & 2.27E+09 & \multicolumn{1}{l|}{4.16E+06} & 2.27E+09 & \multicolumn{1}{l|}{4.17E+06} & 2.27E+09 & \multicolumn{1}{l|}{4.18E+06} & \cellcolor{yellow!10}2.30E+09 & \multicolumn{1}{l|}{\cellcolor{yellow!10}4.19E+06} & 2.31E+09 & 4.14E+06 \\
\multicolumn{1}{|l}{}                                                & NS   & 3.43E+10 & \multicolumn{1}{l|}{9.34E+05} & 2.69E+10 & \multicolumn{1}{l|}{6.31E+05} & 4.18E+10 & \multicolumn{1}{l|}{2.92E+05} & \cellcolor{yellow!10}2.81E+11 & \multicolumn{1}{l|}{\cellcolor{yellow!10}5.22E+05} & 2.10E+13 & 2.56E+06 \\ \hline
\multicolumn{1}{|l}{}                                                & full & 4.73E+13 & \multicolumn{1}{l|}{2.24E+08} & 4.79E+13 & \multicolumn{1}{l|}{2.21E+08} & 4.86E+13 & \multicolumn{1}{l|}{2.06E+08} & 4.76E+13 & \multicolumn{1}{l|}{1.61E+08} & 4.86E+13 & 1.27E+08 \\
\multicolumn{1}{|c}{$\rho_{\LL}=10^4\,\frac{\text{kg}}{\text{m}^d}$} & CH   & 2.28E+09 & \multicolumn{1}{l|}{3.63E+06} & 2.28E+09 & \multicolumn{1}{l|}{3.63E+06} & 2.28E+09 & \multicolumn{1}{l|}{3.84E+06} & 2.28E+09 & \multicolumn{1}{l|}{3.83E+06} & 2.27E+09 & 4.13E+06 \\
\multicolumn{1}{|l}{}                                                & NS   & 6.93E+11 & \multicolumn{1}{l|}{1.17E+06} & 6.92E+11 & \multicolumn{1}{l|}{1.26E+06} & 6.96E+11 & \multicolumn{1}{l|}{1.12E+06} & 1.15E+12 & \multicolumn{1}{l|}{6.29E+05} & 2.28E+13 & 1.90E+06 \\ \hline
\end{tabular}
}
\end{table}

Test case~1F serves to assess the sensitivity of the condition numbers to variations in the liquid and ambient viscosities. We consider liquid and ambient viscosities ranging 
between $10^{-5}\,\mathrm{kg}\,\mathrm{m}^{2-d}\,\mathrm{s}^{-1}$ and $10^{-1}\,\mathrm{kg}\,\mathrm{m}^{2-d}\,\mathrm{s}^{-1}$. The lower limit resembles the viscosity of air at atmospheric conditions, while the upper limit resembles the viscosity of low-viscosity oil. Other parameters are as listed in Table~\TAB{cases}. We mention that the convergence of the Newton procedure with $\varepsilon$\nobreakdash-continuation (see Section~\SEC{eps_cont}) is robust for the considered range of viscosities, but the robustness deteriorates for viscosity values that are significantly lower than the lower bound $10^{-5}\,\mathrm{kg}\,\mathrm{m}^{2-d}\,\mathrm{s}^{-1}$. In addition, the magnitude of discretization errors in the form of spurious currents near the interface, increases as the viscosity diminishes. The condition number of the aggregated NSCH system does not display a significant dependence on the viscosities, neither for the original system nor for the equilibrated system.
The condition number of the CH subsystem is independent of the viscosities, in accordance with the fact that~$\eta_{\LL}$ and~$\eta_{\AA}$ do not appear in the CH subsystem. The condition number of the NS subsystem depends weakly on the viscosities for $\eta_{\LL}$ and~$\eta_{\AA}$ between
$10^{-5}\,\mathrm{kg}\,\mathrm{m}^{2-d}\,\mathrm{s}^{-1}$ and $10^{-2}\,\mathrm{kg}\,\mathrm{m}^{2-d}\,\mathrm{s}^{-1}$. However, for larger values of~$\eta_{\LL}$ and~$\eta_{\AA}$, the condition number displays a sharp increase. Equilibration renders the condition number of the NS subsystem essentially independent of~$\eta_{\LL}$ and~$\eta_{\AA}$.

\begin{table}[!htbp]
\caption{Lower bounds of the 1-norm condition numbers as computed by MATLAB's \texttt{condest} function for test case \#1F in Table~\TAB{cases}. Parameters of interest are liquid and ambient viscosities $\eta_{\LL}$ and $\eta_{\AA}$.}
\label{tab:condition_4}
\resizebox{\columnwidth}{!}{
\begin{tabular}{ll|cccccccccc|}
\cline{3-12}
 & & \multicolumn{2}{c|}{$\eta_{\AA}=10^{-5}\,\frac{\text{kg}}{\text{m}^{d-2}\,\text{s}}$}  & \multicolumn{2}{c|}{$\eta_{\AA}=10^{-4}\,\frac{\text{kg}}{\text{m}^{d-2}\,\text{s}}$}  & \multicolumn{2}{c|}{$\eta_{\AA}=10^{-3}\,\frac{\text{kg}}{\text{m}^{d-2}\,\text{s}}$}  & \multicolumn{2}{c|}{$\eta_{\AA}=10^{-2}\,\frac{\text{kg}}{\text{m}^{d-2}\,\text{s}}$}  & \multicolumn{2}{c|}{$\eta_{\AA}=10^{-1}\,\frac{\text{kg}}{\text{m}^{d-2}\,\text{s}}$} \\
 & & orig & \multicolumn{1}{c|}{equi}                                      & orig & \multicolumn{1}{c|}{equi}                                      & orig & \multicolumn{1}{c|}{equi}                                      & orig & \multicolumn{1}{c|}{equi}                                      & orig & \multicolumn{1}{c|}{equi}                                 \\ \hline
\multicolumn{1}{|l}{}                                                                 & full & 4.64E+13 & \multicolumn{1}{l|}{2.53E+07} & 4.93E+13 & \multicolumn{1}{l|}{3.57E+07} & 6.59E+13 & \multicolumn{1}{l|}{6.81E+07} & 6.75E+13 & \multicolumn{1}{l|}{1.15E+08} & 2.04E+14 & 1.57E+08 \\
\multicolumn{1}{|c}{$\eta_{\LL}=10^{-5}\,\frac{\text{kg}}{\text{m}^{d-2}\,\text{s}}$} & CH   & 2.31E+09 & \multicolumn{1}{l|}{4.13E+06} & 2.27E+09 & \multicolumn{1}{l|}{4.12E+06} & 2.27E+09 & \multicolumn{1}{l|}{4.12E+06} & 2.28E+09 & \multicolumn{1}{l|}{4.11E+06} & 2.26E+09 & 4.11E+06 \\
\multicolumn{1}{|l}{}                                                                 & NS   & 2.21E+11 & \multicolumn{1}{l|}{3.91E+06} & 2.21E+11 & \multicolumn{1}{l|}{2.95E+06} & 2.20E+11 & \multicolumn{1}{l|}{2.33E+06} & 2.72E+12 & \multicolumn{1}{l|}{2.48E+06} & 2.04E+14 & 1.71E+07 \\ \hline
\multicolumn{1}{|l}{}                                                                 & full & 4.82E+13 & \multicolumn{1}{l|}{3.36E+07} & 4.81E+13 & \multicolumn{1}{l|}{3.81E+07} & 4.68E+13 & \multicolumn{1}{l|}{9.04E+07} & 6.39E+13 & \multicolumn{1}{l|}{1.56E+08} & 2.11E+14 & 2.08E+08 \\
\multicolumn{1}{|c}{$\eta_{\LL}=10^{-4}\,\frac{\text{kg}}{\text{m}^{d-2}\,\text{s}}$} & CH   & 2.26E+09 & \multicolumn{1}{l|}{4.10E+06} & 2.27E+09 & \multicolumn{1}{l|}{4.14E+06} & 2.27E+09 & \multicolumn{1}{l|}{4.12E+06} & 2.27E+09 & \multicolumn{1}{l|}{4.11E+06} & 2.28E+09 & 4.06E+06 \\
\multicolumn{1}{|l}{}                                                                 & NS   & 2.21E+11 & \multicolumn{1}{l|}{3.69E+06} & 2.28E+11 & \multicolumn{1}{l|}{1.90E+06} & 2.27E+11 & \multicolumn{1}{l|}{1.21E+06} & 3.10E+12 & \multicolumn{1}{l|}{1.15E+06} & 2.11E+14 & 8.25E+06 \\ \hline
\multicolumn{1}{|l}{}                                                                 & full & 5.18E+13 & \multicolumn{1}{l|}{6.89E+07} & 4.83E+13 & \multicolumn{1}{l|}{9.43E+07} & \cellcolor{yellow!10}4.99E+13 & \multicolumn{1}{l|}{\cellcolor{yellow!10}1.55E+08} & 4.82E+13 & \multicolumn{1}{l|}{2.40E+08} & 2.44E+14 & 2.80E+08 \\
\multicolumn{1}{|c}{$\eta_{\LL}=10^{-3}\,\frac{\text{kg}}{\text{m}^{d-2}\,\text{s}}$} & CH   & 2.28E+09 & \multicolumn{1}{l|}{4.13E+06} & 2.24E+09 & \multicolumn{1}{l|}{4.13E+06} & \cellcolor{yellow!10}2.30E+09 & \multicolumn{1}{l|}{\cellcolor{yellow!10}4.19E+06} & 2.31E+09 & \multicolumn{1}{l|}{4.13E+06} & 2.30E+09 & 4.10E+06 \\
\multicolumn{1}{|l}{}                                                                 & NS   & 2.57E+11 & \multicolumn{1}{l|}{3.14E+06} & 2.62E+11 & \multicolumn{1}{l|}{1.60E+06} & \cellcolor{yellow!10}2.81E+11 & \multicolumn{1}{l|}{\cellcolor{yellow!10}5.22E+05} & 3.17E+12 & \multicolumn{1}{l|}{5.34E+05} & 2.44E+14 & 3.25E+06 \\ \hline
\multicolumn{1}{|l}{}                                                                 & full & 4.98E+13 & \multicolumn{1}{l|}{1.57E+08} & 4.84E+13 & \multicolumn{1}{l|}{2.06E+08} & 4.91E+13 & \multicolumn{1}{l|}{3.49E+08} & 4.78E+13 & \multicolumn{1}{l|}{3.24E+08} & 2.47E+14 & 3.28E+08 \\
\multicolumn{1}{|c}{$\eta_{\LL}=10^{-2}\,\frac{\text{kg}}{\text{m}^{d-2}\,\text{s}}$} & CH   & 2.28E+09 & \multicolumn{1}{l|}{4.16E+06} & 2.29E+09 & \multicolumn{1}{l|}{4.16E+06} & 2.34E+09 & \multicolumn{1}{l|}{4.20E+06} & 2.30E+09 & \multicolumn{1}{l|}{4.21E+06} & 2.31E+09 & 4.13E+06 \\
\multicolumn{1}{|l}{}                                                                 & NS   & 2.03E+12 & \multicolumn{1}{l|}{4.47E+06} & 2.00E+12 & \multicolumn{1}{l|}{2.06E+06} & 2.33E+12 & \multicolumn{1}{l|}{6.57E+05} & 3.51E+12 & \multicolumn{1}{l|}{4.85E+05} & 2.47E+14 & 2.35E+06 \\ \hline
\multicolumn{1}{|l}{}                                                                 & full & 1.50E+14 & \multicolumn{1}{l|}{9.12E+07} & 1.51E+14 & \multicolumn{1}{l|}{1.51E+08} & 1.52E+14 & \multicolumn{1}{l|}{2.89E+08} & 1.82E+14 & \multicolumn{1}{l|}{5.39E+08} & 2.93E+14 & 3.63E+08 \\
\multicolumn{1}{|c}{$\eta_{\LL}=10^{-1}\,\frac{\text{kg}}{\text{m}^{d-2}\,\text{s}}$} & CH   & 2.29E+09 & \multicolumn{1}{l|}{4.16E+06} & 2.32E+09 & \multicolumn{1}{l|}{4.16E+06} & 2.29E+09 & \multicolumn{1}{l|}{4.17E+06} & 2.40E+09 & \multicolumn{1}{l|}{4.31E+06} & 2.30E+09 & 4.20E+06 \\
\multicolumn{1}{|l}{}                                                                 & NS   & 1.50E+14 & \multicolumn{1}{l|}{1.23E+07} & 1.51E+14 & \multicolumn{1}{l|}{7.06E+06} & 1.52E+14 & \multicolumn{1}{l|}{3.62E+06} & 1.82E+14 & \multicolumn{1}{l|}{4.58E+06} & 2.91E+14 & 1.90E+06 \\ \hline
\end{tabular}
}
\end{table}

\begin{figure}[h!!!!!!!!!!]
\centering
\includegraphics[width=\textwidth, trim = 0mm 1mm 10mm 9mm, clip]{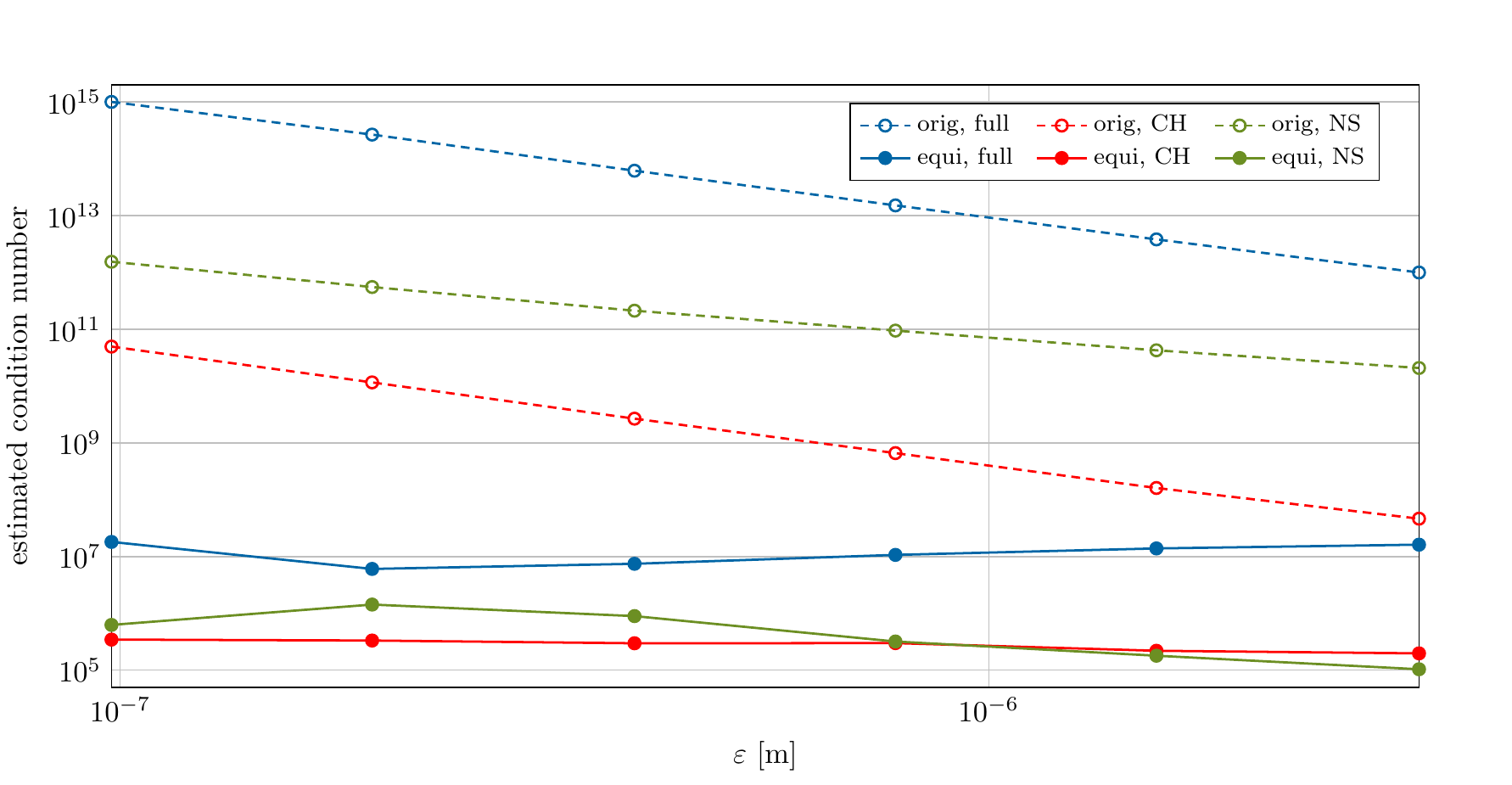}
\caption{Lower bounds of the 1-norm condition numbers as computed by MATLAB's \texttt{condest} function for test case 1G in Table~\TAB{cases}. The interface thickness $\varepsilon \coloneqq \varepsilon_k$, mobility coefficient $m \coloneqq m_k$ and the maximum refinement level $L_\MAX \coloneqq L_{\MAX,k}$ are adapted simultaneously according to~\EQ{simultaneous}.}
\label{fig:eps_mob_L}
\end{figure}

Finally, we consider in test case~1G the effect of simultaneous variation of $\varepsilon,m$ and $L_{\MAX}$, corresponding to a scenario in which the interface-thickness parameter $\varepsilon$ is decreased to approximate the sharp-interface limit more closely, the mobility parameter $m$ is reduced proportional to $\varepsilon^3$ to fix the diffusive time scale $T_{\DIFF}$ (see~\APP{time_scale}), and the maximum refinement level 
$L_{\MAX}$ is raised to ensure adequate resolution of the diffuse interface at length scale~$\varepsilon$. In particular, we regard
\begin{equation}
\label{eq:simultaneous}
\varepsilon_k=2^{-k}\!\times\!100\,\text{\textmu{}m},
\quad
m_k=\frac{\varepsilon_k^3}{\sigma\,T_{\DIFF}},
\qquad
L_{\MAX,k}=k-5,
\qquad
k=5,6,\ldots,10
\end{equation}
with $T_{\DIFF}=2.686\!\times\!{}10^{-8}\,\textrm{s}$. Other parameters are as listed in Table~\TAB{cases}.
Figure~\FIG{eps_mob_L} plots the estimated condition number of the tangent matrix of the aggregated NSCH system and of the NS and CH subsystems versus the interface thickness, both for the original systems and their equilibrated counter parts. The condition numbers in Figure~\FIG{eps_mob_L} pertain to the super-refined approximations. The largest value of~$\varepsilon$ in  Figure~\FIG{eps_mob_L} correspond to $L_{\MAX}=0$. In this case, the super-refinement consists of a uniform mesh with mesh width 
$h=h_0/2=\frac{5}{2}\,\text{\textmu{}m}$. Figure~\FIG{eps_mob_L} conveys that the condition number of the aggregated NSCH system increases approximately as $\varepsilon^{-2}$ as~$\varepsilon$ is reduced and~$m$ and~$L_{\MAX}$ are adapted accordingly. The condition number of the CH subsystem scales similarly with decreasing~$\varepsilon$. The condition number of the NS subsystem is approximately proportional to~$\varepsilon^{-3/2}$. Equilibration is effective, and renders the condition numbers of the NSCH system and the CH subsystem essentially independent of~$\varepsilon$. 
The condition number of the equilibrated NS subsystem increases approximately proportional to~$\varepsilon^{-1}$ as $\varepsilon$ decreases, and~$m$ and $L_{\MAX}$ are adapted simultaneously. It is noteworthy that the results for the equilibrated condition numbers in Figure~\FIG{eps_mob_L} are consistent with the results in Table~\TAB{condition_2} and Figure~\FIG{results_m}, in the sense that according 
to Table~\TAB{condition_2} the condition number of the equilibrated NSCH, NS and CH systems scale as $\varepsilon^{-3},\varepsilon^{-1}$ and~$\varepsilon^{-3}$, respectively, under simultaneous~$\varepsilon$ and $L_{\MAX}$ refinement, while according to Figure~\FIG{results_m}
the condition numbers of the equilibrated NSCH, NS and CH systems scale as $m^1,m^0$ and~$m^1$, respectively, as~$m$ decreases.

\subsection{preconditioned-GMRES convergence}
\label{sec:GMRES_conv}
To analyze the convergence behavior of GMRES with the $\boldsymbol{P}_{\UPP}$ and $\boldsymbol{P}_{\LOW}$ preconditioners introduced in Section~\SEC{partitioned}, we regard the residual reduction versus the number of GMRES iterations for the oscillating-droplet test case with parameter values as listed in Table~\TAB{cases}, test case~2. Specifically, the parameters for the liquid and ambient fluid are selected such that the test case pertains to the oscillation of a water droplet suspended in air. We consider a relatively thin diffuse interface, corresponding to an interface-thickness parameter~$\varepsilon=2^{-10}\!\times\!100\,\text{\textmu{}m}$. This implies that $\varepsilon$ is approximately $150\times$ smaller than the radius of the circular droplet with the same ($2$-dimensional) volume as the ellipse in~\EQ{quartellipse}.
The maximum refinement level, $L_{\MAX}=7$, is chosen such that the ratio of the interface thickness parameter, $\varepsilon$, to the minimal mesh width, $h$, is $\varepsilon/h=2^{-2}\!\times\!10$, and the interface is adequately resolved. 

\begin{figure}[h!!!!!!!!!!]
\centering
\includegraphics[width=\textwidth, trim = 0mm 1mm 10mm 9mm, clip]{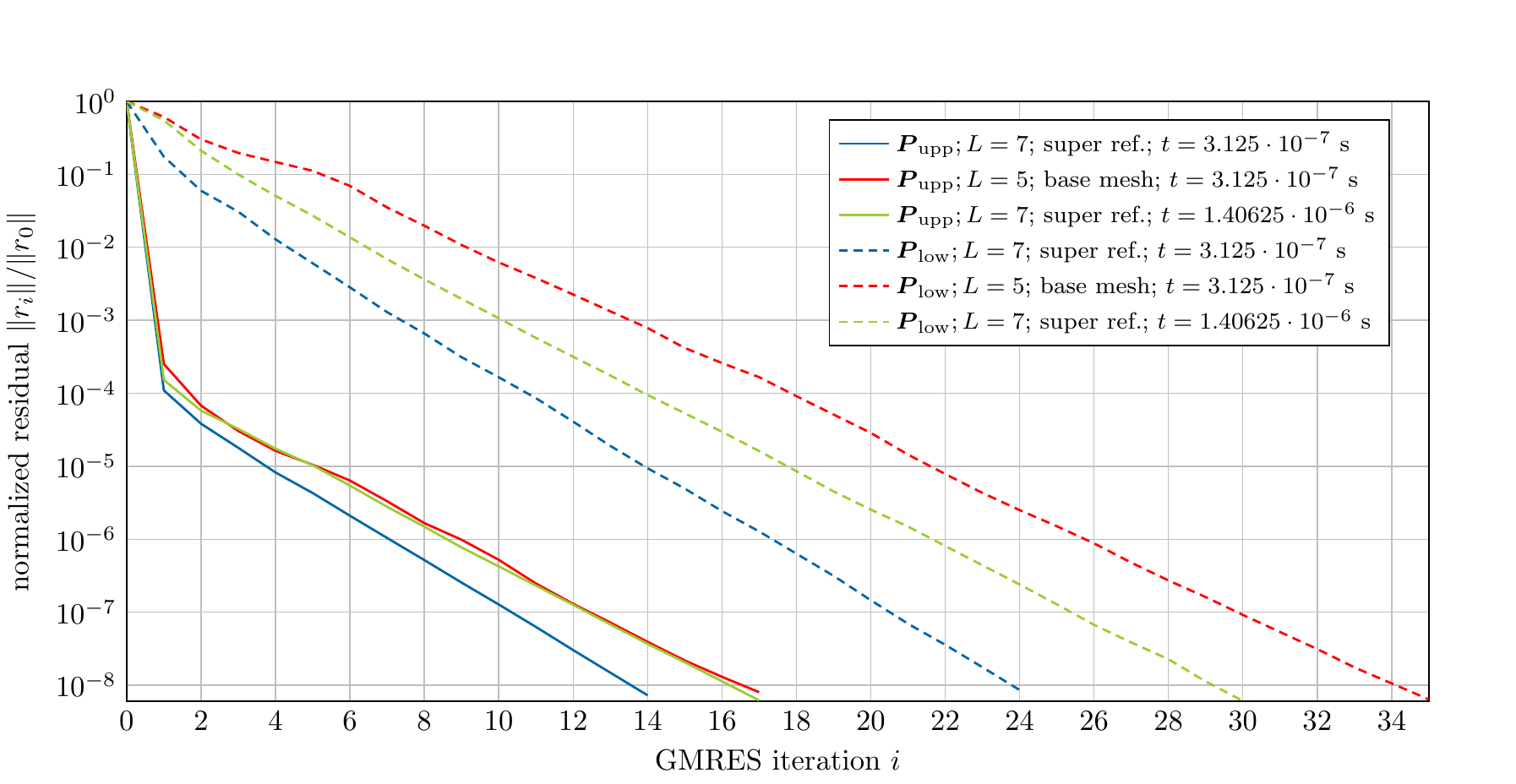}
\caption{Convergence behavior of $\boldsymbol{P}_{\UPP}$ and $\boldsymbol{P}_{\LOW}$ preconditioned GMRES applied to the linear tangent problem within the final Newton iteration of the super-refinement of refinement level $L = L_{\MAX} = 7$ of the first successful time step of test case \#2 in Table \TAB{cases}.}
\label{fig:GMRES}
\end{figure}

The convergence behavior of the preconditioned GMRES method displays some minor variability during the Newton iterations, during the refinement iterations, and in time. To illustrate the main aspects of the convergence behavior, Figure~\FIG{GMRES} plots
the residual reduction $\|r_i\|/\|r_0\|$, versus the number of GMRES iterations, $i$, for the $\boldsymbol{P}_{\UPP}$ and $\boldsymbol{P}_{\LOW}$ preconditioners, for the super-refinement of the maximum refinement level in the first time step and in time step~$25$, and for refinement level $L=5$ in the first time step. For completeness, we mention that the results in Figure~\FIG{GMRES} pertain to the final step in the Newton procedure, but similar behavior is observed in other Newton iterations. One can observe that the upper-triangular preconditioner yields much more effective convergence behavior than the lower-triangular preconditioner (see also Remark~\ref{rem:PuppPlow}), on account of the fact that $\boldsymbol{P}_{\UPP}$\nobreakdash-preconditioned GMRES exhibits a significant decrease in the residual in the first iteration. In subsequent iterations, the residual reduction per iteration is very similar for $\boldsymbol{P}_{\UPP}$\nobreakdash-preconditioned GMRES and $\boldsymbol{P}_{\LOW}$\nobreakdash-preconditioned GMRES, and both are essentially constant within each tangent problem, resulting in almost straight curves 
in Figure~\FIG{GMRES} for $i>1$. Figure~\FIG{GMRES} conveys that, excepting the first iteration, the residual reduction per GMRES iteration lies between~$0.5$ and~$0.6$. 

\begin{remark}
\label{rem:GMRES_tolerance}
As a stopping criterion for the preconditioned GMRES procedure, we require that the residual reduction exceeds~$10^8$. Such a stringent condition on the approximate solution of the linear tangent problems in~\EQ{Newton}$_1$ is necessary to ensure proper (quadratic) convergence of the Newton process~\EQ{Newton}. If the stopping criterion for the GMRES procedure is too loose, the Newton method generally displays suboptimal convergence, or convergence can stagnate.
\end{remark}

\begin{remark}
In the first time step of test case~2, time-step reduction (by a factor of 4) is invoked, to ensure adequate convergence of the Newton procedure;
see Section~\SEC{partitioned_monolithic}. We conjecture that the poor convergence of the Newton method in the first time step with the original time step according to Table~\TAB{cases} is caused by the fact that the $\tanh{}$ profile of the initial data for the order parameter in~\EQ{varphi0} is not equilibrated for the elliptic initial configuration~\EQ{quartellipse}, which can lead to strong transient behavior. The results reported in Figure~\FIG{GMRES} for the first time step in fact pertain to the first time step which is completed at the original time step~$\tau_{\MAX}$ in Tabel~\TAB{cases}.
\end{remark}

\subsection{Comparison to sharp-interface model}
\label{sec:num_thin_int}
A fundamental notion underlying diffuse-interface models for binary-fluid flows, is that in the sharp-interface limit $\varepsilon \! \to{} \! +0$ the solution approaches the solution of a corresponding sharp-interface model (provided that the limit exists), viz.\ a model comprising the incompressible Navier--Stokes equations on two complementary domains, separated by a manifold of codimension~1 which carries surface energy; see~\cite{Abels:2013bd,Abels:2018ly}. In actual computations of binary-fluid flows based on diffuse-interface models, the sharp-interface limit is inaccessible, i.e.\ one cannot pass to the limit $\varepsilon \! \to{} \! +0$, and one generally assumes that  the diffuse-interface model forms an accurate representation of the sharp-interface model if the interface-thickness parameter~$\varepsilon$ is sufficiently small relative to other relevant length scales in the problem under consideration, e.g.\ the minimal radius of curvature of the fluid--fluid meniscus. It is noteworthy that analyses of sharp-interface limits of diffuse-interface models via matched-asymptotic expansions (see e.g.~\cite{Abels:2018ly}) do not generally provide insight into the deviation between the two models, as the deviation pertains to higher-order terms.

In this section we compare results of the AGG NSCH diffuse-interface model with a relatively thin diffuse interface to corresponding sharp-interface results for an oscillating-droplet test case. The characteristic parameters are listed in Table~\TAB{cases} under test case~2. Note that this setup is the same as used for the investigation of the convergence properties of preconditioned GMRES in Section~\SEC{GMRES_conv}. For the ellipsoidal initial configuration of the droplet according to~\EQ{quartellipse}, the minimal radius of curvature occurs on the major axis of the ellipse and is~$5\,\text{\textmu{}m}$. The interface-thickness parameter~$\varepsilon$ is therefore approximately $70\times$ smaller than this minimal radius of curvature, suggesting that the sharp-interface limit is adequately approximated. Let us note that the disparity in length scales is enabled by the adaptive refinement procedure described in Section~\SEC{error_adaptive}. Rayleigh's classical theory~\cite{Strutt:1879xi} on the oscillation of inviscid droplets conveys that the period of oscillation of the considered droplet is approximately:
\begin{equation}
\label{eq:Rayleigh}
T_{\textrm{osc}}=\frac{2\pi}{\omega_n}\qquad{}\omega_n^2=(n^3-n)\frac{\sigma_{\LA}}{\rho_{\LL}R^3}
\end{equation}
with $n=2$ for the considered ellipsoidal mode of oscillation and $R=\sqrt{200}\,\text{\textmu{}m}$ the radius of the circular droplet with the same ($2$\nobreakdash-dimensional) volume as the ellipse in~\EQ{quartellipse}. Equation~\EQ{Rayleigh} implies $T_{\textrm{osc}}\approx{}16\,\text{\textmu{}s}$. The period of oscillation, $T_{\textrm{osc}}$, is approximately $600\times$ larger than the diffusive time scale $T_{\DIFF}=\varepsilon^3/\sigma_{\LA}m$. 

We compare the results of the AGG NSCH model with those of a sharp-interface model based on the methodology described in~\cite{Hack:2021rf}. In the sharp-interface model, the liquid and ambient fluids are represented by two systems of incompressible Navier--Stokes equations on complementary domains, $\Omega_{\LL}$ and $\Omega_{\AA}=\Omega\setminus\Omega_{\LL}$, respectively. The two NS systems are coupled at their mutual interface by dynamic and kinematic boundary conditions. The motion of the liquid and ambient domains is accommodated by means of an Arbitrary-Lagrangian-Eulerian formulation. The NS equations in ALE form are approximated by means of a finite-element method, using standard Taylor--Hood P2\nobreakdash-P1 elements for velocity and pressure. The sharp-interface model employs a time step $\tau_{\textsc{si}} = 2^{-8}\!\times\!10\,\text{\textmu{}s}$.

\begin{figure}[h!!!!!!!!!!]
\centering
\includegraphics[width=\textwidth, trim = 0mm 1mm 10mm 9mm, clip]{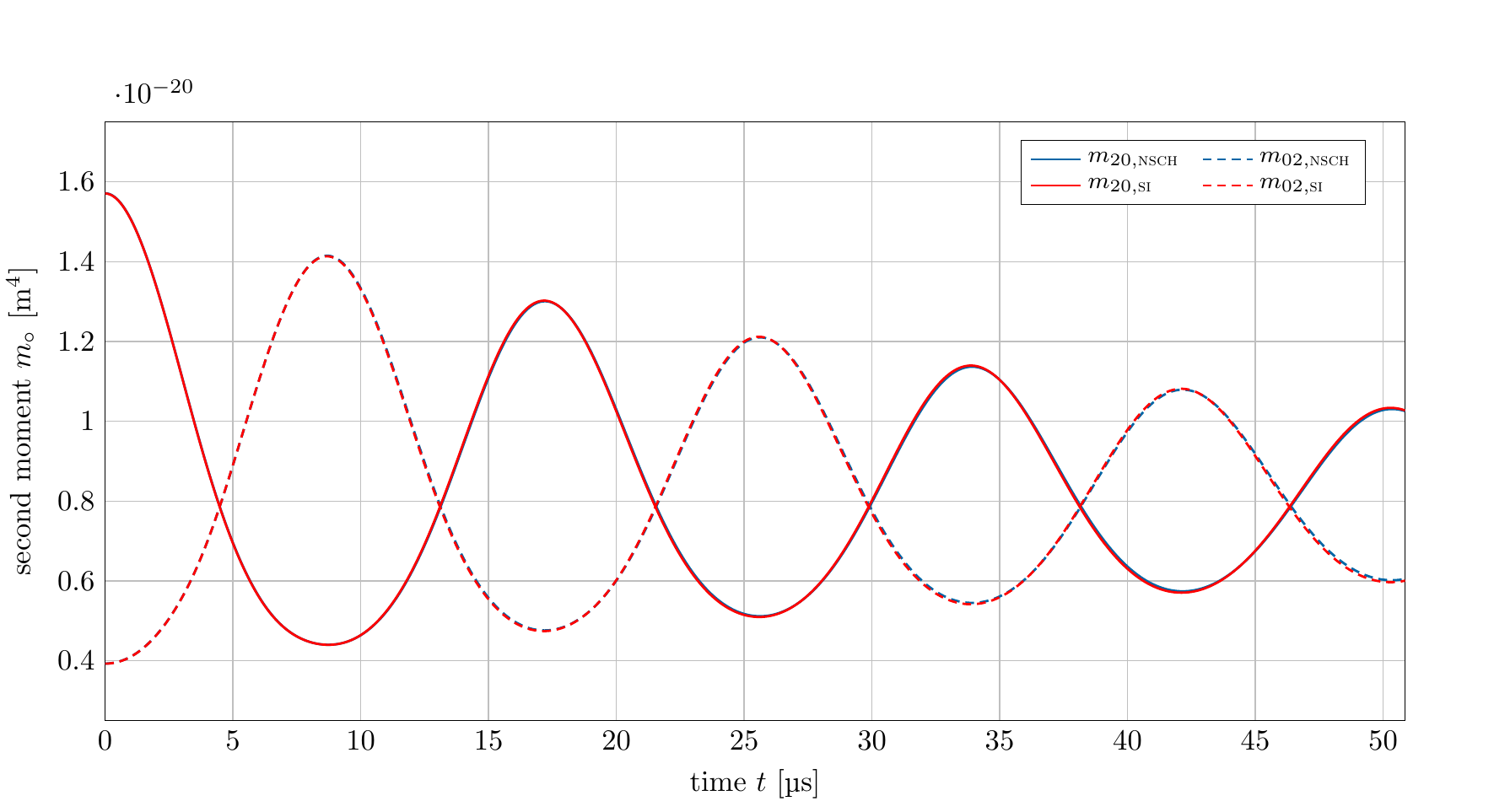}
\caption{Second area moments $m_{20}$ (continuous) and $m_{02}$ (dashed) versus time $t$ for the NSCH (blue) and sharp-interface (red) models.}
\label{fig:moments}
\end{figure}

Intrinsically, diffuse-interface models do not provide an explicit representation of the interface. A comparison of the diffuse-interface results to the sharp-interface results in terms of the location of the interface is therefore non-trivial. We present a comparison in terms of geometric quantities that are well defined for both the diffuse-interface model and the sharp-interface model, viz.\ the second area moments~$m_{20}$ and~$m_{02}$ with respect to the axes  of symmetry:
\begin{alignat}{2}
m_{20,\textsc{nsch}}(t) &= \int_\Omega \frac{1}{2}(1+\varphi(t,\boldsymbol{x}))\, x_1^2 \, \textrm{d} \boldsymbol{x}, 
&\qquad
m_{20,\textsc{si}}(t) &= \int_{\Omega_{\LL}(t)} x_1^2 \, \textrm{d} \boldsymbol{x}, 
\\
m_{02,\textsc{nsch}}(t) &= \int_\Omega \frac{1}{2}(1+\varphi(t,\boldsymbol{x}))\, x_2^2 \, \textrm{d} \boldsymbol{x},
&\qquad
m_{02,\textsc{si}}(t) &= \int_{\Omega_{\LL}(t)} x_2^2 \, \textrm{d} \boldsymbol x.
\end{alignat}
Note that these values differ by a factor 4 from those of a full droplet, because only a quarter of the droplet is considered.
Figure~\FIG{moments} plots the second area moments of the NSCH and SI models versus time. The results in Figure~\FIG{moments} illustrate that the motion of the liquid-ambient interface, as expressed by the second area moments, is essentially indistinguishable for the diffuse-interface and sharp-interface models. The difference in oscillation period between both models, measured after two oscillations, is 
approximately $2^{-8}\!\times\!10\,\text{\textmu{}s}$, which is smaller than the time-step sizes used in both the NSCH and 
the sharp-interface simulations. 

Both the diffuse-interface and the sharp-interface binary-fluid model satisfy an energy-dissipation inequality. For the sharp-interface model,
the energy-dissipation inequality is:
\begin{equation}
\label{eq:sidissipation}
\frac{\mathrm{d}}{\mathrm{d}t}\big(E_{\textrm{kin},\textsc{si}}(t)+E_{\textrm{int},\textsc{si}}(t)\big)
\leq
-\int_{\Omega_{\LL}(t)}\nabla \boldsymbol u:\boldsymbol{\tau}_{\LL}\,\mathrm{d}\boldsymbol{x}
-\int_{\Omega_{\AA}(t)}\nabla \boldsymbol u:\boldsymbol{\tau}_{\AA}\,\mathrm{d}\boldsymbol{x}
+\textrm{bnd}
\end{equation}
where `$\textrm{bnd}$' denotes terms that are only supported on the external boundary~$\Gamma_{\EXT}$ and
\begin{equation}
E_{\textrm{kin},\textsc{si}}(t)=\int_{\Omega_{\LL}(t)}\frac{1}{2}\rho_{\LL}|\boldsymbol{u}|^2\,\mathrm{d}\boldsymbol{x}
+\int_{\Omega_{\AA}(t)}\frac{1}{2}\rho_{\AA}|\boldsymbol{u}|^2\,\mathrm{d}\boldsymbol{x},
\qquad
E_{\textrm{int},\textsc{si}}(t)
=\sigma_{\LA}\operatorname{meas}(\Gamma_{\LA}(t))
\end{equation}
represent kinetic energy and interface energy, respectively. The viscous-stress tensors~$\boldsymbol{\tau}_{\LL}$ and~$\boldsymbol{\tau}_{\AA}$ in~\EQ{sidissipation} are identical in form to~\EQ{tau} with the viscosity parameter according to~$\eta_{\LL}$ and~$\eta_{\AA}$, respectively.
One can infer that the integrals appearing on the right-hand side of~\EQ{sidissipation} are non-negative, from which it follows 
that in the sharp-interface model, the energy decays by viscous dissipation, up to boundary terms. The energy-dissipation inequality corresponding to the NSCH model is 
\begin{equation}
\label{eq:nschdissipation}
\frac{\mathrm{d}}{\mathrm{d}t}\big(E_{\textrm{kin},\textsc{nsch}}(t)+E_{\textrm{int},\textsc{nsch}}(t)\big)
\leq
-\int_{\Omega}\nabla \boldsymbol u:\boldsymbol{\tau}\,\mathrm{d}\boldsymbol{x}
-\int_{\Omega}m|\nabla\mu|^2\,\mathrm{d}\boldsymbol{x}
+\textrm{bnd}
\end{equation}
with
\begin{equation}
\label{eq:kinintnsch}
E_{\textrm{kin},\textsc{nsch}}(t)=\int_{\Omega}\frac{1}{2}\rho|\boldsymbol{u}|^2\,\mathrm{d}\boldsymbol{x},
\qquad
E_{\textrm{int},\textsc{nsch}}(t) = \int_\Omega \left( \frac{ \varepsilon\sigma  }{ 2 }  | \nabla \varphi |^2+\frac{ \sigma }{ \varepsilon } \Psi \right) \textrm{d} \boldsymbol x,
\end{equation}
where $\rho$ is the mixture density according to~\EQ{rhophi}. In the diffuse-interface model, the interface energy is represented by a Ginzburg--Landau-type functional, corresponding to the second functional in~\EQ{kinintnsch}. Equation~\EQ{nschdissipation} conveys that in the diffuse-interface model, the energy decays by viscous dissipation and by gradients in the chemical potential, up to boundary terms. 

\begin{figure}[h!!!!!!!!!!]
\centering
\includegraphics[width=\textwidth, trim = 0mm 1mm 10mm 9mm, clip]{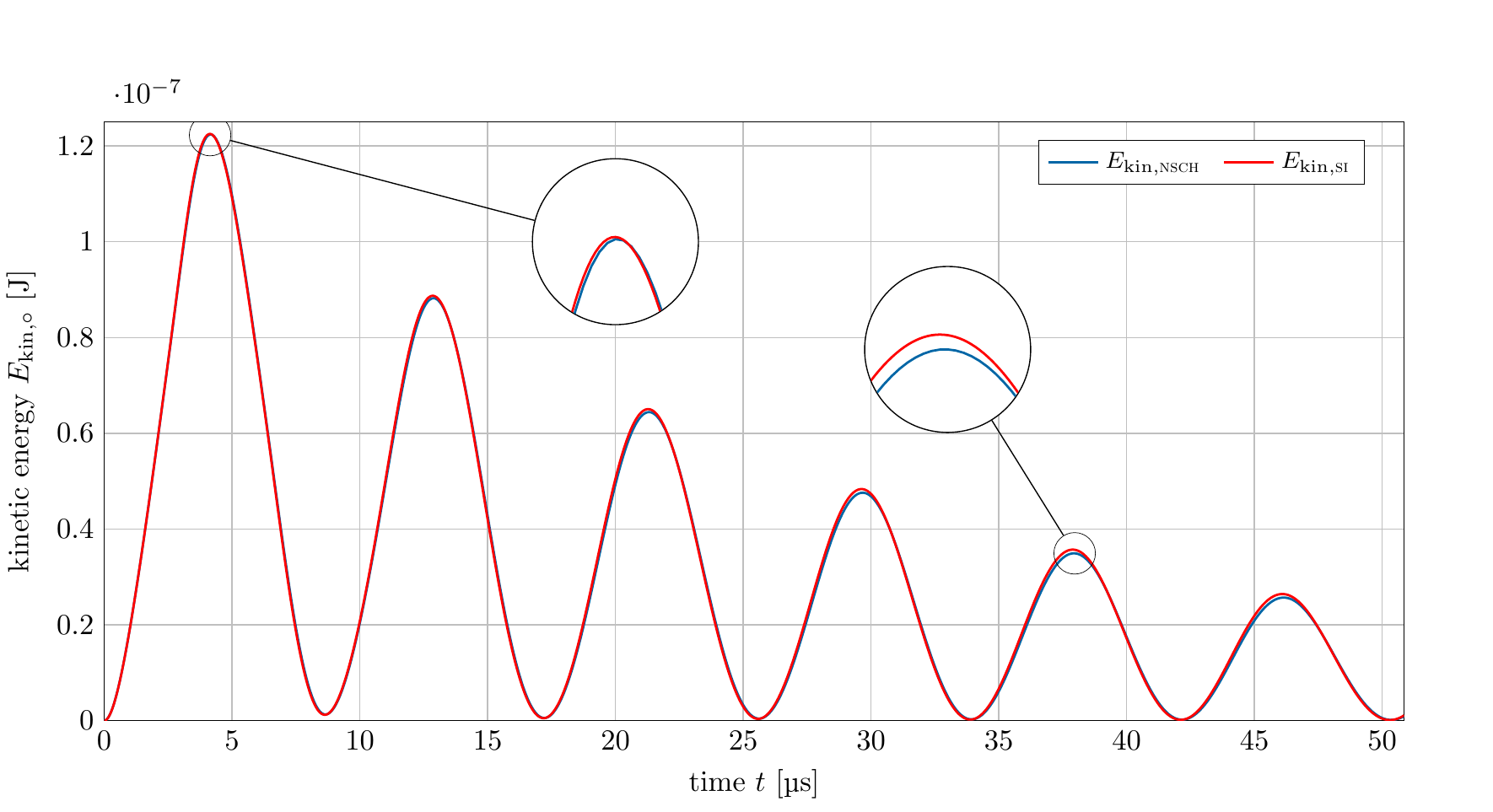}
\caption{The kinetic energies $E_{\textrm{kin},\circ}$ [J] over time $t$ [\textmu s] for the NSCH (blue) and sharp-interface (red) models. A magnification on the local extrema after both $0.25$ and $2.25$ oscillations are provided.}
\label{fig:kin_energy}
\end{figure}

\begin{figure}[h!!!!!!!!!!]
\centering
\includegraphics[width=\textwidth, trim = 0mm 1mm 10mm 9mm, clip]{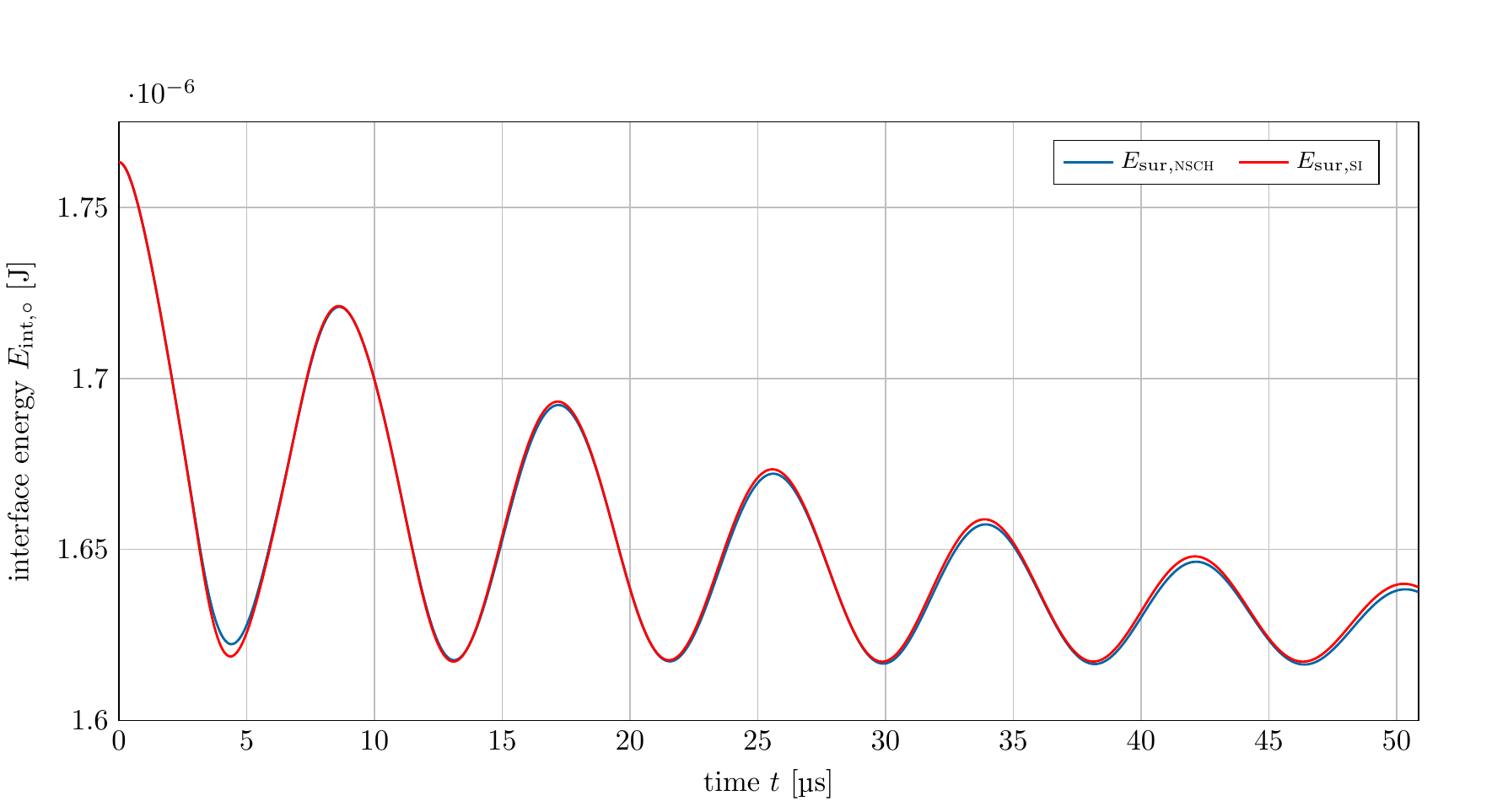}
\caption{The Ginzburg--Landau energy $E_{\textrm{int},\textsc{nsch}}$ and interface energy $E_{\textrm{int},\textsc{si}}$ [J] over time $t$ [\textmu s] for the NSCH (blue) and sharp-interface (red) models, respectively.}
\label{fig:sur_energy}
\end{figure}

\begin{figure}[h!!!!!!!!!!]
\centering
\includegraphics[width=\textwidth, trim = 0mm 1mm 10mm 9mm, clip]{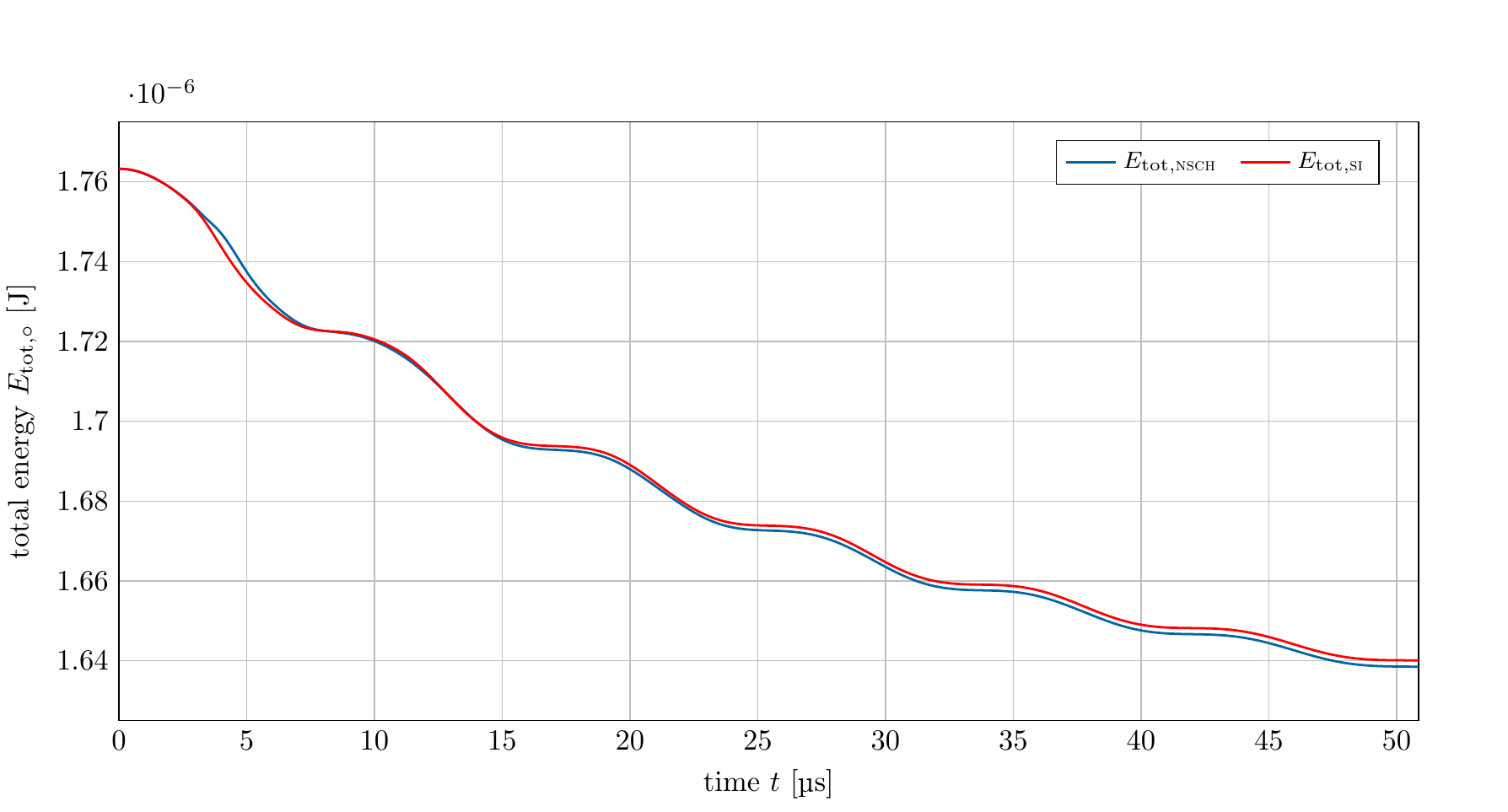}
\caption{The total energies $E_{\textrm{tot},\circ}$ [J] over time $t$ [\textmu s] for the NSCH (blue) and sharp-interface (red) models.}
\label{fig:total_energy}
\end{figure}

To compare the energy-dissipation properties of the AGG NSCH model and the sharp-interface model, Figures~\FIG{kin_energy} and~\FIG{sur_energy} plot the evolution of the kinetic energy and the interface energy pertaining to both models, respectively. In addition, Figure~\FIG{total_energy} displays the evolution of the total energy for both models. From Figures~\FIG{moments} and~\FIG{kin_energy}, one can infer that the kinetic energy attains minima at the extrema of the second area moments, i.e.\ when the deformation of the droplet attains its (local in time) maximum. Figure~\FIG{sur_energy} conveys that the interface energy is maximal at the minima of the kinetic energy, which illustrates the exchange of kinetic and interface energy during the oscillation of the droplet, up to dissipation. The minimal  interface energy is approximately $1.617\,\text{\textmu{}J}$, which coincides with (one quarter of) the surface energy of a spherical droplet with the same volume as the ellipse in~\EQ{quartellipse}. From Figs.~\FIG{kin_energy} and~\FIG{total_energy}, one can observe that the energy dissipation is largest  
near the maxima of the kinetic energy, coinciding with the spherical configuration of the droplet. At these instants, the velocities in 
the liquid droplet are largest, inducing the largest viscous dissipation. Overall, the total energy as well as both its components exhibit very similar behavior in the NSCH model and in the sharp-interface model, except that the NSCH model exhibits slightly more dissipation. 

From Figure~\FIG{sur_energy} it can also be observed that the interface energy of the NSCH model deviates more pronouncedly from that of the sharp-interface model near the first minimum at approximately~$4.4\,\text{\textmu{}s}$, and a corresponding deviation occurs in the total energy in Figure~\FIG{total_energy}. This deviation is caused by an initial transient, which is temporally under-resolved in the approximation of the NSCH model at the considered time step.

The AGG NSCH system~\EQ{strong} is thermodynamically consistent in accordance with the free-energy-dissipation relation~\EQ{nschdissipation}. Despite the use of double-well splitting, neither the backward Euler scheme~\EQ{error_system_BE} nor the Crank--Nicolson scheme~\EQ{error_system_CN} do however preserve this property uniformly, i.e.\ independent of the time step. The results in Figure~\FIG{total_energy} convey that the discrete approximation provided by the Crank--Nicolson scheme with time 
step~$\tau_{\MAX}$ according to Table~\TAB{cases} complies with the energy-dissipation relation. It is to be mentioned, however, that the energy-dissipation relation of the Crank-Nicolson scheme is conditional on the time step, and it is violated at a twice larger time step (results not displayed).

\begin{figure}[h!!!!!!!!!!]
\begin{center}
\includegraphics[width=\textwidth]{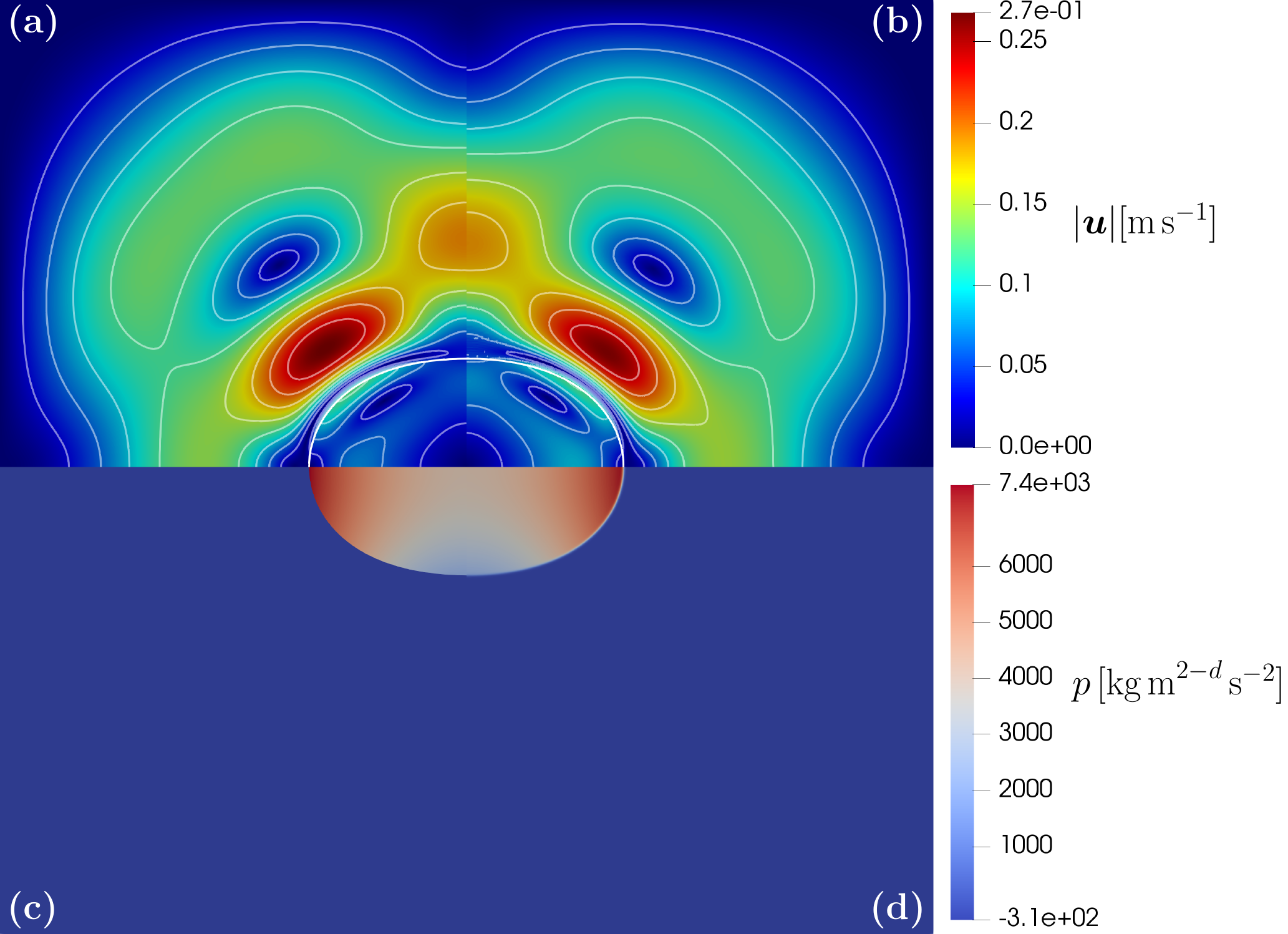}
\end{center}
\caption{Snapshot of the velocity magnitude (top) and pressure (bottom) profiles for the SI (left) and NSCH (right) models, after two oscillations at $t=3.38 \cdot 10^{-5}$ s. The colorbars are identical for both models; the contour lines for the velocity magnitude are divided into 9 intervals over the colorbar range. $\Gamma_\LA$ in (a) and a $\varphi=0$ contour line in (b) have been added to highlight the interface locations.}
\label{fig:qual_vel_pres}
\end{figure}

Finally, we present a qualitative comparison of the the droplet configuration, velocity magnitude~$|\boldsymbol u|$ and pressure field~$p$.
Figure~\FIG{qual_vel_pres} plots the computed~$|\boldsymbol{u}|$ and~$p$ fields at $t\approx{}33.8\,\text{\textmu{}s}$, near the instant at which the second area moments reach their extrema after two full oscillations of the droplet; see Figure~\FIG{moments}. In the NSCH results, the representation of the droplet geometry is indicated by a plot of the~$\varphi=0$ contour line. The droplet configurations pertaining to the NSCH and to the SI model are almost indistinguishable. The pressure fields of the SI and NSCH models are also virtually identical. The velocity-magnitude fields in the SI and NSCH models match closely: the local extrema are in nearly identical locations, both in the droplet and in the ambience. Minor deviations between the NSCH and SI results can be observed, which are accentuated by the contour lines. It is noteworthy that at the considered instant, the velocity in the droplet almost vanishes, while the ambient velocity is relatively large. Because the ambient density is~$10^3$ times smaller than the liquid density, however, the kinetic energy nearly vanishes; cf. Figure~\FIG{kin_energy}.

%
%===============================================================================================================%
%===============================================================================================================%
%===============================================================================================================%

\section{Conclusion}
\label{sec:concl}
In this work, we developed and analyzed an adaptive-approximation method for binary-fluid-flow problems, 
based on the Abels--Garcke--Gr\"un (AGG) Navier--Stokes-Cahn--Hilliard (NSCH) diffuse-interface model. To circumvent instabilities in the formulation for large density and viscosity contrasts due to range-violations of the phase field, we applied an Arrhenius interpolation for the mixture viscosity and a soft-clipped linear interpolation for the mixture density.

To resolve the spatial scale disparity between the diffuse-interface thickness and other representative length scales in the problem, we introduced an adaptive-refinement procedure directed by a two-level hierarchical a-posteriori error estimate. We regarded approximations based on truncated hierarchical B\nobreakdash-splines, restricted to standard $C^0$\nobreakdash-continuous spline bases. The adaptive-approximation framework constructs a sequence of hierarchically refined approximation spaces within each time step by means of the \texttt{Solve} $\rightarrow$ \texttt{Estimate} $\rightarrow$ \texttt{Mark} $\rightarrow$ \texttt{Refine} (SEMR) procedure. 

We introduced several enhancements of the basic SEMR process to improve the robustness and efficiency of the adaptive-refinement method and the corresponding solution procedure:
\begin{itemize}
\item On coarse levels of the adaptive-refinement process, we apply a Backward Euler (BE) scheme with a second order splitting of the double-well potential. On the finest levels, we use a Crank--Nicolson (CN) scheme to recover second-order-in-time accuracy of the approximation. On the coarse levels, the computational complexity of the BE scheme is significantly lower than that of the CN scheme. In addition, the BE scheme is generally more stable on coarse meshes. The temporal accuracy of the BE scheme is too limited to apply it directly as a time integrator for diffuse-interface binary-fluid-flow problems with significant interface dynamics. However, the BE scheme is adequate to guide the initial refinement steps, and to provide an initial estimate for the Newton procedure on the next hierarchical level of refinement.
\item We introduced an $\varepsilon$\nobreakdash-continuation approach, in which both the interface-thickness parameter, $\varepsilon$, and the mobility, $m\coloneqq{}m_{\varepsilon}$, are adapted during the adaptive-refinement process, in such a manner that the diffuse-interface is adequately resolved at each refinement level, for the interface-thickness parameter at that level. The mobility is adapted in accordance with the interface thickness, such that the diffuse interface is suitably equilibrated within a time step. We established that to this purpose, the mobility needs to scale with the (level-dependent) diffuse-interface thickness 
as~$m_{\varepsilon}\propto{}\varepsilon^3$. The $\varepsilon$\nobreakdash-continuation approach serves two purposes. First, it enhances the robustness of the solution procedure on coarse levels of approximation, where the diffuse interface corresponding to the original interface-thickness parameter would be severely under-resolved, leading to non-robustness of the Newton procedure. Second, it mitigates the conditions on the maximum time step in relation to robustness of the Newton procedure.
\item We proposed a partitioned solution method for the linear tangent problems in Newton's method for the coupled NSCH system, to circumvent 
ill-conditioning of the monolithic system. The partitioned solution method corresponds to preconditioned GMRES, wherein the 
preconditioner is the original tangent matrix with either the CH-NS or the NS-CH submatrix suppressed. Suppression of either of these submatrices effectively decouples the NSCH system into its NS and CH constituents, which can then be solved separately. 
We conjectured that the preconditioner in which the NS-CH submatrix is retained ($\boldsymbol{P}_{\UPP}$) is more effective than the preconditioner without the NS-CH submatrix ($\boldsymbol{P}_{\LOW}$), because the former provides an implicit treatment of the surface-tension terms. This conjecture was confirmed by the numerical experiments. 
\end{itemize}

To determine the main properties of the AGG NSCH model, and of the enhanced adaptive approximation method, we conducted numerical experiments for an oscillating droplet in an ambient fluid. We empirically investigated the dependence of the condition number of the linear tangent matrices in Newton's method on the main problem parameters, both for the monolithic NSCH system and for the NS and CH submatrices. To account for the fact that contemporary direct solution methods generally apply equilibration, we also considered the effect of equilibration on the condition numbers of these (sub\nobreakdash-)matrices. We found that equilibration is in general very effective, reducing the conditions numbers in many cases by orders of magnitude. Moreover, equilibration in many cases moderates or eliminates parameter dependence of the condition numbers. The effectiveness of equilibration however typically deteriorates as the system size increases, e.g.\ under mesh refinement. 
In the relevant range of the system-parameter values, we found that the condition numbers of the equilibrated NSCH system and of the equilibrated NS and CH subsystems depend on the minimal element size in the adaptively refined mesh,~$h_{\MIN}$, the interface-thickness parameter,~$\varepsilon$, the mobility,~$m$, the time-step size,~$\tau$, the liquid and ambient densities,~$\rho_{\LL}$ and~$\rho_{\AA}$, and the liquid and ambient viscosities,~$\eta_{\LL}$
and~$\eta_{\AA}$, as:
\begin{equation}
\begin{aligned}
\operatorname{cond}_{\textsc{nsch}}&\propto{}\varepsilon^{\frac{1}{2}}\,m^1\,\rho_{\LL}^0\,\rho_{\AA}^0\,\eta_{\LL}^0\,\eta_{\AA}^0\,h_{\MIN}^{-\frac{7}{2}}\,\tau^1
\\
\operatorname{cond}_{\textsc{ns}}&\propto{}\varepsilon^{0}\,m^0\,\rho_{\LL}^0\,\rho_{\AA}^0\,\eta_{\LL}^0\,\eta_{\AA}^0\,h_{\MIN}^{-1}\,\tau^{-\frac{1}{2}}
\\
\operatorname{cond}_{\textsc{ch}}&\propto{}\varepsilon^{\frac{1}{2}}\,m^1\,\rho_{\LL}^0\,\rho_{\AA}^0\,\eta_{\LL}^0\,\eta_{\AA}^0\,h_{\MIN}^{-\frac{7}{2}}\,\tau^1
\end{aligned}
\end{equation}
For simultaneous variations of the interface thickness~$\varepsilon$, the minimal element size $h_{\MIN}\propto\varepsilon$, and the mobility $m\propto{}\varepsilon^3$, corresponding to a practical scenario in which the interface-thickness parameter $\varepsilon$ is decreased to approximate the sharp-interface limit more closely, the mobility parameter $m$ is reduced proportional to $\varepsilon^3$ to fix the diffusive time scale, and the minimal element size is decreased proportional to~$\varepsilon$ to ensure adequate resolution of the diffuse interface,
the condition numbers of the equilibrated NSCH, NS and CH (sub-)matrices are proportional to~$\varepsilon^0$, $\varepsilon$ and~$\varepsilon^0$, respectively.

The proposed partitioned solution method is generally very effective in solving the linear tangent problems in Newton's method. GMRES with 
the~$\boldsymbol{P}_{\UPP}$ preconditioner  is more efficient than with the $\boldsymbol{P}_{\LOW}$ preconditioner. 
The~$\boldsymbol{P}_{\LOW}$ preconditioner provides an essentially uniform convergence rate. The~$\boldsymbol{P}_{\UPP}$ preconditioner  exhibits a similar convergence rate in all iterations except the first: in the first iteration, it yields a residual reduction of several orders of magnitude.

We presented a comparison of the numerical results of the AGG NSCH model with a thin diffuse interface with those of a corresponding sharp-interface model for the droplet-oscillation test case. The diffuse-interface results and the sharp-interface results are in excellent agreement, with regard to both the droplet-oscillation dynamics and the energy-dissipation characteristics. At the considered time step, the NSCH model with the CN time discretization was observed to be thermodynamically consistent in the sense that the pertinent free-energy-dissipation relation is satisfied, but this thermodynamic consistency is conditional on the time step and is violated at larger time steps.

The proposed enhancements of the SEMR method and the corresponding solution strategy enable the solution of the AGG NSCH diffuse-interface model for binary-fluid flow problems for physically relevant parameter values, i.e.\ with realistic values of the surface tension and of the liquid and ambient densities and viscosities. By virtue of the adaptive-refinement strategy, it is feasible to consider diffuse-interface thicknesses that are several orders of magnitude smaller than other characteristic length scales in the problem under consideration, thus enabling the investigation of sharp-interface limits. The proposed $\varepsilon$\nobreakdash-continuation strategy is congruent with the adaptive-refinement procedure, because the continuation iterations can be merged with the adaptive-refinement iterations in the SEMR process. The fact that 
$\varepsilon$\nobreakdash-continuation also effectively increases the admissible time step in the time-integration procedure, suggests that the continuation approach could also be beneficial in non-adaptive methods with implicit time-integration schemes. 

In the presented numerical simulations, the $\varepsilon$\nobreakdash-continuation strategy and the partitioned solution method improve the robustness of the solution procedure to the extent that the accuracy and stability (in the sense of thermodynamic consistency) of the considered time-integration scheme becomes a limiting factor. This suggests that
significant gains in efficiency can be achieved by means of high-order implicit time-integration schemes for NSCH systems, in conjunction with the proposed $\varepsilon$\nobreakdash-continuation procedure. Nonetheless, we anticipate that scenarios exhibiting temporal multi-scale behavior, in particular, due to topological changes in the fluid--fluid interface caused by break-up or coalescence, necessitate adaptive refinement of the time step.

\section*{Acknowledgments}
This research was partly conducted within the Industrial Partnership Program {\it Fundamental Fluid Dynamics Challenges in 
Inkjet Printing\/} ({\it FIP\/}), a joint research program of Canon Production Printing, Eindhoven University of Technology,
University of Twente, and the Netherlands Organization for Scientific Research (NWO). T.H.B.\ Demont and C.\ Diddens gratefully acknowledge financial support through the FIP program. All simulations have been performed using the open source software package Nutils \cite{nutils}  (www.nutils.org).%\mycite{nutils}  (www.nutils.org).

\appendix

\section{Characteristic diffusive time scale of the CH equations}
\label{app:time_scale}
We present a concise derivation of the diffusive time scale $T_{\DIFF}$ of the Cahn--Hilliard equations, which characterizes the equilibration rate of the liquid-ambient diffuse interface. We consider the generic (i.e.\ with unit parameters) Cahn--Hilliard equations:
\begin{equation}
\label{eq:CHgeneric}
\begin{aligned}
\partial_t \bar{\varphi} (t,x) & = \Delta \bar{\mu}(t,x),
\\
\bar{\mu}(t,x)                & = \Psi'(\bar{\varphi} (t,x) ) - \Delta \bar{\varphi} (t,x).
\end{aligned}
\end{equation}
We define rescaled $\varphi_{\varepsilon,m,\sigma}$ and $\mu_{\varepsilon,m,\sigma}$ as
\begin{equation}
\label{eq:phimuscaled}
\begin{aligned}
\varphi_{\varepsilon,m,\sigma} (t,x) & \coloneqq \bar{\varphi} (t/T, x/X ),
\\
\mu_{\varepsilon,m,\sigma} (t,x) & \coloneqq M \bar{\mu}(t/T, x/X ),
\end{aligned}
\end{equation}
for a suitable time scale $T\coloneqq{}T_{\varepsilon,m,\sigma}$, length scale $X\coloneqq{}X_{\varepsilon,m,\sigma}$, 
and scaling parameter $M\coloneqq{}M_{\varepsilon,m,\sigma}$. We next determine $T$, $X$, and $M$ such that the scaled functions
in~\EQ{phimuscaled} satisfy the parametrized Cahn--Hilliard equations:
\begin{equation}
\label{eq:CHparam}
\begin{aligned}
\partial_t \varphi_{\varepsilon,m,\sigma} (t,x) & = m \Delta \mu_{\varepsilon,m,\sigma} (t,x), 
\\
\mu_{\varepsilon,m,\sigma} (t,x) & = \sigma \varepsilon^{-1} \Psi' \big(\varphi_{\varepsilon,m,\sigma} (t,x)\big)
- \sigma \varepsilon \Delta \varphi_{\varepsilon,m,\sigma} (t,x).
\end{aligned}
\end{equation}
In order for the rescaled functions in~\EQ{phimuscaled} to satisfy~\EQ{CHparam}, it must hold that
\begin{multline}
\label{eq:Mmu}
M\bar{\mu}(t/T,x/X) 
= \sigma \varepsilon^{-1}\Psi'\big(\bar{\varphi}(t/T,x/X)\big) - \sigma\varepsilon\Delta_x\bar{\varphi} (t/T, x/X)
\\
= \sigma \varepsilon^{-1}\Psi'\big(\bar{\varphi}(t/T,x/X)\big) - \sigma\varepsilon{}X^{-2}\Delta\bar{\varphi} (t/T, x/X)
\end{multline}
where $\Delta\bar{\varphi}$ is to be interpreted as the Laplacian of~$\bar{\varphi}$, acting on its second argument, and $\Delta_x$ is the Laplacian acting on the~$x$ argument.
The second identity in~\EQ{Mmu} follows directly from the chain rule. Moreover, it must hold that
\begin{equation}
\label{eq:T-1phi}
T^{-1}\bar{\varphi}{}'(t/T, x/ X)
=
\partial_t \varphi_{\varepsilon,m,\sigma} (t,x) 
=
m\Delta_x\mu_{\varepsilon,m,\sigma} (t,x)  
=mMX^{-2}\Delta \bar{\mu}(t/T,x/X)
\end{equation}
where $(\cdot)'$ denotes the derive to a function's first argument.
From~\EQ{Mmu} and~\EQ{T-1phi}, we infer that the scaled functions~\EQ{phimuscaled} indeed satisfy the parametrized CH equations~\EQ{CHparam} if
\begin{equation}
\label{eq:MMT}
M=\sigma\varepsilon^{-1},
\qquad
M=\sigma\varepsilon{}X^{-2},
\qquad
T^{-1}=mMX^{-2},
\end{equation}
and hence $M=\sigma\varepsilon^{-1}$, $X=\varepsilon$ and~$T=\varepsilon^3\sigma^{-1}m^{-1}$. We identify the time scale~$T$ as the characteristic time scale of the diffusive process in the CH equations and, accordingly, denote it by $T_{\DIFF}$.\\

\section{Diffuse interface expansion dynamics of the CH equations}
\label{app:interface_expansion}
To elucidate the dynamics of the interface expansion as in the $\varepsilon$-continuation process, as discussed in Remark~\ref{rem:epscont}, we regard a generic interface equilibration test case between two equilibrium solutions corresponding to two different $\varepsilon$ values. We consider the generic Navier--Stokes--Cahn--Hilliard equations analogous to the generic CH equations~\EQ{CHgeneric}:
\begin{equation}
\label{eq:NSCHgeneric}
\left.
\begin{split}
\partial_t \left( \bar{\boldsymbol u} \right) + \nabla \cdot \left( \bar{\boldsymbol u} \otimes \bar{\boldsymbol u} \right) + \nabla \bar{p} - \nabla \cdot \bar{\boldsymbol \tau} - \nabla \cdot \bar{\boldsymbol \zeta} & = 0\\
\nabla \cdot \bar{\boldsymbol u} & = 0\\
\partial_t \bar{\varphi} + \nabla \cdot \left( \bar{\varphi} \bar{\boldsymbol u} \right) - \Delta \bar{\mu} & = 0\\
\bar{\mu} + \Delta \bar{\varphi} - \Psi' & = 0
\end{split}
\right\}.
\end{equation}
The viscous stress tensor is given by
\begin{equation}
\bar{\boldsymbol \tau} \coloneqq 2 \bar{\varepsilon}_T + \Lambda \textrm{tr} ( \bar{\varepsilon}_T ) \boldsymbol I,
\end{equation}
where the symmetric strain-rate tensor $\bar{\varepsilon}_T \coloneqq \frac{1}{2} ( \nabla \bar{\boldsymbol u} + ( \bar{\boldsymbol u} )^T )$ is applied over the velocity field $\bar{\boldsymbol u}$. Furthermore, the capillary-stress tensor is given by:
\begin{equation}
\bar{\boldsymbol \zeta} = - \nabla \bar{\varphi} \otimes \nabla \bar{\varphi} + \boldsymbol I \left( \frac{1}{2} | \nabla \bar{\varphi} |^2 + \Psi \right).
\end{equation}
Note that the corresponding diffusive time scale to this system of equations is equal to $T_{\DIFF} = 1$. We consider a 1D domain $\Omega \coloneqq (-10,10)$, with boundary conditions on $\partial\Omega=\{-10,10\}$ according to~\EQ{EXT}, and initial conditions:
\begin{equation}
\left\{
\begin{aligned}
\bar{\varphi}( 0,x ) &= \bar{\varphi}_0 ( x) \coloneqq - \tanh \Big( \frac{x}{\sqrt{2}\varepsilon_{\ast}} \Big),\\
\bar{u} \left( 0,x \right) & = \bar{u}_0(x) \coloneqq 0,
\end{aligned}
\right.
\end{equation}
with $\varepsilon_{\ast}=2^{-3}$. Let us note that the initial data $\bar{\varphi}_0,\bar{u}_0$ corresponds for an equilibrium solution of the NSCH system~\EQ{strong} on $\Omega=\mathds{R}$ for $\varepsilon_{\ast} = 2^{-3}$, but not for the considered interface thickness $\varepsilon=1$ in~\EQ{NSCHgeneric}. The considered setting therefore essentially pertains to the expansion of the diffuse interface from length scale $\varepsilon_{\ast}=2^{-3}$ to length scale $\epsilon=1$. We assume that the effect of the domain truncation in combination with the boundary conditions~\EQ{EXT} on the interface equilibration process is negligible and, accordingly, that the phase field $\bar{\varphi}$ approaches the equilibrium solution $\bar{\varphi}_\infty ( x) \coloneqq - \tanh ({x}/{\sqrt{2}} )$ as $t \rightarrow \infty$. 

To examine the interface-equilibration process, we consider a finite-element approximation on a uniform mesh with mesh width $h = 10^{-2}$, without adaptive spatial refinement and $\varepsilon$-continuation. The time step size is restricted to $\tau_{\MAX} = 10^{-1}$, with an initial time-step-reduction factor $m=12$, such that $\tau_0 = 2^{-12} \tau_{\MAX}$. The small initial time-step size is selected to resolve the initial fast dynamics of the interface with sufficient accuracy. 

\begin{figure}[h!!!!!!!!!!]
\centering
\includegraphics[width=\textwidth, trim = 0mm 1mm 10mm 9mm, clip]{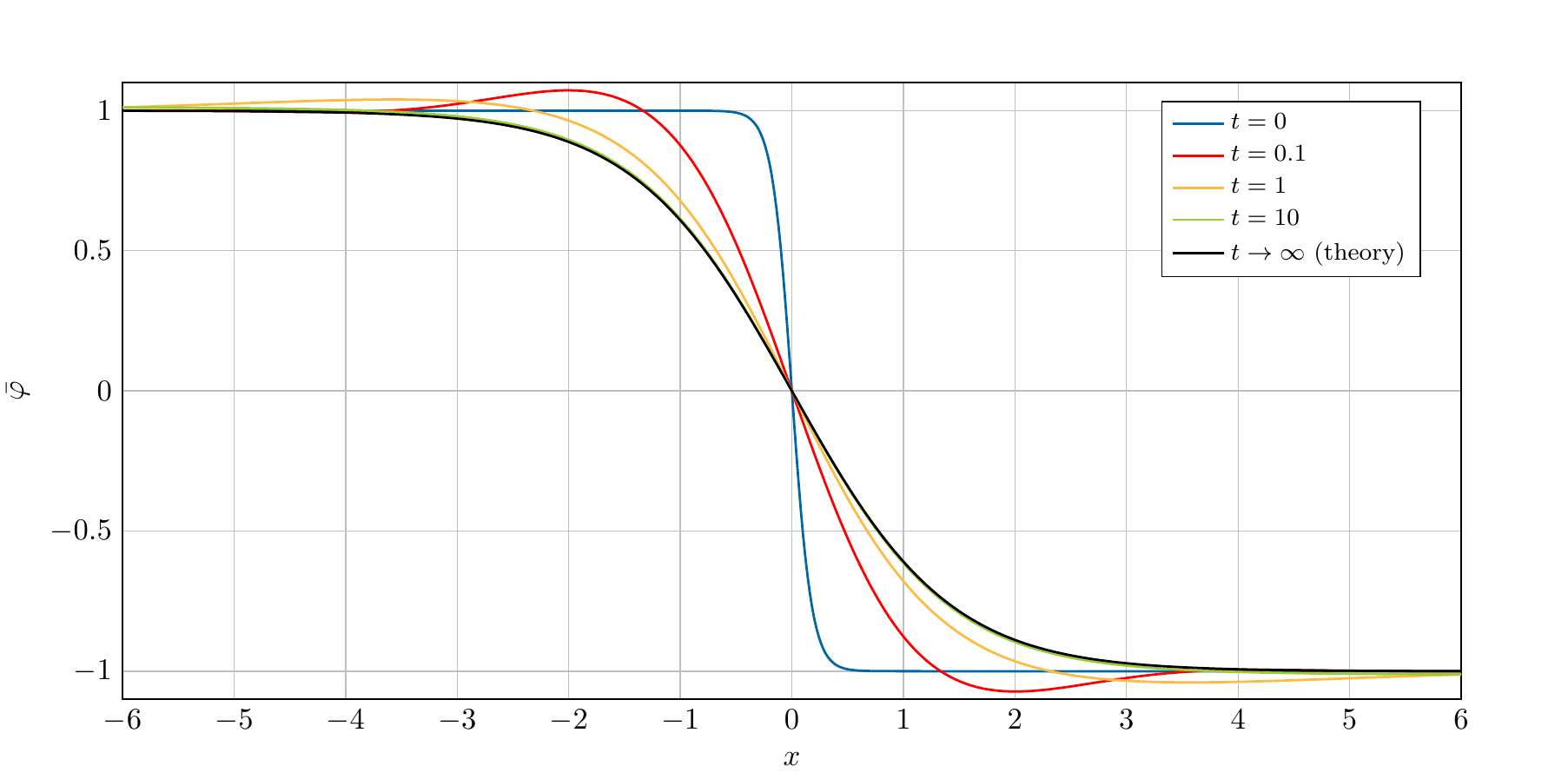}
\caption{Computed phase field $\bar{\varphi}(t,x)$ of \EQ{NSCHgeneric} for $t=10^{-1},1,10$, and initial condition $\bar{\varphi}_0$ and equilibrium profile $\bar{\varphi}_{\infty}(x)$.}
\label{fig:phi_equi}
\end{figure}

The computed phase field $\bar{\varphi}(t,x)$ is presented in Figure \FIG{phi_equi} for $t=10^{\{-1,0,1\}}$. For reference, the equilibrium profile $\bar{\varphi}_\infty ( x) \coloneqq - \tanh ({x}/{\sqrt{2}} )$ is displayed as well. One can observe that $\bar{\varphi}(t,x)$ indeed displays under- and overshoots, exceeding the bounds $\pm{}1$. During the evolution, the velocity field $\bar{u}(t,x)$ and the corresponding kinetic energy $E_{\textrm{kin}}$ remain, up to numerical errors, equal to $0$.
This implies that the interface-expansion evolution is essentially independent of the fluid properties and, hence, that $(\bar{\varphi} , \bar{\mu})$ corresponds to a solution of the generic CH equations \EQ{CHgeneric}. The results in Figure~\FIG{phi_equi} moreover convey that the diffuse interface is reasonably well equilibrated for $t\geq{}1$, and indistinguishable from the equilibrium profile for $t\geq{}10$.

\bibliographystyle{plain}
\bibliography{bibfile}\label{sec:bibliography}

\end{document}